\documentclass[10pt,nocut]{article}

\usepackage{./packages}
\usepackage{./notations}

\usepackage{tikz,tkz-tab}
\usepackage{float}

\usepackage{titlesec}
\makeatletter
\newcommand{\setappendix}{Appendix~\thesection:~~}
\newcommand{\setsection}{\thesection~~}
\titleformat{\section}{\bfseries\LARGE}{%
	\ifnum\pdfstrcmp{\@currenvir}{appendices}=0
		\setappendix
	\else
		\setsection
\fi}{0em}{}
\makeatother
\usepackage[titletoc]{appendix}

\makeatletter
\def\blfootnote{\xdef\@thefnmark{}\@footnotetext}
\makeatother

\title{The distribution of the Lasso:\\
Uniform control over sparse  balls and adaptive parameter tuning}

\author{L\'eo Miolane$^{\star}$ and Andrea Montanari$^{\dagger}$}
\date{\today}

\begin{document}
\blfootnote{%
	\hspace{-5.2mm}${\star}$ 
   D\'epartement d'Informatique de l'ENS, \'Ecole Normale Sup\'erieure, CNRS, PSL Research University \& Inria, Paris, France.
   \\
   ${\dagger}$ Department of Electrical Engineering and Department of Statistics, Stanford University.
}
\maketitle

\begin{abstract}
	The Lasso is a popular regression method for high-dimensional problems in which the number of parameters $\theta_1,\dots,\theta_N$,
	is larger than the number $n$ of samples: $N>n$. 
	A useful heuristics relates the statistical properties of the Lasso estimator to that of a simple  soft-thresholding denoiser,
	in a denoising problem in which the parameters  $(\theta_i)_{i\le N}$ are observed in Gaussian noise, with a carefully tuned variance.
	Earlier work confirmed this picture in the limit $n,N\to\infty$, pointwise in the parameters $\theta$, and in the value of the regularization parameter.

	Here, we consider a standard random design model and prove exponential concentration of
	its empirical distribution around the prediction provided by the Gaussian denoising model. Crucially, our results are uniform with respect to $\theta$
	belonging to $\ell_q$ balls, $q\in [0,1]$, and with respect to the regularization parameter. This allows to derive sharp results for the performances of various
	data-driven procedures to tune the regularization.

	Our proofs make use of Gaussian comparison inequalities, and in particular of a version of Gordon's minimax theorem developed by
	Thrampoulidis, Oymak, and Hassibi, which controls the optimum value of the Lasso optimization problem.
	Crucially, we prove a stability property of the minimizer in Wasserstein distance, that allows to characterize properties of the minimizer itself.
\end{abstract}


\section{Introduction}

Given data $(x_i,y_i)$, $1\le i\le n$, with $x_i\in\R^N$, $y_i\in \R$, the Lasso~\cite{tibshirani1996regression,chen1995examples} fits a linear model by minimizing the cost function
\begin{align}
	\mathcal{L}_{\lambda}(\theta)&= \frac{1}{2n}\sum_{i=1}^{n}\big(y_i-\<x_i, \theta\>\big)^2 + \frac{\lambda}{n} |\theta| \nonumber\\
								 & = \frac{1}{2n}\left\| y-X \theta  \right\|^2 + \frac{\lambda}{n} |\theta| \, . \label{eq:LassoCost}
\end{align}
Here $X\in \R^{n\times N}$ is the matrix with rows $x_1,\dots,x_n$, $y = (y_1,\dots,y_n)$, $\|v\|$ denotes the $\ell_2$ norm of vector $v$,
and $|v|$ its $\ell_1$ norm. To fix normalizations, we will assume
that the columns of $X$ have $\ell_2$ norm $1+o(1)$. (Note that this normalization is different from the one that is sometimes adopted in the literature, but the two are completely
equivalent.)

A large body of theoretical work supports the use of $\ell_1$ regularization in the high-dimensional regime $n\lesssim N$,
when only a small subset of the coefficients $\theta$ are expected to be large. Broadly speaking, we can distinguish two types of 
theoretical approaches. A first line of work makes deterministic assumptions about the design matrix $X$, such as the restricted isometry
property and its generalizations~\cite{candes2005decoding,buhlmann2011statistics}. Under such conditions, minimax optimal estimation rates  as well as oracle inequalities 
have been proved in a remarkable sequence of papers~\cite{candes2007dantzig,bickel2009simultaneous,van2009conditions,negahban2012unified,raskutti2011minimax}. 
As an example, assume that that the linear model is correct. Namely, 
\begin{equation}
	y =  X \theta^{\star} +  \sigma z \,,\label{eq:LinearModel}
\end{equation}
for $\sigma \geq 0$, $z\sim\cN(0,\id_n)$, and $\theta^\star$ a vector with $s_0$
non-zero entries. Then, a theorem of Bickel, Ritov and Tsybakov~\cite{bickel2009simultaneous} implies that,
with high probability,
\begin{align}
	\lambda \ge \sigma \sqrt{c_0 \log N}\;\;\; \Rightarrow \;\;\;\| \what{\theta}_{\lambda} - \theta^{\star} \|^2 \le C s_0 \lambda^2\, ,\label{eq:RIP_Result}
\end{align}
for some constants $c_0$, $C$ that depend on the specific assumptions on the design. (The normalization of~\cite{bickel2009simultaneous} is recovered by setting
$\sigma^2=\sigma^2_{\#}/n$, where $\sigma_{\#}^2$ is the noise variance of~\cite{bickel2009simultaneous}.)

\begin{figure*}[h!]
	\centering
	\hspace*{-1.3cm}
	\includegraphics[width=1.1\linewidth]{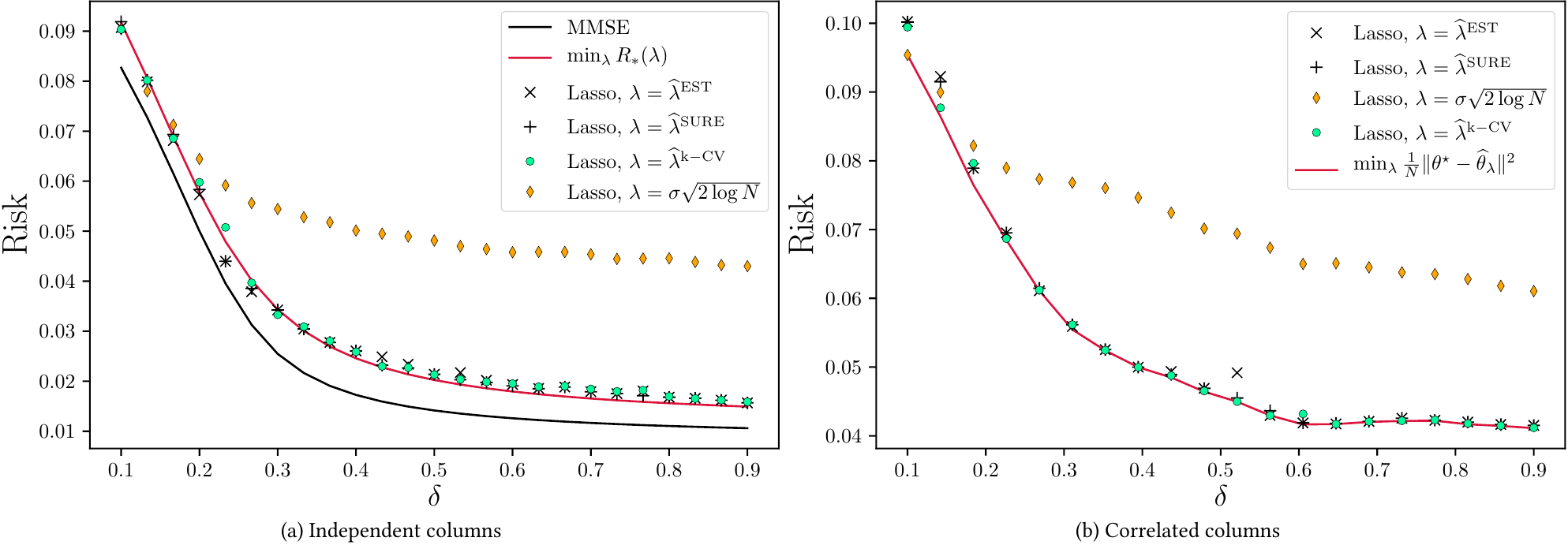}
	\caption{Estimation risk of the Lasso for different choices of $\lambda$, as a function of $\delta$. $N=8000$.
		In both plots, $\sigma = 0.2$. The true coefficients vector $\theta^{\star}$ is chosen to be $sN$-sparse with $s=0.1$. The entries on the support of $\theta^{\star}$ are drawn i.i.d.\ $\cN(0,1)$. 
	 Cross-validation  is carried out using $4$ folds. SURE  is computed using the estimator $\what{\sigma}$ for the plot on the left, and the true value of $\sigma$ on the right.
	 \\
	\textit{Left:} A standard random design with $(X_{ij})\sim_{iid}\cN(0,1/n)$.
	\\
	\textit{Right:} The rows of the design matrix $X$ are i.i.d. Gaussian, with correlation structure given by an autoregressive process, see Eq.~(\ref{eq:Auto}). 
Here we used $\phi = 2$.
}
	\label{fig:comp_lambda_delta}
\end{figure*}
Unfortunately, this analysis provides limited insight into the choice of the regularization parameter $\lambda$
which --in practice-- can impact significantly the estimation accuracy. As an example, Fig.~\ref{eq:RIP_Result} reports the result of a small simulation in which 
we compare four different methods of selecting $\lambda$. The bound of Eq.~(\ref{eq:RIP_Result}) suggests to set $\lambda = \sigma\sqrt{c_0\log N}$.
For the standard random design used in the left frame, the optimal
constant is expected to be $c_0=2$~\cite{donoho1994minimax,donoho2011noise}. 
We compare this method to three procedures that adapt the choice of $\lambda$ to the data: cross validation (CV), Stein's Unbiased Risk Estimate (SURE),
and a procedure that minimizes an estimate of the risk (EST). We refer to the next sections for further details on these methods. 
Note that all of these adaptive procedures significantly outperform the 
`theory driven' $\lambda$: over a broad range of sample sizes $n$, the resulting estimation error is $2$ to $3$ times smaller.
Further,  the error achieved by these methods is quite close to the Bayes optimum.

These empirical observations are not captured by the bound (\ref{eq:RIP_Result}), or by similar results.

An alternative style of analysis postulates an idealized model for the data and derives asymptotically exact results. 
Throughout this paper we will consider the simplest of such models, by assuming that design matrix
to have i.i.d. entries $X_{ij}\sim\cN(0,1/n)$. While this assumption is likely to be violated in practice, it allows
to derive useful insights that are mathematically consistent, and susceptible of being generalized to a broader context.
This type of analysis was first carried out in the context of the Lasso in~\cite{bayati2012lasso} and then extended to a number of other 
problems, see e.g.~\cite{karoui2013asymptotic,thrampoulidis2015regularized,donoho2016high,thrampoulidis2018precise,el2018impact,sur2018modern}. 
As an example, Figure~\ref{fig:comp_lambda_delta} reports the predictions of this analysis
for the risk of the three adaptive procedure for selecting $\lambda$. The agreement with the numerical simulations is excellent. 

Unfortunately, the results in~\cite{bayati2012lasso} (and in follow-up work) do not allow to derive in a mathematically rigorous way curves such as the ones in
Figure~\ref{fig:comp_lambda_delta}. In fact earlier results hold `pointwise' over $\lambda$ and hence do not apply to adaptive procedures to select $\lambda$.
Further they provide asymptotic estimates  `pointwise' over $\theta$, and hence do not allow to compute --for instance-- minimax risk.

In order to clarify these points, it is useful to overview informally the picture emerging from 
\cite{bayati2012lasso,donoho2009message}. Fix $\theta\in\reals^N$, $\lambda\in\reals_{>0}$, 
and let $\eta(x;b) = (|x|-b)_+\sign(x)$ be the soft thresholding function.
By the KKT conditions the Lasso estimator $\htheta_{\lambda}$ satisfies
\begin{align}
	\htheta_{\lambda}= \eta\big(\detheta_{\lambda};\alpha\tau\big)\, ,\;\;\;\;\;\;\;\;
	\detheta_{\lambda} = \htheta_{\lambda} + \frac{\alpha\tau}{\lambda} X^{\sT}(y - X \what{\theta}_{\lambda})\, ,
\end{align}
where the vector $\detheta_{\lambda}$ is also referred to as the `debiased Lasso'~\cite{zhang2014confidence,van2014asymptotically,javanmard2014confidence}.
The above identity holds for arbitrary $\alpha,\tau>0$. 
However,~\cite{bayati2012lasso} predicts that the distribution of the debiased estimator $\detheta_{\lambda}$
simplifies dramatically for specific choices of these parameters.

Namely, let $\Theta$ be a random variable with distribution given by the empirical distribution of $(\theta_i)_{i\le N}$
(i.e., $\Theta=\theta_i$ with probability $1/N$, for $i\in\{1,\dots,N\}$) and let $Z\sim\cN(0,1)$ be independent of $\Theta$.
Define $\alpha_*,\tau_*$ to be the solution of the following system of equations (we refer to Section~\ref{sec:Definitions} for a discussion of existence and uniqueness):
\begin{equation}
	\begin{cases}
		\tau^2 &\!\!\!\!= \, \sigma^2 + \frac{1}{\delta} \E \Big[(\eta(\Theta + \tau Z, \alpha\tau)- \Theta)^2 \Big]\, , \\
		\lambda &\!\!\!\!= \, \alpha \tau \left(1 - \frac{1}{\delta} \P\left(\big|\Theta + \tau Z\big|>\alpha\tau  \right) \right).
	\end{cases}
\end{equation}
When $\alpha,\tau$ are selected in this way, $\detheta_{\lambda}$ is approximately normal with mean $\theta^\star$ (the true parameters vector) and 
variance $\tau_*^2$: $\detheta \approx \cN(\theta^{\star},\tau_*^2\id)$. More precisely, for any  test function ${\sf f}:\reals\times\reals\to\reals$,
with $|{\sf f}(x)-{\sf f}(y)|\le L(1+\|x\|+\|y\|)\|x-y\|$,
almost surely,
\begin{align}
	\lim_{N\to\infty} \frac{1}{N}\sum_{i=1}^N{\sf f}(\theta^\star_i,\detheta_{\lambda,i}) = \E\big\{{\sf f}(\Theta,\Theta+\tau_* Z)\big\}\, ,\;\;\;
	\lim_{N\to\infty} \frac{1}{N}\sum_{i=1}^N{\sf f}(\theta^\star_i,\htheta_{\lambda,i}) = \E\big\{{\sf f}(\Theta,\eta(\Theta+\tau_* Z;\alpha_*\tau_*))\big\}\, .\label{eq:Earlier}
\end{align}
This is an asymptotic result, which holds along sequences of problems with: 
$(i)$ Converging aspect ratio $n/N\to\delta\in (0,\infty)$; $(ii)$ Fixed regularization $\lambda\in (0,\infty)$; $(iii)$ Parameter vectors  $\theta^\star=\theta^\star(n)$ 
whose empirical distribution converges (weakly) to a limit law $p_{\Theta}$. 
As emphasized above, this does not allow deduce the behavior of the Lasso with adaptive choices of $\lambda$ (there could be deviations from the above limits for exceptional 
values of $\lambda$), or to compute the minimax risk (there could be deviations for exceptional vectors $\theta^\star$).

The importance of establishing uniform convergence with respect to the regularization parameter $\lambda$ 
was recently emphasized by Mousavi, Maleki, and Baraniuk~\cite{mousavi2017consistent}. Among other results, these authors derive a uniform convergence 
statement  for the related approximate message passing (AMP) algorithm. However, in order to establish uniform convergence,  they have to construct an
ad-hoc smoothing of the quantity of interest, which is roughly equivalent to discretizing the corresponding tuning parameter.

In this paper, we obtain uniform (in $\lambda$) convergence results for the Lasso, hence providing a sound mathematical basis to the comparison of various adaptive procedures,
as well as to the study of minimax risk. 

The rest of the paper is organized as follows. Section~\ref{sec:Related} reviews related work. We state our main theoretical results in Section~\ref{sec:Results}.
In Section~\ref{sec:Applications} we apply these results to two types of statistical questions: estimating the risk and noise level, and selecting $\lambda$ through adaptive 
procedures. Further, we illustrate our results in numerical simulations. Finally, Section~\ref{sec:Proofs} outlines the main proof ideas, with most technical 
legwork deferred to the appendices.

\section{Related work}
\label{sec:Related}

There is --by now-- a substantial literature on determining exact asymptotics in high-dimensional statistical models,
and a number of mathematical techniques have been developed for this
task. We will only provide a few pointers focusing on
high-dimensional regression problems.

The original proof of~\cite{bayati2012lasso} was based on an asymptotically exact analysis of an approximate message passing (AMP)
algorithm~\cite{bayati2011amp} that was first proposed in~\cite{donoho2009message} to minimize the Lasso cost function. 
Variants of AMP have been developed in a number of contexts, opening the way to the analysis of various statistical estimation problems.
A short list includes generalized linear models~\cite{rangan2011generalized}, phase retrieval~\cite{schniter2015compressive,ma2018optimization}, 
robust regression~\cite{donoho2016high}, logistic regression~\cite{sur2018modern}, generalized compressed sensing~\cite{berthier2017state}.
This approach is technically less direct than others, but has the advantage of providing an efficient algorithm, and is
and not necessarily limited to convex problems (see~\cite{montanari2016non} for a non-convex example). 

As mentioned above, our work was partially motivated by the recent results of Mousavi, Maleki, and Baraniuk~\cite{mousavi2017consistent}
that establish a form of uniformity for the AMP estimates --but not for the Lasso solution.  It would be interesting
to understand whether the approach of~\cite{mousavi2017consistent} could also be used to obtain uniform results for the Lasso
or other statistical estimators.

Here we follow a different route that exploits powerful Gaussian comparison inequalities first proved
by Gordon~\cite{gordon1985some,gordon1988milman}. Gordon inequality allows to bound the distribution of a minimax 
value, i.e. the value of a random variable $G_*= \min_{i\le N}\max_{j\le M}G_{ij}$, where $(G_{ij})_{i\le N,j\le M}$ is a Gaussian process,
in terms of a similar quantity for a `simpler' Gaussian process.
The use of Gordon's inequality in this context was pioneered by Stojnic~\cite{stojnic2013framework} and
then developed by a number of authors in the context of regularized regression~\cite{thrampoulidis2015regularized},
M-estimation ~\cite{thrampoulidis2018precise}, generalized compressed sensing~\cite{amelunxen2014living}, binary
compressed sensing~\cite{stojnic2010recovery} and so on. The key idea is to write the optimization problem of interest as a minimax
problem, and then apply a suitable version of Gordon's inequality. A matching bound is obtained by convex duality and then 
a second application of Gordon's inequality. In particular, convexity of the cost function of interest is a crucial ingredient.

While the Gaussian comparison inequality provides direct access to the value of the optimization problem, 
understanding the properties of the estimator can be more challenging. In this paper we identify a property (that we call \emph{local stability})
that allows to transfer information on the minimum (the Lasso cost)  into information about the minimizer (the Lasso estimator).
We believe this strategy can be applied to other examples beyond the Lasso.

Independently, a different approach based on leave-one-out techniques was developed by El Karoui
in the context of ridge-regularized robust regression~\cite{karoui2013asymptotic,el2018impact}.

Finally, a parallel line of research determines exact asymptotics for Bayes optimal estimation,
under a model in which the coordinates of $\theta$  are i.i.d. with common distribution $p_{\Theta}$. 
In particular, the asymptotic Bayes optimal error for linear regression with random designs was recently determined in
\cite{barbier2016CS,reeves2016replica}. Of course --in general-- Bayes optimal estimation requires knowledge of the 
distribution $p_{\Theta}$, and is not computationally efficient. We will use this Bayes-optimal error as a benchmark of our adaptive 
procedures. Generalizations of these results were also obtained in~\cite{barbier2018optimal} for other regression problems. 
A successful approach to these models uses smart interpolation techniques that generalize ideas in spin-glass theory.

\section{Main results}
\label{sec:Results}

\subsection{Definitions}
\label{sec:Definitions}

As stated above, we consider the standard linear model (\ref{eq:LinearModel}) where
$y = X \theta^{\star} + \sigma z$
, with noise $z\sim\cN(0,\id_n)$, and $X$ a Gaussian design:
$(X_{i,j})_{i\le n,j\le N} \iid \cN(0,1/n)$. The Lasso estimator is defined by
\begin{equation}\label{eq:def_lasso_est}
	\htheta_{\lambda} = \arg\min_{\theta\in\reals^N}\mathcal{L}_{\lambda}(\theta)\,.
\end{equation}
(The minimizer is almost surely  unique since the columns of $X$ are in generic positions.)
We  set $\delta =n/N$ to be the number of samples per dimension.
We are interested in uniform estimation over sparse vectors $\theta^\star$. Following~\cite{donoho1994minimax,johnstone2002function} we formalize this notion using $\ell_p$-balls
(which are convex sets only for $p\ge 1$).
\begin{definition}
	Define for $p,\xi>0$ the $\ell_p$-ball
	$$
	\cF_p(\xi) = \left\{ x \in \R^N \,
		\middle| \, \frac{1}{N} \sum_{i=1}^N |x_i|^p \leq \xi^p
	\right\} \, ,
	$$
	and for $s \in [0,1]$
	$$
	\cF_0(s) = \left\{ x \in \R^N \, \middle| \, \|x\|_0 \leq s N \right\} .
	$$
\end{definition}
By Jensen's inequality we have for $p\geq p'>0$, $\cF_{p}(\xi) \subset \cF_{p'}(\xi)$.

Let $\phi(x) = \frac{e^{-x^2/2}}{\sqrt{2\pi}}$ be the standard Gaussian density and $\Phi(x) = \int_{-\infty}^x \phi(t) dt$ be the associated cumulative function.
In the case of $\ell_0$ balls (sparse vectors), a crucial role is played by the following sparsity level.
\begin{definition}
	Define the critical sparsity as
	$$
	s_{\rm max}(\delta) = \delta \max_{\alpha \geq 0} \left\{
		\frac{1 - \frac{2}{\delta}\big((1+\alpha^2) \Phi(-\alpha) - \alpha \phi(\alpha)\big)}{1 + \alpha^2 - 2 \big((1+\alpha^2) \Phi(-\alpha) - \alpha \phi(\alpha)\big)} 
	\right\}.
	$$
\end{definition}
The critical sparsity curve first appears in the seminal work by Donoho and Tanner on compressed sensing ~\cite{donoho2005neighborliness,donoho2006high}.
These authors consider the noiseless case ($z=0$) of model (\ref{eq:LinearModel}) and reconstruction via $\ell_1$ minimization (which corresponds to the $\lambda\to 0$
limit of the Lasso). They prove that  $\ell_1$ minimization reconstructs exactly $\theta^\star$ with high probability, if 
$\|\theta^{\star}\|_0\le N (s_{\max}(\delta)-\eps)$, and fails with high probability if $\|\theta^{\star}\|_0\ge N(s_{\max}(\delta)+\eps)$ (for any $\eps>0$).
A second interpretation of the critical sparsity $s_{\rm max}(\delta)$ was given in~\cite{donoho2011noise,tropp2015convex,thrampoulidis2015regularized}. For 
$\|\theta^{\star}\|_0\le N (s_{\max}(\delta)-\eps)$, the Lasso achieves stable reconstruction. Namely, there exists $M = M(s,\delta)<\infty$ for
$s<s_{\max}(\delta)$, such that, if $\|\theta^\star\|_0\le Ns$,  $\|\htheta_{\lambda}-\theta^\star\|_2\le M(s,\delta)\sigma^2$.  
Our results provide a third interpretation: uniform limit laws for the Lasso will be obtained on $\ell_0$ balls only for $s<s_{\max}(\delta)$. 

A crucial role in our results is provided by the following max-min problem:
\begin{align} \label{eq:max_min_scalar}
	&\max_{\beta \geq 0} 
	\min_{\tau \geq \sigma}\;\;\;\;\psi_{\lambda}(\beta,\tau) \, ,\\
	&\;\psi_{\lambda}(\beta,\tau)  \equiv \left(\frac{\sigma^2}{\tau} + \tau\right) \frac{\beta}{2}
	- \frac{1}{2} \beta^2
	+ 
	\frac{1}{\delta}
	\E
	\min_{w \in \R}\left\{
		\frac{w^2}{2\tau} \beta
		-	\beta Z w
	+ \lambda |w+\Theta| - \lambda |\Theta| \right\}.\nonumber
\end{align}
The expectation above is with respect to $(\Theta,Z) \sim \what{\mu}_{\theta^{\star}} \otimes \cN(0,1)$, where $\what{\mu}_{\theta^{\star}}$ denotes the empirical distribution of the entries of the vector $\theta^{\star}$:
$$
\what{\mu}_{\theta^{\star}} = \frac{1}{N} \sum_{i=1}^N \delta_{\theta^{\star}_i} \,.
$$

\begin{proposition}\label{th:scalar_max_min}
	The max-min~\eqref{eq:max_min_scalar} is achieved at a unique couple $(\beta_*(\lambda),\tau_*(\lambda))$. Moreover, $(\tau_*(\lambda),\beta_*(\lambda))$ is also the unique couple $(\beta,\tau) \in (0,+\infty)^2$ that verify
	\begin{equation}\label{eq:state_evolution}
		\begin{cases}
			\tau^2 &\!\!\!\!= \, \sigma^2 + \frac{1}{\delta} \E \Big[(\eta(\Theta + \tau Z, \tau \frac{\lambda}{\beta})- \Theta)^2 \Big] \\
			\beta &\!\!\!\!= \, \tau \left(1 - \frac{1}{\delta} \E \left[ \eta'(\Theta + \tau Z,\frac{\tau \lambda}{\beta}) \right] \right).
		\end{cases}
	\end{equation}
	We will also use the notation $\alpha_*(\lambda) = \lambda / \beta_*(\lambda)$ and 
	\begin{equation}\label{eq:def_s_star}
	s_*(\lambda) = \E \big[ \eta'(\Theta + \tau_*(\lambda) Z,\tau_*(\lambda) \alpha_*(\lambda)) \big] = \P\big(|\Theta + \tau_*(\lambda) Z| \geq \alpha_*(\lambda) \tau_*(\lambda)\big) \,.
	\end{equation}
\end{proposition}

We will sometimes omit the dependency on $\lambda$ and write simply $\alpha_*, \beta_*, \tau_*,s_*$. The distribution $\mu^*_{\lambda}$ defined below will correspond (see Theorem~\ref{th:unif_lambda_law} in the next section) to the limit of the empirical distribution of the entries of $(\what{\theta}_{\lambda},\theta^{\star})$.
\begin{definition}\label{def:mu_star}
	We denote by $\mu_{\lambda}^*$ the law of the couple
	$
\left( \eta\big(\Theta + \tau_*(\lambda) Z , \, \alpha_*(\lambda) \tau_*(\lambda)\big) , \ \Theta \right)
	$,
	where $(\Theta,Z) \sim \what{\mu}_{\theta^{\star}} \otimes \cN(0,1)$. 
\end{definition}

\subsection{Results}

We fix from now on $0 < \lambda_{\rm min} \leq \lambda_{\rm max}$ and $\mathcal{D} \subset \R^N$ that can be either $\cF_p(\xi)$ for some $\xi,p>0$, or
$\cF_0(s)$ for some $s < s_{\rm max}(\delta)$.
Our uniformity domain is defined by $\Omega = \big(\delta,\sigma, \mathcal{D},\lambda_{\rm min},\lambda_{\rm max}\big)$.
Namely, we will control $\what{\theta}_{\lambda}$ uniformly with respect to $\theta^{\star} \in \mathcal{D}$ and $\lambda \in [\lambda_{\rm min},\lambda_{\rm max}]$,
with $n/N=\delta$. 
We will call \textit{constant} any quantity that only depends on $\Omega$. In absence of further specifications, $C,c$ will be constants (that depend only on $\Omega$) that are allowed to change from one line to another.

Our first result shows that the empirical distribution of the entries $\{(\what{\theta}_{\lambda,i},\theta_i^{\star})\}_{i\le N}$ is uniformly close to the model
$\mu^*_{\lambda}$. We quantify deviations using the Wasserstein distance. Recall that, given two probability measures $\mu,\nu$ on $\reals^d$
with finite second moment,  their Wasserstein distance of order $2$ is 
\begin{align}
	W_2(\mu,\nu) & = \Big(\inf_{\gamma\in \cC(\mu,\nu)}  \int \|x-y\|_2 ^2 \; \gamma(\de x,\de y) \,\Big)^{1/2} , \label{eq:Wasserstein_def}
\end{align}
where the infimum is taken over all couplings of $\mu$ and $\nu$. Note that $W_2$ metrizes the convergence in Eq.~(\ref{eq:Earlier}).
Namely $\lim_{n\to\infty}W_2(\mu_n,\mu_*) = 0$ if and only if, for any test function ${\sf f}:\reals\times\reals\to\reals$,
with $|{\sf f}(x)-{\sf f}(y)|\le L(1+\|x\|+\|y\|)\|x-y\|$, we have
$\lim_{n\to\infty} \int {\sf f}(x) \mu_n(\de x) = \int {\sf f}(x) \mu_*(\de x)$~\cite{villani2008optimal}. It provides therefore a natural way to
extend earlier results to a non-asymptotic regime.
\begin{theorem}\label{th:unif_lambda_law}
	Assume that $\mathcal{D} = \cF_{p}(\xi)$ for some $\xi>0$ and $p>0$.
	Then there exists constants $C,c >0$ that only depend on $\Omega$, such that
	for all $\epsilon \in (0,\frac{1}{2}]$
	$$
	\sup_{\theta^{\star} \in \mathcal{D}} \ \P \left(
		\sup_{\lambda \in [\lambda_{\rm min},\lambda_{\rm max}]} 
		W_2\big(\what{\mu}_{(\what{\theta}_{\lambda},\theta^{\star})},\mu^*_{\lambda}\big)^2 \geq \epsilon
	\right) 
	\leq 
	C \epsilon^{-\max(1,a)-1} N^{(1/p - 1)_+} \exp\left(-cN \epsilon^2 \epsilon^a \log(\epsilon)^{-2} \right) \,,
	$$
	where $a = \frac{1}{2} + \frac{1}{p}$.
\end{theorem}
Theorem~\ref{th:unif_lambda_law} is proved in Section~\ref{sec:proof_unif}.

\begin{remark}
	It is worth emphasizing in what sense Theorem~\ref{th:unif_lambda_law} is uniform with respect to $\lambda\in [\lambda_{\min},\lambda_{\max}]$ and to 
	$\theta^{\star} \in \mathcal{D}$:
	\begin{itemize}
		\item \emph{Uniformity with respect to $\lambda$.} We bound (in probability) the maximum (over $\lambda$) deviation between the empirical distribution 
			$\what{\mu}_{(\what{\theta}_{\lambda},\theta^{\star})}$ and the predicted distribution $\mu^*_{\lambda}$.   (The supremum over $\lambda$  is `inside' the probability.)
		\item \emph{Uniformity with respect to $\theta^\star$.} We bound the maximum probability (over $\theta^\star$) of a deviation between
			$\what{\mu}_{(\what{\theta}_{\lambda},\theta^{\star})}$ and $\mu^*_{\lambda}$.   (The supremum over $\theta^\star$  is `outside' the probability.)
	\end{itemize} 
	The reader might wonder whether it is possible to strengthen this result and bound the maximum deviation over $\theta^{\star}$ (`move the supremum over
	$\theta^\star$ inside'). The answer is negative. In particular, we can choose the support of $\theta^{\star}$ to coincide with a submatrix of $X$ with atypically small minimum 
	singular value. This will result in larger estimation error $\|\htheta_{\lambda}-\theta^{\star}\|_2$, and hence in a large Wasserstein distance
	$W_2(\what{\mu}_{(\what{\theta}_{\lambda},\theta^{\star})},\mu^*_{\lambda})$. 
\end{remark}

\begin{remark}\label{rmk:ell0ball}
	Note that Theorem~\ref{th:unif_lambda_law}  does not hold for $\ell_0$ balls. This is probably a fundamental problem, since controlling $W_2$ distance uniformly over
	$\ell_0$ balls is impossible even in the simple sequence model (or, equivalently, for orthogonal designs $X$). 
	Namely, consider the case in which we observe $y_i = \theta^\star_i+z_i$, $i\le N$, where $(z_i)_{i\le N}\iid \cN(0,\tau_*^2)$, and we try to estimate
	$\theta^\star$ by computing $\htheta_{\lambda,i} = \eta(y_i;\lambda)$. Then there are vectors $\theta^\star\in\cF_0(s)$ such that
	the empirical law $\what{\mu}_{(\htheta_{\lambda},\theta^{\star})}$ does not concentrate in Wasserstein distance  around its expectation $\mu^*_{\lambda}$, i.e. the law of 
	$(\Theta,\eta(\Theta+Z;\lambda)$) for $G\sim\cN(0,\tau_*)$.

	In order to see this, it is sufficient to consider the vector 
	$$
	\theta^\star = (N, 2N, \dots, k N,0,\dots,0)\, .
	$$
	In Appendix~\ref{app:ell0ball_remark}, we prove  that (for this choice of $\theta^\star$)
	there exists a constant $c_0$ such that $W_2(\what{\mu}_{(\htheta_{\lambda},\theta^{\star})},\mu^*_{\lambda})\ge\sqrt{k/N}$ 
	with probability at least $1-e^{-c_0 k}$ for all $N$ large enough.
\end{remark}

We can think of several possibilities to overcome this intrinsic non-uniformity over $\ell_0$ balls. One option 
would be to consider a weaker notion of distance between probability measures. Here we follow a different route, and prove uniform
estimates over $\ell_0$ balls  for several specific quantities of interest. In order to state these results, we introduce the following quantities,
which correspond to the risk and the prediction error (and are expressed in terms of the solution $(\tau_*,\beta_*)$ of (\ref{eq:state_evolution}))
\begin{align}\label{eq:def_R_star}
	R_*(\lambda) &= \delta \big(\tau_*(\lambda)^2 - \sigma^2\big) \, ,\\
	\label{eq:def_P_star}
	P_*(\lambda)& = \beta_*(\lambda)^2 + \frac{2 \sigma^2}{\delta} s_*(\lambda) - \frac{\sigma^2}{\delta} \,.
\end{align} 
\begin{theorem}\label{th:unif_lambda_risk}
	Assume here that $\mathcal{D}$ is either $\cF_0(s)$ or $\cF_p(\xi)$ for some $0 \leq s < s_{\rm max}(\delta)$ and $\xi >0, p >0$.
	There exists constants $C,c >0$ that only depend on $\Omega$, such that
	for all $\epsilon \in (0,1]$
	\begin{align}
		\label{eq:unif_risk}
		\sup_{\theta^{\star} \in \mathcal{D}} \ \P \left(
			\sup_{\lambda \in [\lambda_{\rm min},\lambda_{\rm max}]} 
			\Big( \frac{1}{N}\| \what{\theta}_{\lambda} - \theta^{\star} \|^2 - R_*(\lambda) \Big)^2 \geq \epsilon
	\right) &
	\leq \frac{C}{\epsilon^2} N^q e^{-cN\epsilon^2} \,,\\
	\label{eq:unif_beta}
		\sup_{\theta^{\star} \in \mathcal{D} } \P \Big( \sup_{\lambda \in [\lambda_{\rm min},\lambda_{\rm max}]}\Big( \frac{1}{n}\| y - X\what{\theta}_{\lambda}\|^2 - \beta_*(\lambda)^2 \Big)^2 \geq \epsilon \Big) &\leq \frac{C}{\epsilon^2} N^q e^{-cN\epsilon^2} \,,\\
		\label{eq:unif_prediction}
		\sup_{\theta^{\star} \in \mathcal{D} } \P \left( \sup_{\lambda \in [\lambda_{\rm min},\lambda_{\rm max}]}\Big( \frac{1}{n}\|X(\theta^{\star} - \what{\theta}_{\lambda})\|^2 - P_*(\lambda) \Big)^2 \geq \epsilon \right) &\leq \frac{C}{\epsilon^2} N^q e^{-cN\epsilon^2} \,,
	\end{align}
	where $q = 0$ if\, $\mathcal{D} = \cF_0(s)$ and $q=(1/p - 1)_+$ if\, $\mathcal{D} = \cF_p(\xi)$.
\end{theorem}

The statement~\eqref{eq:unif_risk} is proved in Appendix~\ref{sec:proof_unif}, while~\eqref{eq:unif_beta}-\eqref{eq:unif_prediction} are proved in Appendix~\ref{sec:gordon_u}.

So far we focused on the Lasso estimator $\what{\theta}_{\lambda}$. 
The \emph{debiased Lasso} estimator is defined as
$$
\what{\theta}_{\lambda}^{d}= \what{\theta}_{\lambda} + \frac{X^{\sT}(y - X \what{\theta}_{\lambda})}{1 - \frac{1}{n}\|\what{\theta}_{\lambda}\|_0}
\,.
$$
This estimator plays a crucial role in the construction of confidence intervals and $p$-values~\cite{zhang2014confidence,van2014asymptotically,javanmard2014confidence,takahashi2018statistical},
and provide an explicit construction of the `direct observations' model in the sense that $\what{\theta}_{\lambda}^{d}$ is approximately distributed 
as $\cN(\theta^\star,\tau_*\id)$.
We let $\mu^{(d)}_{\lambda}$ be the law of the couple $\big(\Theta + \tau_*(\lambda) Z, \ \Theta \big)$,
where $(\Theta,Z) \sim \hat{\mu}_{\theta^{\star}} \otimes \cN(0,1)$.
\begin{theorem}\label{th:law_debiased}
	Let $\what{\mu}_{(\what{\theta}_{\lambda}^d, \theta^{\star})}$ denote the empirical distribution (on $\R^2$) of the entries of $(\what{\theta}_{\lambda}^d, \theta^{\star})$.
	There exists constants $c,C > 0$ such that for all $\epsilon \in (0,1]$,
	$$
	\sup_{\theta^{\star} \in \cF_4(\xi)} \P \Big(\sup_{\lambda \in [\lambda_{\rm min},\lambda_{\rm max}]} W_2(\what{\mu}_{(\what{\theta}_{\lambda}^d,\theta^{\star})},\mu^{(d)}_{\lambda}) \geq \epsilon \Big) \leq \frac{C}{\epsilon^{11}} e^{-c N \epsilon^{17}}.
	$$
\end{theorem}
Theorem~\ref{th:law_debiased} is proved in Section~\ref{sec:proof_law_debiased}.

\section{Applications}
\label{sec:Applications}

\subsection{Estimation of the risk and the noise level}

In order to select the regularization parameter and to evaluate the quality of the Lasso solution $\htheta_{\lambda}$,
it is useful to estimate the risk and noise level. The paper~\cite{bayati2013estimating} developed a suite of estimators of these quantities based
on the asymptotic theory of~\cite{bayati2012lasso}. The same paper also proposed generalizations of these estimators to correlated designs.
Here we revisit these estimators and prove stronger guarantees. First, we obtain quantitative bound on the consistency rate
of our estimators. Second, our results are uniform over $\lambda$, which justifies using these estimators to select $\lambda$.

Let us start with the estimation of $\tau_*(\lambda)$ which plays a crucial role in the asymptotic theory.
We define
$$
\what{\tau}(\lambda) = \sqrt{n}\frac{\|y - X\what{\theta}_{\lambda} \|}{n-\|\what{\theta}_{\lambda}\|_0} \,.
$$
We will see with Theorem~\ref{th:sparsity_uniform} presented in Appendix~\ref{sec:proof_sparsity} that 
$$
\lim_{N,n\to\infty} \frac{1}{N}\|\what{\theta}_{\lambda}\|_0 = \P(|\Theta+\tau_*Z|\ge \tau_*\lambda/\beta_*)\equiv s_*(\lambda)\, .
$$
Further, by Theorem~\ref{th:unif_lambda_risk},
we have $\frac{1}{\sqrt{n}}\| y - X \what{\theta}_{\lambda}\| = \beta_*(\lambda) + o_n(1)$.
Recall that by~\eqref{eq:state_evolution} we have $\beta_*(\lambda) = \tau_*(\lambda)\big(1 - \frac{1}{\delta} s_*(\lambda)\big)$.
We deduce $\what{\tau}(\lambda) = \tau_*(\lambda) + o_n(1)$. More precisely we have the following consistency result. 
\begin{corollary}\label{cor:estim_tau}
	Assume here that $\mathcal{D}$ is either $\cF_0(s)$ or $\cF_p(\xi)$ for some $0 \leq s < s_{\rm max}(\delta)$ and $\xi >0, p >0$.
	There exists constants $C,c>0$ that only depend on $\Omega$ such that
	for all $\epsilon \in (0,1]$
	$$
	\sup_{\theta^{\star} \in \mathcal{D}} \
	\P \left( \sup_{\lambda \in [\lambda_{\rm min},\lambda_{\rm max}]}\left| \what{\tau}(\lambda) - \tau_*(\lambda) \right| \geq \epsilon  \right) \leq 
	C \epsilon^{-6} N^q\exp\left(-cN \epsilon^6  \right) \,,
	$$
	where $q = 0$ if\, $\mathcal{D} = \cF_0(s)$ and $q=(1/p - 1)_+$ if\, $\mathcal{D} = \cF_p(\xi)$.
\end{corollary}

We next consider estimating the $\ell_2$  error of the Lasso. Following~\cite{bayati2012lasso}, we define
$$
\widehat{R}(\lambda) = \what{\tau}(\lambda)^2 \Big(\frac{2}{N}\|\what{\theta}_{\lambda}\|_0 -1 \Big) + \frac{\big\|X^{\sT}(y-X\what{\theta}_{\lambda})\big\|^2}{N\big(1-\frac{1}{n}\|\what{\theta}_{\lambda}\|_0\big)^2} \,.
$$
\begin{corollary}\label{cor:estim_risk}
	Assume here that $\mathcal{D}$ is either $\cF_0(s)$ or $\cF_p(\xi)$ for some $0 \leq s < s_{\rm max}(\delta)$ and $\xi >0, p >0$. There exists constants $C,c >0$ such that for all $\epsilon \in (0,1]$,
	$$
	\sup_{\theta^{\star} \in \mathcal{D} } \P \Big( \sup_{\lambda \in [\lambda_{\rm min},\lambda_{\rm max}]}\Big| \what{R}(\lambda) - \frac{1}{N} \| \what{\theta}_{\lambda} - \theta^{\star} \|^2 \Big| \geq \epsilon \Big) \leq \frac{C}{\epsilon^6} N^qe^{-cN\epsilon^6} \,,
	$$
	where $q = 0$ if\, $\mathcal{D} = \cF_0(s)$ and $q=(1/p - 1)_+$ if\, $\mathcal{D} = \cF_p(\xi)$.
\end{corollary}
Corollary~\ref{cor:estim_risk} is proved in Appendix~\ref{sec:proof_estim_risk}.
Since by Corollary~\ref{cor:estim_risk}, Corollary~\ref{cor:estim_tau}, Theorem~\ref{th:unif_lambda_risk} we have with high probability $\what{R}(\lambda) \simeq \frac{1}{N} \| \what{\theta}_{\lambda} - \theta^{\star} \|^2 \simeq \delta (\tau_*(\lambda)^2 - \sigma^2) \simeq \delta (\what{\tau}(\lambda)^2 -\sigma^2)$, the estimator
\begin{equation}\label{eq:def_sigma_hat}
	\what{\sigma}^2(\lambda) = \what{\tau}(\lambda)^2 - \frac{N}{n} \what{R}(\lambda)
	= \what{\tau}(\lambda)^2 \Big(1 + \frac{N}{n} - \frac{2}{n} \|\what{\theta}_{\lambda}\|_0 \Big)- \frac{\big\|X^{\sT}(y-X\what{\theta}_{\lambda})\big\|^2}{n\big(1-\frac{1}{n}\|\what{\theta}_{\lambda}\|_0\big)^2}
\end{equation}
is a consistent estimator of the noise level $\sigma^2$.
\begin{corollary}\label{th:estim_sigma}
	There exists constants $C,c>0$ that only depend on $\Omega$, such that
	for all $\epsilon \in (0,1]$
	$$
	\sup_{\theta^{\star} \in \mathcal{D}} \ \P \left(
		\sup_{\lambda \in [\lambda_{\rm min},\lambda_{\rm max}]}
		\big| \what{\sigma}^2(\lambda) - \sigma^2 \big| > \epsilon 
	\right) 
	\leq \frac{C}{\epsilon^6} N^q e^{-cN\epsilon^6} \,.
	$$
\end{corollary}

Finally, we consider the prediction error  $\|X \theta^{\star} - X \what{\theta}_{\lambda} \|$.
Stein Unbiased Risk Estimator (SURE) provides a general method to estimate the prediction error, see e.g. 
\cite{stein1981estimation,efron2004estimation,tibshirani2012degrees}. In the present case, it takes the form
\begin{equation} \label{eq:def_SURE}
	\what{P}^{\rm SURE}(\lambda) = \frac{1}{n} \| y - X\what{\theta}_{\lambda} \|^2 + \frac{2 \sigma^2}{n} \|\what{\theta}_{\lambda} \|_0 \,.
\end{equation}
Tibshirani and Taylor~\cite{tibshirani2012degrees} proved that $\what{P}^{\rm SURE}(\lambda)$ is an unbiased estimator of the prediction error, namely
\begin{align}
	\E\{\what{P}^{\rm SURE}(\lambda) \}= \frac{1}{n} \|X \theta^{\star} - X \what{\theta}_{\lambda} \|^2 + \sigma^2\, .
\end{align}
The next result establishes consistency, uniformly over $\lambda$ and $\theta^\star$, with quantitative concentration estimates. 
\begin{corollary}\label{cor:sure}
	Assume here that $\mathcal{D}$ is either $\cF_0(s)$ or $\cF_p(\xi)$ for some $0 \leq s < s_{\rm max}(\delta)$ and $\xi >0, p >0$.
	There exists constants $C,c>0$ that only depend on $\Omega$ such that
	for all $\epsilon \in (0,1]$
	$$
	\sup_{\theta^{\star} \in \mathcal{D}} \P \Big(
		\sup_{\lambda \in [\lambda_{\rm min},\lambda_{\rm max}]} \Big| \frac{1}{n} \|X \theta^{\star} - X \what{\theta}_{\lambda} \|^2 + \sigma^2 - \what{P}^{\rm SURE}(\lambda) \Big| \geq \epsilon
	\Big) \leq \frac{C}{\epsilon^6} N^q e^{-cN\epsilon^6} \,,
	$$
	where $q = 0$ if\, $\mathcal{D} = \cF_0(s)$ and $q=(1/p - 1)_+$ if\, $\mathcal{D} = \cF_p(\xi)$.

	The same result holds if $\sigma$ in~\eqref{eq:def_SURE} is replaced by an estimator of the noise level satisfying the same consistency condition 
	as  $\what{\sigma}$ defined by~\eqref{eq:def_sigma_hat} (cf. Corollary~\ref{th:estim_sigma}). 
\end{corollary}
This corollary follows simply from Theorem~\ref{th:sparsity_uniform} and Theorem~\ref{th:unif_lambda_risk}.

\begin{remark}
	Notice that exact unbiasedness of $\what{P}^{\rm SURE}(\lambda)$ only holds if the noise $z$ in the linear model (\ref{eq:LinearModel}) is Gaussian~\cite{tibshirani2012degrees}.
	In contrast, it is not hard to generalize the proofs in the present paper to include other noise distributions. 
\end{remark}

\subsection{Adaptive selection of \texorpdfstring{$\lambda$}{lambda}}

As anticipated, we can use our uniform bounds to select $\lambda$
through an adaptive procedure. We discuss here three such procedures, that have already been illustrated in Figure~\ref{fig:comp_lambda_delta}:
$(i)$~Selecting $\lambda$ by minimizing the estimate $\what{\tau}(\lambda)$, we denote this by  $\what{\lambda}^{\rm EST}$; $(ii)$~Select $\lambda$ as to minimize
Stein's Unbiased Risk Estimate $\what{P}^{\rm SURE}(\lambda)$,
$\what{\lambda}^{\rm SURE}$; $(iii)$~Select $\lambda$ by $k$-fold
cross-validation, $\what{\lambda}^{k\text{-CV}}$. We will next describe these procedures in greater detail, and state the corresponding guarantees.

\noindent\emph{Minimization of $\, \what{\tau}(\lambda)$.} Since the $\ell_2$ risk of the Lasso is by Theorem~\ref{th:unif_lambda_risk} approximately equal to
$R_*(\lambda) = \delta(\tau_*(\lambda)^2 - \sigma^2)$ and since by 
Corollary~\ref{cor:estim_tau}, $\what{\tau}$ is a consistent estimator (uniformly in $\lambda$) of $\tau_*$, a natural procedure for selecting $\lambda$ is to minimize $\what{\tau}$.
We then define
$$
\what{\lambda}^{\rm EST} = \underset{\lambda \in [\lambda_{\rm min},\lambda_{\rm max}]}{\text{arg\,min}} \what{\tau}(\lambda) \,.
$$

The next result is an immediate consequence of Theorem~\ref{th:unif_lambda_risk} and Corollary~\ref{cor:estim_tau}:
\begin{proposition}
	Assume here that $\mathcal{D}$ is either $\cF_0(s)$ or $\cF_p(\xi)$ for some $0 \leq s < s_{\rm max}(\delta)$ and $\xi >0, p >0$.
	There exists constants $C,c>0$ that only depend on $\Omega$ such that
	for all $\epsilon \in (0,1]$
	$$
	\inf_{\theta^{\star} \in \mathcal{D}}
	\P\left(
		\frac{1}{N} \| \what{\theta}_{\, \what{\lambda}^{\rm EST}} - \theta^{\star} \|^2  
		\, \leq \, \inf_{\lambda \in [\lambda_{\rm min},\lambda_{\rm max}]} \Big\{ \frac{1}{N} \| \what{\theta}_{\lambda} - \theta^{\star} \|^2 \Big\}
		+ \epsilon
	\right) \geq
	1- C \epsilon^{-6} N^q \exp\left(-cN \epsilon^6  \right),
	$$
	where $q = 0$ if\, $\mathcal{D} = \cF_0(s)$ and $q=(1/p - 1)_+$ if\, $\mathcal{D} = \cF_p(\xi)$.
\end{proposition}

\noindent\emph{Minimization of SURE.} We define 
$$
\what{\lambda}^{\rm SURE} = 
\underset{\lambda \in [\lambda_{\rm min},\lambda_{\rm max}]}{\text{arg\,min}} \what{P}^{\rm SURE}(\lambda) \,.
$$
Here, it is understood that we can use either $\sigma$ or $\what{\sigma}(\lambda)$, cf. Eq.~\eqref{eq:def_sigma_hat},
in the definition of $\what{P}^{\rm SURE}$. We deduce from Corollary~\ref{cor:sure}:

\begin{proposition}\label{prop:lambda_sure}
	Assume here that $\mathcal{D}$ is either $\cF_0(s)$ or $\cF_p(\xi)$ for some $0 \leq s < s_{\rm max}(\delta)$ and $\xi >0, p>0$.
	There exists constants $C,c>0$ that only depend on $\Omega$ such that
	for all $\epsilon \in (0,1]$
	$$
	\inf_{\theta^{\star} \in \mathcal{D}}
	\P\left(
		\frac{1}{n} \| X \what{\theta}_{\, \what{\lambda}^{\rm SURE}} - X \theta^{\star} \|^2  
		\, \leq \, \inf_{\lambda \in [\lambda_{\rm min},\lambda_{\rm max}]} \Big\{ \frac{1}{n} \| X \what{\theta}_{\lambda} - X \theta^{\star} \|^2 \Big\}
		+ \epsilon
	\right) \geq
	1- C \epsilon^{-6} N^q \exp\left(-cN \epsilon^6  \right),
	$$
	where $q = 0$ if\, $\mathcal{D} = \cF_0(s)$ and $q=(1/p - 1)_+$ if\, $\mathcal{D} = \cF_p(\xi)$.
\end{proposition}

\noindent\emph{Cross-validation.} We analyze now $k$-fold Cross Validation. Let $k\ge 2$ and define $n_k = n(k-1)/k$.
We partition the rows of $X$ in $k$ groups: we obtain $k$-submatrices of size $(n/k) \times N$ that we denote $X^{(1)}, \dots, X^{(k)}$. Let us also write for $i \in \{1,\dots,k\}$, $X^{(\text{-}i)}$ for the submatrix of $X$ obtained by removing the rows $X^{(i)}$. We denote by $y^{(i)}$, $z^{(i)}$ and $y^{(\sminus i)}$, $z^{(\sminus i)}$ the corresponding subvectors of $y$ and $z$.

The estimator $\what{R}^{k\text{-CV}}$ of the risk using $k$-fold cross validation if defined as follows.
For $i =1, \dots, k$ solve the Lasso problem 
$$
\what{\theta}_{\lambda}^{i} = \underset{\theta \in \R^N}{\text{arg\,min}} \left\{
	\frac{1}{2 n_k} \Big\|y^{(\sminus i)} - X^{(\sminus i)}\theta\Big\|^2 + \frac{\lambda}{n} |\theta|
\right\} \,,
$$
and then compute
$$
\what{R}^{k\text{-CV}}(\lambda) = \frac{1}{N} \sum_{i=1}^k \Big\|y^{(i)} - X^{(i)}\what{\theta}_{\lambda}^{i}\Big\|^2 \,.
$$
Finally, we set $\lambda$ as follows
$$
\what{\lambda}^{k\text{-CV}} = 
\underset{\lambda \in [\lambda_{\rm min},\lambda_{\rm max}]}{\text{arg\,min}} \what{R}^{k\text{-CV}}(\lambda) \,.
$$

The next Proposition shows that $\what{R}^{k \text{-CV}}(\lambda)$ is equal to the true risk (shifted by $\delta \sigma^2$) up to $O(k^{-1/2})$.

\begin{proposition}\label{prop:riskCV}
	There exists constants $c,C > 0$ that depend only on $\Omega$, such that for all $k \geq 2$ such that $s_{\rm max}\big((k-1) \delta / k\big) > s$ in the case where $\mathcal{D} = \cF_0(s)$, we have
	$$
	\sup_{\theta^{\star} \in \mathcal{D}}
	\P \left(
		\sup_{\lambda \in [\lambda_{\rm min},\lambda_{\rm max}]} 
		\Big| \what{R}^{k\text{-CV}}(\lambda) - \frac{1}{N} \| \what{\theta}_{\lambda} - \theta^{\star} \|^2 - \delta \sigma^2 \Big| \geq \frac{C}{\sqrt{k}}
	\right)
	\leq C k^6 N^q e^{-c N / k^6} \,,
	$$
	where $q = 0$ if\, $\mathcal{D} = \cF_0(s)$ and $q=(1/p - 1)_+$ if\, $\mathcal{D} = \cF_p(\xi)$.
\end{proposition}
Proposition~\ref{prop:riskCV} is proved in Appendix~\ref{app:proof_riskCV}.
It follows from Proposition~\ref{prop:riskCV} that with high probability, 
$$\frac{1}{N} \|\what{\theta}_{\what{\lambda}^{k \text{- CV}}} - \theta^{\star} \|^2 \leq \inf_{\lambda \in [\lambda_{\rm min},\lambda_{\rm max}]} \frac{1}{N} \| \what{\theta}_{\lambda} - \theta^{\star} \|^2 + O(k^{-1/2}).$$

\subsection{Numerical experiments}



In this Section we compare numerically various different choices for the regularization parameter $\lambda$, namely $\what{\lambda}^{\rm EST}$, $\what{\lambda}^{\rm SURE}$ and $\what{\lambda}^{k\text{-CV}}$, presented in the previous section.
For these experiments we take the components $\theta_1^{\star}, \dots, \theta_N^{\star}$ to be i.i.d.\ from
$$
P_0 = s \cN(0,1) + (1-s) \delta_0 \,.
$$
Within this probabilistic model, we can compare achieved by our various choice of $\lambda$ to the
Bayes optimal error  (Minimal Mean Squared Error):
$$
\MMSE_N = \min_{\what{\theta}} \E \Big[\big\| \theta^{\star} - \what{\theta}(y,X) \big\|^2 \Big]
= \E \Big[ \big\| \theta^{\star} - \E [\theta^{\star}| y,X] \big\|^2\Big] \,,
$$
where the minimum is taken over all estimators $\what{\theta}$ (i.e.\ measurable functions of $X,y$).
The limit of the MMSE has been recently computed by~\cite{barbier2016CS} and~\cite{reeves2016replica}. Recall, that given two random variables
$U,V$, their mutual information is the Kullback-Leibler divergence between their joint distribution and the product of the marginals:
$I(U;V) \equiv D_{\mbox{\tiny \rm KL}}(p_{U,V}\|p_U\times p_V)$.

\begin{theorem}[Information-theoretic limit, from~\cite{barbier2016CS,reeves2016replica}]
	Define the function
	$$
	\Psi_{\delta,\sigma}(m) = I_{P_0} \Big(\frac{\sigma^{-2}}{1+m}\Big) + \frac{\delta}{2} \Big(\log(1+m)- \frac{m}{1+m}\Big) \,,
	$$
	where $I_{P_0}(r) = I(\Theta; \sqrt{r} \Theta + Z)$ for $(\Theta,Z) \sim P_0 \otimes \cN(0,1)$. Then, for almost every $\delta,\sigma >0$ the function $\Psi_{\delta,\sigma}$ admits a unique maximizer $m^*(\delta,\sigma)$ on $\R_{\geq 0}$ and 
	$$
	\MMSE_N \xrightarrow[N \to \infty]{} \delta \sigma^2 m^*(\delta,\sigma) \,.
	$$
\end{theorem}

Figure~\ref{fig:comp_lambda_delta} reports the risk achieved by the various choices of $\lambda$ as a function of the
number of samples per dimension $\delta$. We also compare the data-driven procedures of the previous section to the theory-driven choice
$\lambda = \sigma \sqrt{2 \log N}$.
In the left frame, we consider uncorrelated random designs: $X_{i,j} \iid \cN(0,1/n)$.  On  the right, we consider i.i.d. Gaussian rows with covariance structure
determined by an auto-regressive model. Explicitly, the columns  $(X_j)_{1 \leq j \leq N}$ of $X$ are generated according to:
\begin{equation}
X_1 = u_0 , \qquad X_{j+1} = \frac{1}{\sqrt{1 + \phi^2}} \big(\phi X_j + u_j \big) \label{eq:Auto}
\end{equation}
where $u_j \iid \cN(0,\id/n)$ and $\phi = 2$.
For both types of designs, $\what{\lambda}^{\rm EST}$, $\what{\lambda}^{\rm SURE}$ and $\what{\lambda}^{k\text{-CV}}$ perform similarly, and substantially
outperform the theoretical choice $\lambda = \sigma \sqrt{2 \log N}$.

For uncorrelated designs, the resulting risk is closely tracked by the asymptotic theory, and is surprisingly close to the 
asymptotic prediction for the Bayes risk $\MMSE_N$.

While our theory does not cover the case of correlated designs, the qualitative behavior is remarkably similar.
We also observed that in this case, the risk estimator $\what{R}(\lambda)$ is not consistent but its minimum is roughly located at the same value
of $\lambda$ as for uncorrelated designs.

\begin{figure*}[h!]
	\centering
	\includegraphics[width=1.0\linewidth]{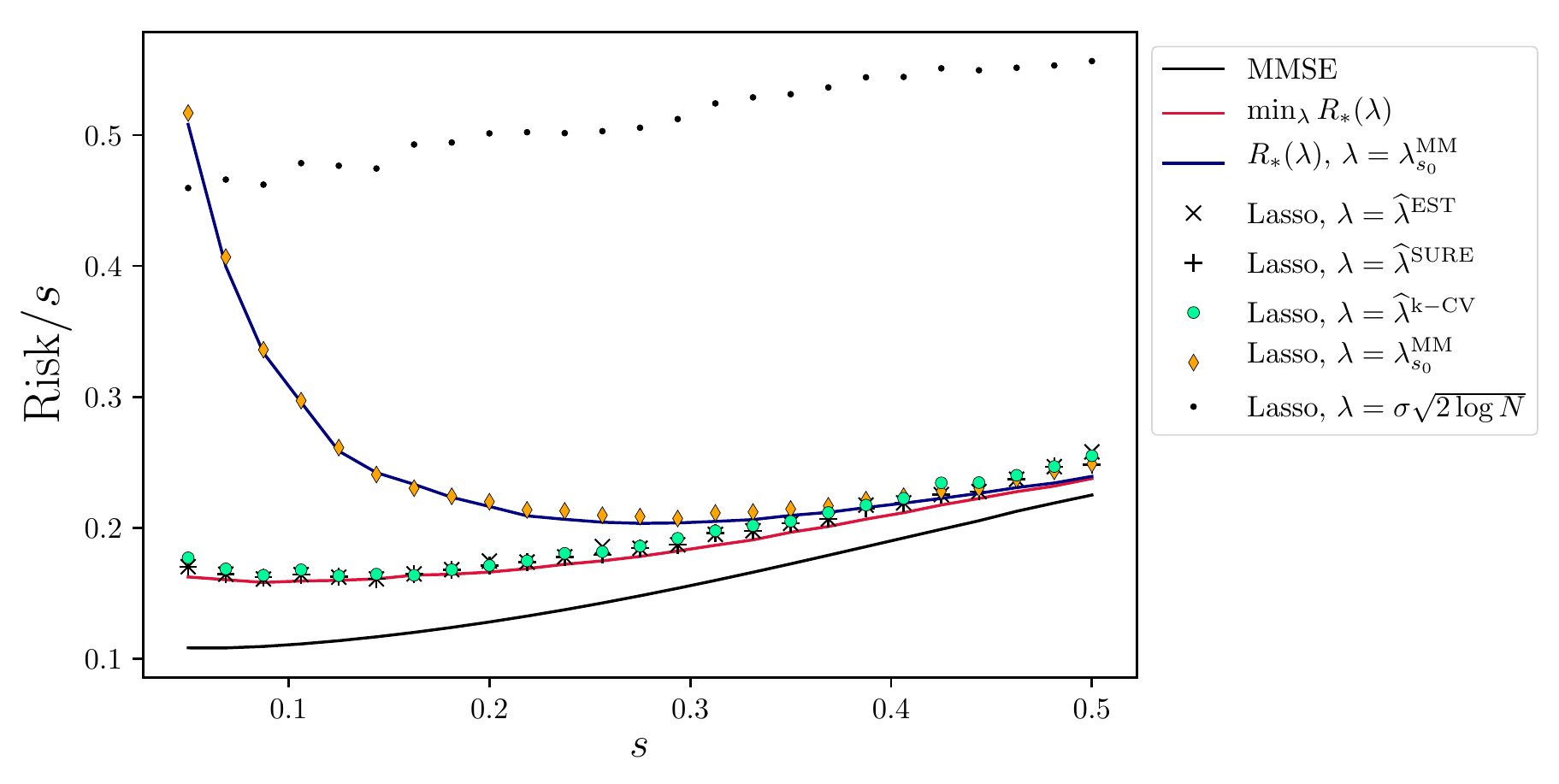}
	\caption{Risk of the Lasso for different choices of $\lambda$.
		$N=10000$, $\sigma = 0.2$, $\delta = 0.8$. Here $\theta^{\star}$ is chosen to be $sN$-sparse, and we vary the sparsity level $s$. 
The entries on the support of $\theta^{\star}$ are  i.i.d.\ $\cN(0,1)$. Cross-validation is carried out using $4$ folds. SURE is computed using the estimator $\what{\sigma}$. The minimax regularization $\lambda^{\rm MM}(s_0)$ is used at the nominal level $s_0 = 0.3$.
	}
	\label{fig:comp_lambda_s}
\end{figure*}

Next we study  adaptivity to sparsity. On Figure~\ref{fig:comp_lambda_s}, we plot the risk as a function of the sparsity of the signal $\theta^{\star}$.
We compare the three adaptive procedures  (namely, $\what{\lambda}^{\rm EST}$, $\what{\lambda}^{\rm SURE}$ and $\what{\lambda}^{k\text{-CV}}$), 
to the following choice
\begin{align*}
	\lambda^{\rm MM}(s_0) &= \alpha_0 \sigma \sqrt{1 - \frac{1}{\delta} M_{s_0}(\alpha_0)}\, ,\\
	M_s(\alpha) &= s (1+\alpha^2) + 2 (1-s) \big(
	(1+\alpha^2) \Phi(-\alpha) - \alpha \phi(\alpha) \big) \,,\\
	\alpha_0& = \arg\min_{\alpha\ge 0} M_{s_0}(\alpha)\, ,
\end{align*}
where $s_0 < s_{\rm max}(\delta)$ is a nominal value for the sparsity (in Figure~\ref{fig:comp_lambda_s}, we use $s_0=0.3$). 
The  value $\lambda^{\rm MM}(s_0)$ is expected to be asymptotically minimax optimal over $\cF_{0}(s_0)$~\cite{donoho2011noise}.

Also in this example, adaptive procedures dramatically outperform the fixed choice $\lambda = \sigma \sqrt{2 \log N}$,
and also the minimax optimal $\lambda$ at the nominal sparsity level.

\section{Proof strategy}
\label{sec:Proofs}

As mentioned above, our proofs are based on Gaussian comparison
inequalities, and in particular on Gordon's min-max theorem~\cite{gordon1985some,gordon1988milman}. In this section we review the
application of this inequality to the Lasso as developed in~\cite{thrampoulidis2015regularized}. We then discuss the 
limitations of earlier work, which does not characterize the empirical distribution of the Lasso estimator $\htheta_{\lambda}$ (or need extra sparsity assumptions~\cite{panahi2017universal})
nor  uniform bounds as in Theorem~\ref{th:unif_lambda_law}.
A key challenge is related to the fact that the Lasso cost function (\ref{eq:LassoCost}) is convex but not strongly convex. Hence, 
a small change in $\lambda$ could cause \emph{a priori} a large change in the minimizer $\htheta_{\lambda}$.

In order to overcome these problems, we establish a property
that we call `local stability.' Namely, if the empirical distribution of $(\htheta_{\lambda},\theta^{\star})$ deviates from 
our prediction, then the value of the optimization problem increases significantly. This implies that the empirical distribution
is stable with respect to perturbations of the cost  (e.g. changes in $\lambda$). Gordon's comparison is again crucial to prove this
stability property.

Finally, we describe how local stability is used to prove the theorems in the previous sections. A full description of the proofs is
provided in the appendices.

\subsection{Tight Gaussian min-max theorem}

It is more convenient (but equivalent) to study $\what{w}_{\lambda} = \what{\theta}_{\lambda} - \theta^{\star}$ instead of $\what{\theta}_{\lambda}$.
The vector $\what{w}_{\lambda}$ is the minimizer of the cost function
\begin{equation}
	\cC_{\lambda}(w) 
	= \frac{1}{2n} \left\| X w - \sigma z \right\|^2 + \frac{\lambda}{n} \big( |w + \theta^{\star}| - |\theta^{\star}| \big) \,.
\end{equation}
Following ~\cite{thrampoulidis2015regularized}, we rewrite the minimization of $\cC_{\lambda}$ as a saddle point problem:
\begin{equation}\label{eq:min_max_strategy}
	\min_{w \in \R^N} \cC_{\lambda}(w)=
	\min_{w \in \R^N} \max_{u \in \R^n} 
	\left\{
		\frac{1}{n} u^{\sT} \big( Xw - \sigma z\big) - \frac{1}{2n} \|u\|^2 + \frac{\lambda}{n}\big( |w + \theta^{\star}| - |\theta^{\star}| \big)
	\right\} \,.
\end{equation}
We apply the following Theorem from~\cite{thrampoulidis2015regularized} which improves over Gordon's Theorem~\cite{gordon1988milman} by exploiting 
convex duality.
\begin{theorem}[Theorem~3 from~\cite{thrampoulidis2015regularized}\,]\label{th:gordon}
	Let $S_w \subset \R^N$ and $S_u \subset \R^n$ be two compact sets and let $Q: S_w \times S_u \to\R$ be a continuous function. Let $G = (G_{i,j}) \iid \cN(0,1)$, $g \sim \cN(0,\bbf{I}_N)$ and $h \sim \cN(0,\bbf{I}_n)$ be independent standard Gaussian vectors. Define 
	$$
	\begin{cases}
		\cC^*(G) = \min\limits_{w \in S_w} \max\limits_{u \in S_u} u^{\sT} G w + Q(w,u)\,, \\
		L^*(g,h) = \min\limits_{w \in S_w} \max\limits_{u \in S_u} \|u\|_2 g^{\sT} w + \|w\|_2 h^{\sT} u + Q(w,u) \,.
	\end{cases}
	$$
	Then we have:
	\begin{itemize}
		\item For all $t \in \R$,
			$$
			\P \Big(\cC^*(G) \leq t \Big) \leq 2 \P \Big(L^*(g,h) \leq t \Big) \,.
			$$
		\item If $S_w$ and $S_u$ are convex and if $Q$ is convex concave, then for all $t \in \R$
			$$
			\P \Big(\cC^*(G) \geq t \Big) \leq 2 \P \Big(L^*(g,h) \geq t \Big) \,.
			$$
	\end{itemize}
\end{theorem}
For the reader's convenience, we provide in Appendix~\ref{sec:proof_gordon} a proof of this theorem.

Because of Gordon's Theorem, it suffices now to study (see Corollary~\ref{cor:gordon} below)
for $(g,g',h) \sim \cN(0,\bbf{I}_N) \otimes \cN(0,1) \otimes \cN(0,\bbf{I}_n)$.
\begin{equation}\label{eq:def_L}
	L_{\lambda}(w) = 
	\frac{1}{2}\left(\sqrt{\frac{\|w\|^2}{n} + \sigma^2}\ \frac{\|h\|}{\sqrt{n}}
	-	\frac{1}{n} g^{\sT} w  + \frac{g'\sigma}{\sqrt{n}} \right)_{\!\!+}^2
	+ \frac{\lambda}{n} |w+\theta^{\star}|
	- \frac{\lambda}{n} |\theta^{\star}| \,.
\end{equation}

\begin{corollary}\label{cor:gordon}
	\begin{enumerate}
		\item[$(a)$] Let $D \subset \R^N$ be a closed set. We have for all $t \in \R$
			$$
			\P\Big(\min_{w \in D} \cC_{\lambda}(w) \leq t\Big) \leq 2 \P \Big(\min_{w \in D} L_{\lambda}(w) \leq t\Big) \,.
			$$
		\item[$(b)$] Let $D \subset \R^N$ be a convex closed set. We have for all $t \in \R$
			$$
			\P\Big(\min_{w \in D} \cC_{\lambda}(w) \geq t\Big) \leq 2 \P \Big(\min_{w \in D} L_{\lambda}(w) \geq t\Big) \,.
			$$
	\end{enumerate}
\end{corollary}
\begin{proof} We  will only prove the first point, since the second follows from the same arguments.
	Define for $(w,u) \in \R^N \times \R^n$
	\begin{align*}
		c_{\lambda}(w,u) &= \frac{1}{n} u^{\sT} Xw - \frac{\sigma}{n} u^{\sT}z - \frac{1}{2n} \|u\|^2 + \frac{\lambda}{n}\big( |w + \theta^{\star}| - |\theta^{\star}| \big) \,,
		\\
		l_{\lambda}(w,u) &= 
		-\frac{1}{n^{3/2}} \|u\| g^{\sT} w 
		+\frac{1}{n} \|u\| g' \sigma 
		+\sqrt{\frac{\|w\|^2}{n} + \sigma^2} \frac{h^{\sT} u}{n} 
		- \frac{1}{2n} \|u\|^2 + 
		\frac{\lambda}{n} \big(  |w+\theta^{\star}| - |\theta^{\star}| \big) \,.
	\end{align*}
	Notice that for all $w \in \R^N$, $L_{\lambda}(w) = \max_{u \in \R^n} l_{\lambda}(w,u)$ and $\cC_{\lambda}(w) = \max_{u \in \R^n} c_{\lambda}(w,u)$.

	Let us suppose that $X,z,g,h,g'$ live on the same probability space and are independent. Let $\epsilon \in (0,1]$. Let $\sigma_{\rm max}(X)$ denote the largest singular value of the matrix $X$. By tightness we can find $K>0$ such that the event
	\begin{equation}\label{eq:event_norm_bounded}
		\Big\{
			\sigma_{\rm max}(X) \leq K, \quad \|z\| \leq K, \quad \|g\| \leq K, \quad \|h\| \leq K, \quad |g'| \leq K
		\Big\}
	\end{equation}
	has probability at least $1-\epsilon$. Let $D \subset \R^N$ be a (non-empty, otherwise the result is trivial) closed set. Let us fix $w_0 \in D$. On the event~\eqref{eq:event_norm_bounded} $\cC_{\lambda}(w_0)$ and $L_{\lambda}(w_0)$ are both upper bounded by some non-random quantity $R$. Let now $w \in D$ such that $\cC_{\lambda}(w) \leq R$. We have then $\frac{\lambda}{n} |w+\theta^{\star}| \leq R + \frac{\lambda}{n}|\theta^{\star}|$,
	which implies that $\|w\|$ is upper bounded by some non-random quantity $R_1$. This implies that, on the event~\eqref{eq:event_norm_bounded}, the minimum of $\cC_{\lambda}$ over $D$ is achieved on $D \cap B(0, R_1)$. Similarly on~\eqref{eq:event_norm_bounded} the minimum of $L_{\lambda}$ over $D$ is achieved on $D \cap B (0,R_2)$, for some non-random quantity $R_2$. Without loss of generalities, one can assume $R_1=R_2$. On the event~\eqref{eq:event_norm_bounded} we have
	$$
	\min_{w \in D} \cC_{\lambda}(w) = 
	\min_{w \in D \cap B(0,R_1)} \cC_{\lambda}(w) = 
	\min_{w \in D \cap B(0,R_1)} \max_{u \in B(0,R_3)}c_{\lambda}(w,u) \,,
	$$
	for some non-random $R_3>0$. This gives that for all $t \in \R$, we have
	$$
	\P \Big(\min_{w \in D} \cC_{\lambda}(w) \leq t \Big) \leq 
	\P \Big(\min_{w \in D \cap B(0,R_1)} \max_{u \in B(0,R_3)}c_{\lambda}(w,u) \leq t \Big) + \epsilon \,,
	$$
	and similarly 
	$$
	\P \Big(\min_{w \in D \cap B(0,R_1)} \max_{u \in B(0,R_3)} l_{\lambda}(w,u) \leq t \Big) \leq 
	\P \Big(\min_{w \in D } L_{\lambda}(w) \leq t \Big) + \epsilon \,.
	$$
	Since the sets $D \cap B(0,R_1)$ and $B(0,R_3)$ are compact, one can apply Theorem~\ref{th:gordon} to $c_{\lambda}$ and $l_{\lambda}$ and obtain:
	$$
	\P \Big(\min_{w \in D} \cC_{\lambda}(w) \leq t \Big) \leq 
	2\P \Big(\min_{w \in D } L_{\lambda}(w) \leq t \Big) + 2 \epsilon\,.
	$$
	The Corollary follows then from the fact one can take $\epsilon$ arbitrarily small.
\end{proof}

\subsection{Local stability}

In order to prove that (for instance) $\what{w}_{\lambda}$ verifies with high probability some property, let's say for instance that the empirical distribution
of $(\htheta_{\lambda}=\theta^{\star}+\what{w}_{\lambda},\theta^{\star})$ is close to $\mu_{\lambda}^\star$,
we define a set $D_{\epsilon} \subset \R^N$ that contains all the vectors that do not verify this property, e.g.\ $D_{\epsilon} = \big\{ w \in \R^N \, \big| \, W_2\big(\what{\mu}_{(\theta^{\star}+w,\theta^{\star})},\mu^*_{\lambda}\big)^2 \ge  \epsilon \big\}$, for some $\epsilon \in (0,1)$. The goal now is to prove that with high probability
$$
\min_{w \in D} \cC_{\lambda}(w) \geq \min_{w \in \R^N} \cC_{\lambda}(w) + \epsilon \,,
$$
for some $\epsilon>0$. Using Gordon's min-max Theorem (Corollary~\ref{cor:gordon}) we will be able to show
\begin{equation}\label{eq:sketch_gordon}
	\P \Big(
		\min_{w \in D_{\epsilon}} \cC_{\lambda}(w) \leq \min_{w \in \R^N} \cC_{\lambda}(w) + \epsilon
	\Big)
	\leq 
	2\P \Big(
		\min_{w \in D_{\epsilon}} L_{\lambda}(w) \leq \min_{w \in \R^N} L_{\lambda}(w) + \epsilon
	\Big) 
	+ o_N(1)\,.
\end{equation}
Informally, this is a consequence of the following two remarks. First, by applying parts $(a)$ and $(b)$ of Corollary~\ref{cor:gordon} to the convex domain $\reals^N$, we 
deduce that  $\min_{w \in \R^N} \cC_{\lambda}(w)\approx \min_{w \in \R^N} L_{\lambda}(w)$. Second, by applying part $(a)$ to the closed domain $D$,
we obtain $\min_{w \in D_{\epsilon}} \cC_{\lambda}(w) \gtrsim \min_{w \in D_{\epsilon}} L_{\lambda}(w)$

It remains now to study the cost function $L_{\lambda}$, which is much simpler. This is done in Appendix~\ref{sec:gordon_w}. 
The key step will be to establish the following `local stability' result (the next statement is an immediate consequence of  Proposition~\ref{prop:w_star} and Theorem~\ref{th:gordon_aux} in the appendices. We prove in fact that the cost function $L_{\lambda}$ is strongly convex on a neighborhood of its minimizer.).
\begin{theorem}\label{th:gordon_aux_simplified}
	The minimizer  $w^*_{\lambda}  =\arg\min_w L_{\lambda}(w)$ exists and is almost surely unique.
	Further, there exists constants $\gamma,c,C > 0$ that only depend on $\Omega$ such that for all $\theta^{\star} \in \mathcal{D}$, all $\lambda \in [\lambda_{\rm min}, \lambda_{\rm max}]$ and all $\epsilon \in (0,1]$
	$$
	\P \Big(
		\exists w\in\R^N, \quad 
		\frac{1}{N} \| w - w^*_{\lambda} \|^2 > \epsilon
		\quad \text{and} \quad
		L_{\lambda}(w) \leq \min\limits_{v \in \R^N} L_{\lambda}(v) + \gamma \epsilon 
	\Big)
	\leq \frac{C}{\epsilon} e^{-cn\epsilon^2} \,.
	$$
\end{theorem}

We do not obtain an equally strong result for the cost function $\cC_{\lambda}(w)$, but we prove 
the following statement, which is sufficient for obtaining uniform control (for the sake of argument, we focus here on the 
domain $\cF_{p}(\xi)$ and control of the empirical distribution).
\begin{theorem}\label{th:key_law}
	Assume that $\mathcal{D} = \cF_{p}(\xi)$ for some $\xi,p>0$.
	There exists constants $C,c,\gamma>0$ that only depend on $\Omega$ such that
	for all $\epsilon \in (0,\frac{1}{2}]$
	\begin{align*}
		\sup_{\lambda \in [\lambda_{\rm min},\lambda_{\rm max}]} \
		\sup_{\theta^{\star} \in \mathcal{D}} \
		\P \Big(\exists \theta \in \R^N, \quad W_2\big(\what{\mu}_{(\theta,\theta^{\star})},\mu^*_{\lambda}\big)^2 \geq \epsilon \quad &\text{and} \quad \mathcal{L}_{\lambda}(\theta) \leq \min \mathcal{L}_{\lambda} + \gamma \epsilon \Big) 
		\\
																																	   &\leq 
																																	   C \epsilon^{-\max(1,a)} \exp\left(-cN \epsilon^2 \epsilon^a \log(\epsilon)^{-2} \right) \,,
	\end{align*}
	where $a = \frac{1}{2} + \frac{1}{p}$.
\end{theorem}
Theorem~\ref{th:key_law} is proved in Appendix~\ref{sec:proof_key}.

\subsection{Sketch of proof of main results}

For the sake of simplicity, we will illustrate the prove strategy by considering the empirical distribution of $\what{w}_{\lambda} =\what{\theta}_{\lambda}-\theta^\star$,
as the argument is similar for other quantities. According to Theorem~\ref{th:unif_lambda_law}, this should be well approximated by $\omu_{\lambda}$ 
that is the law of $\what{\Theta}-\Theta$, when $(\what{\Theta},\Theta)\sim \mu_{\lambda}^*$, cf. Definition~\ref{def:mu_star}.

As anticipated, Eq.~(\ref{eq:sketch_gordon}) and Theorem~\ref{th:gordon_aux_simplified},
allow to control $W_2(\what{\mu}_{\what{w}_{\lambda}}, \omu_{\lambda})$
for a fixed $\lambda$ ($\what{\mu}_{\what{w}_{\lambda}}$ denotes the empirical distribution of the entries of $\what{w}_{\lambda}$). 
Namely, we can define $D_{\eps}$ to be
the set of vectors $w$ such that $W_2(\what{\mu}_w,\omu_{\lambda})\ge \eps>0$. We then prove that the minimizer $w^*_{\lambda}$ of $L_{\lambda}$ has empirical distribution
close to $\omu_{\lambda}$, and therefore by Theorem~\ref{th:gordon_aux_simplified},
$L_{\lambda}(w) > L_{\lambda}(w^*_{\lambda}) + \gamma \epsilon$ for all $w \in D_{\epsilon}$, with high probability. 
This imply that the right-hand side of~\eqref{eq:sketch_gordon} is very small and we deduce that, with high probability, all minimizers or near minimizers of
$\cC_{\lambda}(w)$ have empirical distribution close to $\omu_{\lambda}$,

We now would like to prove Theorem~\ref{th:unif_lambda_law} and show that with high probability $\what{\mu}_{\what{w}_{\lambda}}\approx \omu_{\lambda}$, 
uniformly in $\lambda \in [\lambda_{\rm min},\lambda_{\rm max}]$. To do so, we apply the above argument for 
$\lambda = \lambda_1, \dots, \lambda_k$, where $\lambda_1, \dots,\lambda_k$ is an $\epsilon$-net of $[\lambda_{\rm min},\lambda_{\rm max}]$. 
This implies that, with high probability  for $\lambda\in\{\lambda_1, \dots, \lambda_k\}$, $W_2(\what{\mu}_{\what{w}_{\lambda_i}},\omu_{\lambda_i})\le \eps$. 
Next, for $\lambda \in [\lambda_i,\lambda_{i+1}]$, we show that 
$$
\cC_{\lambda_i}(\what{w}_{\lambda}) = \min_{w \in \R^N} \cC_{\lambda_i}(w) + O(|\lambda_{i+1}- \lambda_i|) \,.
$$
Consequently if $|\lambda_{i+1}-\lambda_i| = O(\epsilon)$
(using again Eq.~(\ref{eq:sketch_gordon}) and Theorem~\ref{th:gordon_aux_simplified}), we obtain 
that $W_2(\what{\mu}_{\what{w}_{\lambda}},\omu_{\lambda_i}) = O(\epsilon)$ and therefore $W_2(\what{\mu}_{\what{w}_{\lambda}},\omu_{\lambda}) = O(\epsilon)$. We conclude that 
$W_2(\what{\mu}_{\what{w}_{\lambda}},\omu_{\lambda}) = O(\epsilon)$ for all $\lambda \in [\lambda_{\rm min},\lambda_{\rm max}]$, 
with high probability, which is the desired claim.

If the strategy exposed above allows to obtain the risk of the Lasso and the empirical distribution of its coordinates, it is not enough to get its sparsity $\|\what{\theta}_{\lambda}\|_0$ or to obtain the empirical distribution of the debiased lasso 
$$
\what{\theta}_{\lambda}^d = \what{\theta}_{\lambda} + \frac{X^{\sT}(y - X \what{\theta}_{\lambda})}{1 - \frac{1}{n}\|\what{\theta}_{\lambda}\|_0}.
$$
Therefore, we will need to analyze the vector
$$
\what{v}_{\lambda} = \frac{1}{\lambda} X^{\sT}(y - X \what{\theta}_{\lambda}) \,,
$$
which is a subgradient of the $\ell_1$-norm at $\what{\theta}_{\lambda}$. We are able to study $\what{v}_{\lambda}$ using Gordon's min-max Theorem because $\what{v}_{\lambda}$ is the unique maximizer of
$$
v \mapsto \min_{w \in \R^N} \Big\{
	\frac{1}{2n} \big\|Xw - \sigma z\big\|^2 + \frac{\lambda}{n} v^{\sT}(w+\theta^{\star})
\Big\} \,.
$$
The detailed analysis is done in Section~\ref{sec:gordon_v}.

\section*{Acknowledgements}

This work was partially supported by grants NSF DMS-1613091, NSF CCF-1714305 and NSF IIS-1741162
and ONR N00014-18-1-2729.

\begin{appendices}

	\section{Study of the scalar optimization problem}\label{sec:scalar}

	In this section we study the scalar optimization problem~\eqref{eq:max_min_scalar}:
	\begin{equation} 
		\max_{\beta \geq 0} 
		\min_{\tau \geq \sigma}
		\left(\frac{\sigma^2}{\tau} + \tau\right) \frac{\beta}{2}
		- \frac{1}{2} \beta^2
		+ 
		\frac{1}{\delta}
		\E
		\min_{w \in \R}\left\{
			\frac{w^2}{2\tau} \beta
			-	\beta Z w
		+ \lambda |w+\Theta| - \lambda |\Theta| \right\} \,,
\end{equation}
where $(\Theta,Z) \sim P_0 \otimes \cN(0,1)$, for some probability distribution $P_0$ with finite first moment: $\E_{P_0} |\Theta| < \infty$.
Of course, we will be mainly interested by the case where $P_0 =\what{\mu}_{\theta^{\star}}$, the empirical distribution of the entries of $\theta^{\star}$.
Define
$$
\psi_{\lambda}(\beta,\tau)
=
\left(\frac{\sigma^2}{\tau} + \tau\right) \frac{\beta}{2}
- \frac{1}{2} \beta^2
+ 
\frac{1}{\delta}
\E
\min_{w \in \R}\left\{
	\frac{w^2}{2\tau} \beta
	-	\beta Z w
+ \lambda |w+\Theta| - \lambda |\Theta| \right\} \,.
$$

\subsection{Basic properties of the scalar optimization problem}

\begin{lemma}[From~\cite{donoho2009message}\,]
	For all $\delta \in (0,1)$, the equation
	$$
	(1+\alpha^2) \Phi(-\alpha) - \alpha \phi(\alpha) = \frac{\delta}{2}
	$$
	admits a unique positive solution $\alpha_{\rm min} = \alpha_{\rm min}(\delta)>0$.
\end{lemma}
\begin{proof}
	Let $\varphi: \alpha \mapsto (1+\alpha^2) \Phi(-\alpha) - \alpha \phi(\alpha)$. $\varphi$ is continuous on $\R_{\geq 0}$, we have $\varphi(0) = \frac{1}{2}$ and $\varphi(+\infty) = 0$. It remains to show that $\varphi$ is strictly decreasing on $\R_{\geq 0}$.
	Compute $\varphi'(\alpha) = 2\alpha \Phi(-\alpha) - 2\phi(\alpha)$ and $\varphi''(\alpha) = 2 \Phi(-\alpha) > 0$.
	Since $\varphi'(+\infty) = 0$, we have that for all $\alpha \geq 0$, $\varphi'(\alpha) < 0$. $\varphi$ is thus strictly decreasing on $\R_{\geq 0}$.
\end{proof}
\\

Let us define
$$
\beta_{\rm max} = \beta_{\rm max}(\delta,\lambda) = \frac{\lambda}{\alpha_{\rm min}(\delta)} \,.
$$
More generally, we will always write $\alpha = \lambda / \beta$.
We prove in this section the following theorem and some auxiliary results.
\begin{theorem}\label{th:min_max_point}
	The max-min~\eqref{eq:max_min_scalar} is achieved at a unique couple $(\beta_*,\tau_*)$ and $0 < \beta_* < \beta_{\rm max}$. Moreover, $(\tau_*,\beta_*)$ is also the unique couple in $(0,+\infty)^2$ that verify
	\begin{equation}\label{eq:system_se}
		\begin{cases}
			\tau^2 &= \sigma^2 + \frac{1}{\delta} \E \Big[(\eta(\Theta + \tau Z, \tau \frac{\lambda}{\beta})- \Theta)^2 \Big] \\
			\beta &= \tau \left(1 - \frac{1}{\delta} \E \left[ \eta'(\Theta + \tau Z,\frac{\tau \lambda}{\beta}) \right] \right) \,.
		\end{cases}
	\end{equation}
\end{theorem}

\begin{lemma}\label{lem:stupid_bound}
	$$
	-\frac{\lambda}{\delta}\E |\Theta| \leq \max_{\beta \geq 0} \min_{\tau \geq \sigma} \psi_{\lambda}(\beta,\tau) \leq \frac{\sigma^2}{2} \,.
	$$
\end{lemma}
\begin{proof}
	We have $\max_{\beta} \min_{\tau} \psi_{\lambda}(\beta,\tau) \geq \min_{\tau} \psi_{\lambda}(0,\tau) = -\frac{\lambda}{\delta} \E |\Theta|$. Then, by taking $w=0$ one get
	$$
	\max_{\beta\geq0} \min_{\tau \geq \sigma} \psi_{\lambda}(\beta,\tau)
	\leq
	\max_{\beta\geq0} \min_{\tau \geq \sigma} 
	\left(\frac{\sigma^2}{\tau} + \tau\right) \frac{\beta}{2}
	- \frac{1}{2} \beta^2 = \frac{\sigma^2}{2} \,.
	$$
\end{proof}

Define for $\alpha \geq 0$ and $y \in \R$,
$$
\ell_{\alpha}(y) = 
\min_{x \in \R} \left\{
	\frac{1}{2}(y-x)^2 + \alpha |x|
\right\} 
$$
and for $Z \sim \cN(0,1), x \in \R$, $\Delta_{\alpha}(x) = \E \big[ \ell_{\alpha}(x+Z) - \alpha |x| \big] \,.$

\begin{lemma}
	\begin{equation} \label{eq:min_Delta}
		\E \left[
			\min_{w \in \R}\left\{
				\frac{w^2}{2\tau} \beta
				-	\beta Z w
			+ \lambda |w+\Theta| - \lambda |\Theta| \right\}
		\right]
		= \tau \beta \E\left[\Delta_{\alpha}\Big(\frac{\Theta}{\tau}\Big)\right]
		-\frac{\beta \tau}{2} \,.
\end{equation}
where $\alpha = \lambda / \beta$.
\end{lemma}
\begin{proof}
	Let $\beta > 0$ and compute
	\begin{align*}
		\E \left[ 
			\min_{w \in \R}\left\{
				\frac{w^2}{2\tau} \beta
				-	\beta Z w
		+ \lambda |w+\Theta| - \lambda |\Theta| \right\} \right]
		&= \frac{-\beta \tau}{2} + \frac{\beta}{\tau}
		\E
		\min_{w \in \R}\left\{
			\frac{1}{2} (w -\tau Z)^2
		+ \frac{\tau\lambda}{\beta} |w+\Theta| - \frac{\tau\lambda}{\beta} |\Theta| \right\}
		\\
		&= \frac{-\beta \tau}{2} + \beta \tau
		\E
		\min_{w \in \R}\left\{
			\frac{1}{2} (w -Z)^2
		+ \alpha \Big|w+\frac{\Theta}{\tau}\Big| - \alpha \Big|\frac{\Theta}{\tau}\Big| \right\} \,,
	\end{align*}
	where $\alpha = \lambda / \beta$.
	Thus
	$$
	\E \left[ 
		\min_{w \in \R}\left\{
			\frac{w^2}{2\tau} \beta
			-	\beta Z w
	+ \lambda |w+\Theta| - \lambda |\Theta| \right\} \right]
	= \frac{-\beta \tau}{2} + \beta \tau
	\E \left[
		\Delta_{\alpha}\left(\frac{\Theta}{\tau}\right)
	\right] \,.
	$$
\end{proof}

\begin{lemma}\label{lem:limit_beta_border}
	\begin{itemize}
		\item If $\beta > \beta_{\rm max}$,\quad
			$\displaystyle
			\psi_{\lambda}(\beta,\tau) \xrightarrow[\tau \to +\infty]{} - \infty \,.
			$
		\item If $\beta = \beta_{\rm max}$,\quad
			$\displaystyle
			\psi_{\lambda}(\beta,\tau) \xrightarrow[\tau \to +\infty]{} - \frac{\beta^2}{2} - \frac{\lambda}{\delta} \E\big|\Theta\big|\,.
			$
	\end{itemize}
\end{lemma}
\begin{proof}
	By~\eqref{eq:min_Delta} and the fact that $\Delta_{\alpha}(0) = \frac{1}{2} + \alpha \phi(\alpha) - (\alpha^2 + 1) \Phi(-\alpha)$ by Lemma~\ref{lem:prop_Delta}, we get that for all $\beta,\tau >0$
	$$
	\psi_{\lambda}(\beta,\tau) = \frac{\sigma^2 \beta}{2 \tau} - \frac{\beta^2}{2} 
	+ \frac{\tau \beta}{\delta}\left(\frac{\delta}{2} + \alpha \phi(\alpha) - (\alpha^2 + 1 ) \Phi(-\alpha)\right) + \frac{\beta}{\delta} \xi_{\alpha}(\tau) \,,
	$$
	where $\alpha = \lambda / \beta$ and
	$$
	\xi_{\alpha}(\tau)=
	\tau\E \left[ \left(\Delta_{\alpha}\left(\frac{\Theta}{\tau}\right)
	- \Delta_{\alpha}(0) \right) \right] \,.
	$$
	Using the definition of $\beta_{\rm max}$: if $\beta > \beta_{\rm max}$ then $\alpha < \alpha_{\rm min}$ and therefore $\frac{\delta}{2} + \alpha \phi(\alpha) - (\alpha^2 + 1)\Phi(-\alpha) <0$. If $\beta = \beta_{\rm max}$, $\frac{\delta}{2} + \alpha \phi(\alpha) - (\alpha^2 + 1)\Phi(-\alpha) = 0$. It remains to compute the limit of $\xi_{\alpha}(\tau)$ as $\tau \to \infty$.

	Using the expression (see Lemma~\ref{lem:prop_Delta}) of the left-and right-derivatives of $\Delta_{\alpha}$ at $0$, we have almost-surely:
	$$
	\tau\left(\Delta_{\alpha}\left(\frac{\Theta}{\tau}\right)
	- \Delta_{\alpha}(0) \right) 
	\xrightarrow[\tau \to \infty]{}
	- \alpha |\Theta| \,.
	$$
	Suppose that $\E|\Theta| < \infty$.
	By Lemma~\ref{lem:prop_Delta}, $\Delta_{\alpha}$ is $\alpha$-Lipschitz. Consequently, for all $\tau > 0$:
	$$
	\left|
	\tau\left(\Delta_{\alpha}\left(\frac{\Theta}{\tau}\right)
	- \Delta_{\alpha}(0) \right) 
	\right| \leq \alpha |\Theta| \,.
	$$
	Since we have assumed that $\E |\Theta| < \infty$ we can apply the dominated convergence theorem to obtain that $\xi_{\alpha}(\tau) \xrightarrow[\tau \to +\infty]{} - \alpha \E |\Theta|$.
\end{proof}
\\

Define
$$
w^*(\alpha,\tau) = \eta\big(\Theta + \tau Z, \alpha \tau \big) - \Theta \,.
$$
$w^*(\alpha,\tau)$ is the minimizer of $w \mapsto \frac{w^2}{2\tau} \beta -	\beta Z w + \lambda |w+\Theta|$ (recall that we always write $\alpha = \lambda / \beta$).

\begin{lemma}\label{lem:tau_star}
	If $\beta \geq \beta_{\rm max}$ the equation
	\begin{equation}\label{eq:tau_fixed_point}
		\tau^2 = \sigma^2 + \frac{1}{\delta} \E \Big[ w^*(\alpha,\tau)^2 \Big] 
		= \sigma^2 + \frac{1}{\delta} \E \Big[(\eta(\Theta + \tau Z, \tau \frac{\lambda}{\beta})- \Theta)^2 \Big]
		\,.
	\end{equation}
	does not admits any solution on $(0,+\infty)$.
	For all $\beta \in (0, \beta_{\rm max})$, the function $\psi_{\lambda}(\beta,\cdot)$ admits a unique minimizer $\tau_*(\beta)$ on $(0,+\infty)$ that is also the unique solution of~\eqref{eq:tau_fixed_point}.
	Moreover, $\alpha \mapsto \tau_*(\alpha)$ is $\cC^{\infty}$ on $(\alpha_{\rm min},+\infty)$ and for all $\alpha > \alpha_{\rm min}$
	$$
	\left|\frac{\partial \tau_*}{\partial \alpha} (\alpha) \right| \leq (\alpha + 1)\frac{\tau_*(\alpha)^3}{\delta \sigma^2}  \,.
	$$
\end{lemma}
\begin{proof}
	Most of this lemma was already proved in~\cite{donoho2009message}, we however provide a full proof for completeness.
	We have to study the fixed point equation
	$$
	\tau^2 = \sigma^2 + \frac{1}{\delta}\E \left[w^*(\alpha,\tau)^2\right] = F_{\alpha}(\tau^2) \,,
	$$
	where $\alpha=\lambda / \beta$. We can compute $F_{\alpha}$ explicitly:
	$$
	F_{\alpha}(\tau^2)= \sigma^2 + \frac{\tau^2}{\delta} 
	\Big(1+\alpha^2 + 
		\E \left[ 
			(x^2 - \alpha^2 - 1)\big(\Phi(\alpha - x) - \Phi(-\alpha - x)\big) - (x+\alpha)\phi(\alpha-x) + (x-\alpha)\phi(x+\alpha)
		\right]
	\Big) \,,
	$$
	where we used the notation $x=\frac{\Theta}{\tau}$.
	We can then compute the derivatives:
	\begin{align*}
		F_{\alpha}'(\tau^2)
		&= \delta^{-1} (1+\alpha^2) \E \left[\Phi(x-\alpha) + \Phi(-x-\alpha)\right] - \delta^{-1} \E\left[(x+ \alpha)\phi(x-\alpha) - (x-\alpha)\phi(-x-\alpha)\right] \,,
		\\
		F_{\alpha}''(\tau^2) &= \frac{-1}{2\delta\tau^2} \E \left[x^3 (\phi(x-\alpha) - \phi(x+\alpha))\right] \leq 0 \,.
	\end{align*}
	$F_{\alpha}$ is therefore concave. By dominated convergence
	$$
	F'_{\alpha}(\tau^2) \xrightarrow[\tau \to +\infty]{} 
	\frac{2}{\delta} \left( (1+\alpha^2) \Phi(-\alpha) - \alpha \phi(\alpha) \right) 
	\begin{cases}
		< 1 & \text{if} \ \beta < \beta_{\rm max} \\
		\geq 1 & \text{if} \ \beta \geq \beta_{\rm max} \,.
	\end{cases}
	$$
	Since $F_{\alpha}(0) = \sigma^2 > 0$ and by concavity of $F_{\alpha}$, the fixed point equation admits a unique solution $\tau_*(\alpha)$ if and only is $\beta \in (0,\beta_{\rm max})$. In that case we have also $F_{\alpha}'(\tau_*(\alpha)) < 1$. 

	Let us now assume that $\beta \in (0,\beta_{\rm max})$.
	We have almost-surely
	$$
	\frac{\partial}{\partial \tau}
	\min_{w \in \R}\left\{
		\frac{w^2}{2\tau} \beta
		+	\beta Z w
	+ \lambda |w+\Theta| \right\}
	=
	-\frac{\beta}{2\tau^2} w^*(\alpha,\tau)^2 \,.
	$$
	Since $|w^*(\alpha,\tau)| \leq \alpha \tau + \tau |Z|$, we have by derivation under the expectation
	$$
	\frac{\partial}{\partial \tau} \psi_{\lambda} (\beta,\tau)=\frac{\beta}{2} - \frac{\beta \sigma^2}{2\tau^2} - \frac{\beta}{2\delta\tau^2}  \E \left[w^*(\alpha,\tau)^2\right]
	= \frac{\beta}{2\tau^2} \left(\tau^2 - \left(\sigma^2 + \frac{1}{\delta}\E \left[w^*(\alpha,\tau)^2\right]\right)\right) \,.
	$$
	Consequently, $\tau_*(\beta)$ is the unique minimizer of $\psi_{\lambda}(\beta,\cdot)$ over $(0,+\infty)$.

	Let us now compute $\frac{\partial \tau_*^2}{\partial \alpha}$.
	Since $F_{\alpha}$ is a $\cC^{\infty}$ function of $\tau^2$, one can apply the implicit function theorem to obtain that the mapping $\alpha \mapsto \tau_*(\alpha)^2$ is $\cC^{\infty}$ and moreover:
	\begin{equation} \label{eq:deriv_tau_star}
		\frac{\partial \tau_*^2}{\partial \alpha} (\alpha) 
		= 
		\frac{\frac{\partial F_{\alpha}}{\partial \alpha}(\tau_*^2(\alpha))}{1-F_{\alpha}'(\tau^2_*(\alpha))} \,.
\end{equation}
Compute
	$$
	\frac{\partial F_{\alpha}}{\partial \alpha} (\tau^2)
	= \frac{2\tau^2}{\delta} \E \left[\alpha (\Phi(-\alpha+x)+ \Phi(-\alpha-x)) - (\phi(\alpha-x) + \phi(\alpha+x))\right] \,.
	$$
	One verify easily that
	\begin{equation} \label{eq:control_der_F_alpha}
		-1 \leq -2 \phi(0) \leq 2 \alpha \Phi(-\alpha) - 2 \phi(\alpha) \leq \frac{\delta}{2\tau^2}\frac{\partial F_{\alpha}}{\partial \alpha} (\tau^2) \leq \alpha \,.
	\end{equation}
	By concavity on has that $F'_{\alpha}(\tau_*^2(\alpha))$ is smaller than the slope of the line between the points of coordinates $(0,\sigma^2)$ and $(\tau_*^2(\alpha),\tau_*^2(\alpha))$:
	\begin{equation} \label{eq:maj_F_slope}
		F_{\alpha}'(\tau_*^2(\alpha))
		\leq 1 - \frac{\sigma^2}{\tau_*^2(\alpha)} \,.
	\end{equation}
	From equations (\ref{eq:deriv_tau_star}-\ref{eq:control_der_F_alpha}-\ref{eq:maj_F_slope}) we get
	$$
	\left|\frac{\partial \tau_*^2}{\partial \alpha} (\alpha) \right| \leq 2 \frac{\tau_*(\alpha)^4}{\sigma^2 \delta} (\alpha + 1) \,.
	$$
	The result follows then from the fact that $\frac{\partial \tau_*^2}{\partial \alpha} (\alpha) = 2 \tau_*(\alpha) \frac{\partial \tau_*}{\partial \alpha} (\alpha)$.
\end{proof}
\\

Define now
$$
\Psi_{\lambda}: \beta \mapsto \min_{\tau \geq \sigma} \psi_{\lambda}(\beta,\tau) \,.
$$
\begin{lemma}\label{lem:der_Psi}
	The function $\Psi_{\lambda}$ is differentiable on $(0,\beta_{\rm max})$ with derivative
	\begin{align}
		\Psi_{\lambda}'(\beta) 
		&= \tau_*(\alpha) - \beta - \frac{1}{\delta} \E \left[ Z w^*(\alpha,\tau_*(\alpha))\right] \label{eq:compute_der_Psi}
		\\
		&= \tau_*(\alpha) \left(1 - \frac{1}{\delta} \E \left[\Phi\Big(\frac{\Theta}{\tau_*(\alpha)} - \alpha\Big) + \Phi\Big(-\frac{\Theta}{\tau_*(\alpha)} - \alpha\Big)\right]\right) - \beta \,.
		\label{eq:compute_der_Psi2}
	\end{align}
\end{lemma}
\begin{proof}
	$\Psi_{\lambda}$ is differentiable on $(0,\beta_{\rm max})$ (because of Lemma~\ref{lem:tau_star}) with derivative
	\begin{align}
		\Psi_{\lambda}'(\beta) &= \frac{1}{2}\left(\frac{\sigma^2}{\tau_*(\alpha)} + \tau_*(\alpha) \right) - \beta + \frac{1}{\delta} \left(\frac{1}{2\tau_*(\alpha)} \E [w^*(\alpha,\tau^*(\alpha))^2] - \E[Z w^*(\alpha,\tau_*(\alpha))]  \right)
		\\
		&= \tau_*(\alpha) - \beta - \frac{1}{\delta} \E \left[ Z w^*(\alpha,\tau_*(\alpha))\right] \,,
	\end{align}
	because of~\eqref{eq:tau_fixed_point}. The second equality follows by Gaussian integration by parts.
\end{proof}
\\

\begin{corollary}\label{cor:max_beta}
	The function $\Psi_{\lambda}$ achieves its maximum over $\R_{\geq 0}$ at a unique $\beta_* \in (0,\beta_{\rm max})$.
\end{corollary}
\begin{proof}
	$\Psi_{\lambda}$ is the minimum of a collection of $1$-strongly concave functions: it is therefore $1$-strongly concave and admits thus a unique maximizer $\beta_*$ over $\R_{\geq 0}$.
	By Lemma~\ref{lem:limit_beta_border} we know that $\beta_* < \beta_{\rm max}$.
	Indeed, notice that $\max_{\beta} \Phi_{\lambda}(\beta) \geq \Psi_{\lambda}(0) = -\frac{\lambda}{\delta}\E |\Theta|$.
	Lemme~\ref{lem:limit_beta_border} gives that $\beta_* \in [0, \beta_{\rm max})$, because $\Psi_{\lambda}(\beta_{\rm max}) \leq \Psi_{\lambda}(0) - \frac{1}{2} \beta_{\rm max}^2 < \Psi_{\lambda}(0)$.
	By dominated convergence:
	$$
	\E \left[\Phi\Big(\frac{\Theta}{\tau_*(\alpha)} - \alpha\Big) + \Phi\Big(-\frac{\Theta}{\tau_*(\alpha)} - \alpha\Big)\right] \xrightarrow[\beta \to 0^+]{} 0 \,.
	$$
	Indeed, when $\beta \to 0^+$, $\alpha = \lambda / \beta \to +\infty$ and $|\frac{\Theta}{\tau_*(\alpha)}| \leq \frac{|\Theta|}{\sigma}$.
	Therefore by Lemma~\ref{lem:der_Psi} we obtain
	$$
	\liminf_{\beta \to 0^+} \Psi_{\lambda}'(\beta) \geq \sigma > 0 \,.
	$$
	By concavity, we deduce that $\beta_* \in (0,\beta_{\rm max})$.
\end{proof}
\\

\begin{proposition}\label{prop:beta_lip}
	The function $\lambda \mapsto \beta_*(\lambda)$ is $\cC^{\infty}$ and is $2\alpha_{\rm min}^{-1}$-Lipschitz over $(0, +\infty)$.
	$\lambda \mapsto \alpha_*(\lambda)$ is $\cC^{\infty}$ over $(0,+\infty)$ and strictly increasing.
\end{proposition}
\begin{proof}
	Let us define $\gamma_*(\lambda) = \beta_*(\lambda) / \lambda$. $\gamma_*(\lambda)$ is the unique maximizer of
	$$
	\gamma \mapsto
	\min_{\tau \geq \sigma}
	\left(\frac{\sigma^2}{\tau} + \tau\right) \frac{\gamma}{2}
	- \frac{\lambda}{2} \gamma^2
	+ 
	\frac{1}{\delta}
	\E
	\min_{w \in \R}\left\{
		\frac{w^2}{2\tau} \gamma
		-	\gamma Z w
	+ |w+\Theta| - |\Theta| \right\}
	= h(\gamma) -\frac{\lambda}{2} \gamma^2 \,,
	$$
	where $h$ is a concave $\cC^{\infty}$ function on $\R_{>0}$. $\gamma_*(\lambda)$ is thus the unique solution of
	$$
	h'(\gamma) - \lambda \gamma = 0 \,,
	$$
	on $\R_{>0}$. $\gamma \mapsto h'(\gamma) - \lambda \gamma$ is $\cC^{\infty}$ with derivative $\gamma \mapsto h''(\gamma) - \lambda <0$. Consequently, the implicit function theorem gives that the mapping $\lambda \in \R_{>0} \mapsto \gamma_*(\lambda)$ is $\cC^{\infty}$ and that
	$$
	\frac{\partial \gamma_*}{\partial \lambda} (\lambda) = \frac{-\gamma_*(\lambda)}{\lambda - h''(\gamma_*(\lambda))} <0 \,.
	$$
	One deduces that $\lambda \mapsto \alpha_*(\lambda) = \gamma_*(\lambda)^{-1}$ is $\cC^{\infty}$ and strictly increasing and that $\lambda \mapsto \beta_*(\lambda) = \lambda \gamma_*(\lambda)$ is $\cC^{\infty}$. Moreover
	$$
	\left|\frac{\partial \beta_*}{\partial \lambda} (\lambda) \right|
	=
	\left|\lambda \frac{\partial \gamma_*}{\partial \lambda} (\lambda) + \gamma_*(\lambda)\right|
	\leq 2 \gamma_*(\lambda)
	\leq \frac{2}{\alpha_{\rm min}} \,.
	$$
\end{proof}

\begin{proof}[of Theorem~\ref{th:min_max_point}]
	By Corollary~\ref{cor:max_beta}, the maximum in $\beta$ in~\eqref{eq:max_min_scalar} is achieved at a unique $\beta_* \in (0,\beta_{\rm max})$. To this $\beta_*$ corresponds a unique $\tau_*(\beta_*)$ that achieves the minimum in~\eqref{eq:max_min_scalar}, by Lemma~\ref{lem:tau_star}. By~\eqref{eq:tau_fixed_point} and~\eqref{eq:compute_der_Psi2} we obtain that $(\tau_*(\beta_*),\beta_*)$ is solution of the system~\eqref{eq:system_se}. Let now $(\tau,\beta) \in (0,+\infty)^2$ be another solution of~\eqref{eq:system_se}. $\tau$ is therefore solution of~\eqref{eq:tau_fixed_point} which gives that $\beta \in (0,\beta_{\rm max})$ and $\tau = \tau_*(\beta)$ by Lemma~\ref{lem:tau_star}. The second equality in~\eqref{eq:system_se} gives that $\Psi_{\lambda}'(\beta)=0$ and thus that $\beta = \beta_*$ by strong concavity of $\Psi_{\lambda}$. We conclude $(\tau,\beta)= (\tau_*(\beta_*),\beta_*)$.
\end{proof}

\subsection{Control on \texorpdfstring{$\beta_*,\tau_*$}{beta*,tau*}}\label{sec:control_tau_beta}

The goal of this section is to show that $\beta_*$ and $\tau_*$ remain bounded when $\theta^{\star}$ varies in $\mathcal{D}$.

\begin{theorem}\label{th:control_beta_tau}
	There exists constants $\beta_{\rm min}, \tau_{\rm max} > 0$ that only depend on $\Omega$ such that for all $\theta^{\star} \in \mathcal{D}$ and all $\lambda \in [\lambda_{\rm min},\lambda_{\rm max}]$,
	$$
	\beta_{\rm min} \leq \beta_*(\lambda) < \beta_{\rm max} \quad \text{and} \quad
	\sigma \leq \tau_*(\lambda) \leq \tau_{\rm max} \,.
	$$
\end{theorem}

To prove Theorem~\ref{th:control_beta_tau}, we separate the case where $\mathcal{D} = \cF_p(\xi)$ (where it follows from Lemma~\ref{lem:psi_H} and Corollary~\ref{cor:control_tau_lp} below) from the case where $\mathcal{D} = \cF_0(s)$ (where it follows from Lemmas~\ref{lem:control_beta_l0} and~\ref{lem:control_tau_l0}).

\subsubsection{Technical lemmas}
\begin{lemma} \label{lem:intro_H}
	We have
	\begin{align*}
		\max_{\beta \geq 0} \min_{\tau \geq \sigma} \psi_{\lambda}(\beta,\tau)
		= \psi_{\lambda}(\beta_*,\tau_*(\beta_*)) 
		&= \frac{1}{2} \beta_*^2 + \frac{\lambda}{\delta}\E \big[ |w^*(\alpha_*,\tau_*(\beta)) + \Theta| - |\Theta| \big]
		\\
		&= \frac{1}{2} \beta_*^2 + \tau_*(\beta_*) \frac{\lambda}{\delta} \E \left[H_{\alpha_*}\!\!\left(\frac{\Theta}{\tau_*(\beta_*)}\right)\right] \,,
	\end{align*}
	where
	$$
	H_{\alpha}(x) = (x-\alpha) \Phi(- \alpha +x) + (-x-\alpha) \Phi(-x-\alpha) + \phi(-x+\alpha) + \phi(x+\alpha) - |x| \,.
	$$
\end{lemma}
\begin{proof}
	Using the optimality condition~\eqref{eq:tau_fixed_point} of $\tau_*(\beta)$, we have for all $\beta \in (0, \beta_{\rm max})$
	$$
	\psi_{\lambda}(\beta,\tau_*(\beta)) = - \frac{1}{2} \beta^2 + \beta \tau_*(\beta) - \frac{\beta}{\delta} \E \left[ Z w^*(\alpha,\tau_*(\beta)) \right] +\frac{\lambda}{\delta} \E \left[ |w^*(\alpha,\tau_*(\beta)) + \Theta| - |\Theta| \right] \,.
	$$
	At $\beta_*$ the optimality condition (see~\eqref{eq:compute_der_Psi}) reads $\beta_* = \tau_*(\beta_*) - \frac{1}{\delta}\E\big[ Z w^*(\alpha_*,\tau_*(\beta_*))\big]$, thus
	\begin{equation} \label{eq:cost_after_cond}
		\psi_{\lambda}(\beta_*,\tau_*(\beta_*)) =  \frac{1}{2} \beta_*^2  +\frac{\lambda}{\delta} \E \big[ |w^*(\alpha_*,\tau_*(\beta)) + \Theta| - |\Theta| \big] \,.
	\end{equation}
	Compute for $\alpha,\tau >0 $
	\begin{equation} \label{eq:fact_tau}
		\E \big|w^*(\alpha,\tau) + \Theta\big| = \E \big|\eta(\Theta + \tau Z, \alpha \tau)\big|
		= \tau \E \Big| \eta\Big(\frac{\Theta}{\tau} + Z,\alpha\Big)\Big| \,.
	\end{equation}
	Now, for $x \in \R$,
	\begin{align*}
		\E | \eta(x+Z,\alpha) |
		&= 
		\int_{\alpha - x}^{+\infty} (x+z-\alpha) \phi(z) dz
		+\int^{-\alpha - x}_{-\infty} (-x-z-\alpha) \phi(z) dz
		\\
		&=(x-\alpha)\Phi(x-\alpha) + (-x-\alpha)\Phi(-x-\alpha) + \phi(\alpha-x) +\phi(\alpha+x)
		\\
		&= H_{\alpha}(x) + |x| \,.
	\end{align*}
	and we obtain the Lemma by putting this together with~\eqref{eq:fact_tau} and~\eqref{eq:cost_after_cond}.
\end{proof}
\\

The next Lemma summarizes the main properties of $H_{\alpha}$.
\begin{lemma}\label{lem:prop_H}
	$H_{\alpha}$ is a continuous, even function and for $x > 0$
	$$
	H'_{\alpha}(x) = \Phi(x-\alpha) - \Phi(-x-\alpha) - 1 \in (-1,0) \,.
	$$
	$H_{\alpha}$ is therefore $1$-Lipschitz. $H_{\alpha}$ admits a maximum at $0$ and
	$$
	H_{\alpha}(0) = 2\phi(\alpha) - 2\alpha \Phi(-\alpha) > 0 \,.
	$$
	Moreover $H_{\alpha} (x) \xrightarrow[x \to +\infty]{} -\alpha$.
\end{lemma}

\subsubsection{On \texorpdfstring{$\ell_p$}{lp}-balls}

\begin{lemma} \label{lem:psi_H}
	Assume that $\E\big[|\Theta|^p\big]\leq \xi^p$ for some $\xi,p >0$.
	Then, there exists a constant $\beta_{\rm min} = \beta_{\rm min}(\delta,\lambda_{\rm min},\xi,p,\sigma)$ such that for all $\lambda \geq \lambda_{\rm min}$,
	$$
	0 < \beta_{\rm min} \leq \beta_*(\lambda) < \beta_{\rm max} \,.
	$$
\end{lemma}
\begin{proof}
	Let $\beta \in (0, \beta_{\rm max})$. By Lemma~\ref{lem:der_Psi} we have
	$$
	\Psi_{\lambda}'(\beta)
	=
	\tau_*(\beta) \left(1 - \frac{1}{\delta} \E \left[\Phi\Big(\frac{\Theta}{\tau_*(\beta)}-\alpha\Big) + \Phi\Big(-\frac{\Theta}{\tau_*(\beta)}-\alpha\Big)\right]\right) - \beta \,.
	$$
	The function $g_{\alpha} : x\mapsto \Phi(x-\alpha) + \Phi(-x-\alpha)$ is even, and increasing over $\R_{\geq 0}$. 
	Let $K>0$ such that $\frac{\xi^p}{K^p \sigma^p} \leq \frac{\delta}{4}$. By Markov's inequality we have
	$$
	\P\left(
		\left|\frac{\Theta}{\tau_*(\beta)}\right| \geq K
	\right)
	\leq
	\P\left(
		\left|\frac{\Theta}{\sigma}\right|^p \geq K^p
	\right)
	\leq \frac{1}{K^p\sigma^p} \E |\Theta|^p
	\leq \frac{\delta}{4} \,.
	$$
	Thus
	\begin{align*}
		\E \left[
			g_{\alpha}\left(\frac{\Theta}{\tau_*(\beta)}\right)
		\right]
		\leq 
		g_{\alpha}(K) + \frac{\delta}{4} \,.
	\end{align*}
	As $\beta \to 0$, $\alpha \geq \lambda_{\rm min} / \beta \to +\infty$.
	Since $g_{\alpha}(K) \xrightarrow[\alpha \to +\infty]{} 0$ there exists $\beta_0=\beta_0(K,\lambda_{\rm min},\delta) > 0$ such that for all $\beta \in (0,\beta_0)$, $g_{\alpha}(K) \leq \frac{\delta}{4}$. Thus for all $\beta \in (0,\beta_0)$,
	$$
	\Psi_{\lambda}'(\beta)
	\geq
	\tau_*(\beta) \left(1 -\frac{\frac{\delta}{4} + \frac{\delta}{4}}{\delta} \right) - \beta
	\geq \frac{\sigma}{2} - \beta \,.
	$$
	Let $\beta_{\rm min} = \min(\frac{\sigma}{2},\beta_0)$. We conclude that for all $\beta \in (0,\beta_{\rm min})$, $\Psi_{\lambda}'(\beta) > 0$.
	By concavity we have then that $\beta_* \geq \beta_{\rm min}$. The other inequality $\beta_* < \beta_{\rm max}$ was already proved in Corollary~\ref{cor:max_beta}.

\end{proof}

\begin{corollary}\label{cor:control_tau_lp}
	Assume that $\E\big[|\Theta|^p\big]\leq \xi^p$ for some $\xi,p >0$.
	Then there exists a constant $\tau_{\rm max}=\tau_{\rm max}(\xi,p,\delta,s,\lambda_{\rm min},\lambda_{\rm max})$ such that for all $\lambda \in [\lambda_{\rm min},\lambda_{\rm max}]$,
	$$
	\sigma \leq \tau_*(\beta_*(\lambda)) \leq \tau_{\rm max} \,.
	$$
\end{corollary}
\begin{proof}
	Let $t \geq \xi$. By Markov's inequality we have $\P(|\Theta| \geq t) \leq \left( \frac{\xi}{t} \right)^p\leq 1$, since $\theta^{\star} \in \cF_p(\xi)$.
	\begin{align*}
		\E\left[H_{\alpha_*}\! \left(\frac{\Theta}{\tau_*(\alpha_*)}\right)\right]
		&=
		\E\left[\bbf{1}(|\Theta|< t) H_{\alpha_*}\! \left(\frac{\Theta}{\tau_*(\alpha_*)}\right)\right]
		+ \E\left[\bbf{1}(|\Theta|\geq t) H_{\alpha_*}\left(\frac{\Theta}{\tau_*(\alpha_*)}\right)\right]
		\\
		&\geq  \left( 1 - \left(\frac{\xi}{t}\right)^{\! p}\right) \left(H_{\alpha_*}(0)- \frac{t}{\tau_*(\alpha_*)}\right)  - \alpha_* \left(\frac{\xi}{t}\right)^{\! p} \,,
	\end{align*}
	because by Lemma~\ref{lem:prop_H}, $H_{\alpha_*}$ is $1$-Lipschitz and for all $x \in \R$, $-\alpha_* \leq H_{\alpha_*}(x) \leq H_{\alpha_*}(0)$. Replacing $H_{\alpha_*}(0)$ by its expression given by Lemma~\ref{lem:prop_H} we get
	\begin{align*}
		\E\left[H_{\alpha_*}\! \left(\frac{\Theta}{\tau_*(\alpha_*)}\right)\right]
		&\geq   2 \big(\phi(\alpha_*) -  \alpha_* \Phi(-\alpha_*) \big)  - \left(\frac{\xi}{t}\right)^{\!p} \! \big(\alpha_* + 2(\phi(\alpha_*) - \alpha_* \Phi(-\alpha_*))\big) 
		- \frac{t}{\tau_*(\alpha_*)} \,.
	\end{align*}
	Since $\alpha_* \leq \lambda_{\rm max} / \beta_{\rm min}$ and $\phi(\alpha_*) -  \alpha_* \Phi(-\alpha_*) >0$ (because $\alpha_* > \alpha_{\rm min}$),
	we can find a constant $t=t(\delta,\sigma,\lambda_{\rm min},\lambda_{\rm max},p,\xi) \geq \xi$ such that
	$$
	(\alpha_* + 2(\phi(\alpha_*) - \alpha_* \Phi(-\alpha_*))) \left(\frac{\xi}{t}\right)^p
	\leq \phi(\alpha_*) - \alpha_* \Phi(-\alpha_*) \,.
	$$
	For this choice of $t$ we have then
	\begin{align*}
		\frac{\tau_*(\alpha_*) \lambda}{\delta} \E\left[H_{\alpha_*}\left(\frac{\Theta}{\tau_*}\right)\right]
		\geq  \frac{\lambda}{\delta} \tau_*(\alpha_*) (\phi(\alpha_*)-\alpha_* \Phi(-\alpha_*)) - \frac{\lambda t}{\delta} \,.
	\end{align*}
	Consequently by Lemma~\ref{lem:stupid_bound} and Lemma~\ref{lem:intro_H} we have
	$$
	\frac{\beta_*^2}{2} + \frac{\lambda}{\delta} \tau_* (\phi(\alpha_*)-\alpha_*\Phi(-\alpha_*)) 
	- \frac{\lambda t}{\delta} \leq \psi_{\lambda} (\beta_*,\tau_*(\alpha_*))
	\leq \frac{\sigma^2}{2}  \,,
	$$
	which finally gives
	$$
	\tau_*(\alpha_*) \leq \frac{\delta\sigma^2\lambda^{-1} + t}{\phi(\alpha_{\rm max}) - \alpha_{\rm max}\Phi(-\alpha_{\rm max})} \,.
	$$
\end{proof}

\subsubsection{On sparse balls}\label{sec:control_sparse_balls}

Define the critical function:
$$
M_s: \alpha \mapsto s (1+\alpha^2) + 2 (1-s) \big(
	(1+\alpha^2) \Phi(-\alpha) - \alpha \phi(\alpha)
\big) \,.
$$
$M_s$ corresponds to the worst mean squared error achievable by soft-thresholding with threshold $\alpha$ to estimate a vector $\theta^{\star} \in \cF_0(s)$ from the observations $y=\theta^{\star} + w$, where $w \sim \cN(0,I_N)$, see~\cite{donoho1994ideal,donoho1994minimax,johnstone2002function}.
\begin{lemma}[From~\cite{donoho2009message}\,]\label{lem:below_s_max}
	Assume that 
	$$
	s < s_{\rm max}(\delta) = 
	\delta \max_{\alpha \geq 0} \left\{
		\frac{1 - \frac{2}{\delta}\big((1+\alpha^2) \Phi(-\alpha) - \alpha \phi(\alpha)\big)}{1 + \alpha^2 - 2 \big((1+\alpha^2) \Phi(-\alpha) - \alpha \phi(\alpha)\big)}
	\right\} \,.
	$$
	Then there exists $\alpha \geq 0$ such that $M_s(\alpha) < \delta$.
\end{lemma}
\begin{proof}
	Let $s < s_{\rm max}(\delta)$. From the definition of $s_{\rm max}(\delta)$, we can find $\alpha \in \R$ such that
	$$
	\delta 
	\frac{1 - \frac{2}{\delta}\big((1+\alpha^2) \Phi(-\alpha) - \alpha \phi(\alpha)\big)}{1 + \alpha^2 - 2 \big((1+\alpha^2) \Phi(-\alpha) - \alpha \phi(\alpha)\big)}
	> s \,,
	$$
	which gives $M_s(\alpha) < \delta$.
\end{proof}
\\

We assume in this section that $s < s_{\rm max}(\delta)$.
Let us compute the derivatives
\begin{align*}
	M_s'(\alpha) &= 2 \left(\alpha s + 2 (1-s) (\alpha \Phi(-\alpha) - \phi(\alpha))\right) \,,
	\\
	M_s''(\alpha) &= 2 \left( s + (1-s) 2 \Phi(-\alpha) \right) > 0 \,.
\end{align*}
Notice that $M_s(\alpha) = \frac{1}{2} \big(\alpha M_s'(\alpha) + M_s''(\alpha)\big)$. 
Let $\alpha_0$ be the unique $\alpha >0$ such that $M_s'(\alpha)= 0$ and let $\alpha_1 < \alpha_2$ be such that $M_s(\alpha_1) = M_s(\alpha_2) = \delta$.
We can then easily plot the variations of $M_s$:
\vspace{3mm}

\begin{tikzpicture}
	\tkzTabInit{$\alpha$ /1, $M_s'$ /1, $M_s$ /1.5, $\frac{1}{2}M_s'' $/1.5} {$0$ ,$\alpha_1$, $\alpha_0$, $\alpha_2$, $+\infty$}
	\tkzTabLine{, ,- , , z, , +, }
	\tkzTabVar{+/ 1, R/, -/ $M_s(\alpha_0)$, R/, +/ $+\infty$}
	\tkzTabIma{1}{3}{2}{$\delta$}
	\tkzTabIma{3}{5}{4}{$\delta$}
	\tkzTabVar{+/ 1, R/, R/, R/, -/ $s$}
	\tkzTabIma{1}{5}{3}{$M_s(\alpha_0)$}
\end{tikzpicture}

\begin{lemma}\label{lem:control_beta_l0}
	Let $s < s_{\rm max}(\delta)$ and assume that $\P(\Theta \neq 0) \leq s$.
	Then, there exists a constant $\beta_{\rm min} = \beta_{\rm min}(\delta,\lambda_{\rm min},s,\sigma)$ such that for all $\lambda \geq \lambda_{\rm min}$
	$$
	0 < \beta_{\rm min} \leq \beta_*(\lambda) <  \beta_{\rm max}(\lambda_{\rm max},\delta) : =\lambda_{\rm max} / \alpha_{\rm min}(\delta)\,.
	$$
\end{lemma}
\begin{proof}
	We already proved in Corollary~\ref{cor:max_beta} that $\beta_*(\lambda) < \lambda / \alpha_{\rm min}$.
	For all $0 < \beta < \lambda / \alpha_{\rm min}$, we have by Lemma~\ref{lem:der_Psi}
	$$
	\Psi_{\lambda}'(\beta)
	=
	\tau_*(\beta) \left(1 - \frac{1}{\delta} \E \left[\Phi\Big(\frac{\Theta}{\tau_*(\beta)}-\alpha\Big) + \Phi\Big(-\frac{\Theta}{\tau_*(\beta)}-\alpha\Big)\right]\right) - \beta \,.
	$$
	The function $g_{\alpha} : x\mapsto \Phi(x-\alpha) + \Phi(-x-\alpha)$ is even, and increasing over $\R_{\geq 0}$. 
	Therefore
	$$
	\Psi_{\lambda}'(\beta)
	\geq
	\frac{\tau_*(\beta)}{\delta} \big(\delta - s - (1-s) 2 \Phi(-\alpha) \big) - \beta \,.
	$$
	Let $\beta_0 = \beta_0(\lambda_{\rm min},\delta,s) > 0$ such that for all $\beta \in (0,\beta_0)$, $2 \Phi(-\alpha) \leq \frac{1}{2}(\delta - s)$. For all $\beta \in (0,\beta_0)$ we have then
	$$
	\Psi_{\lambda}'(\beta) \geq \frac{\sigma (\delta - s)}{2 \delta} - \beta \,.
	$$
	Let $\beta_{\rm min} = \min(\frac{\sigma (\delta - s)}{2 \delta},\beta_0)$: for all $\beta \in (0,\beta_{\rm min})$, $\Psi_{\lambda}'(\beta) > 0$.
	By concavity we conclude that $\beta_* \geq \beta_{\rm min}$.
\end{proof}

\begin{lemma}\label{lem:control_tau_l0}
	Let $s < s_{\rm max}(\delta)$ and assume that $\P(\Theta \neq 0) \leq s$. Then for all $\beta,\tau,\lambda >0$ we have
	$$
	\psi_{\lambda}(\beta,\tau) \geq \frac{\beta \sigma^2}{2 \tau} - \frac{\beta^2}{2} + \frac{\tau \beta}{2 \delta} \big(\delta - M_s(\alpha)\big) \,.
	$$
\end{lemma}
\begin{proof}By~\eqref{eq:min_Delta} we have for all $\beta,\tau >0$
	$$
	\psi_{\lambda}(\beta,\tau) = \frac{\beta}{2}\left(\frac{\sigma^2}{\tau} + \tau\right) - \frac{\beta^2}{2} + \frac{\tau \beta}{\delta} \E\left[\Delta_{\alpha}\Big(\frac{\Theta}{\tau}\Big) - \frac{1}{2} \right] \,.
	$$
	Since by Lemma~\ref{lem:prop_Delta}, $\Delta_{\alpha}$ is even and non-increasing over $\R_{\geq 0}$, we have
	\begin{align*}
		\E\left[\Delta_{\alpha}\Big(\frac{\Theta}{\tau}\Big) - \frac{1}{2} \right]
		&\geq s \Delta_{\alpha}(+\infty) + (1-s) \Delta_{\alpha}(0) -\frac{1}{2}
		\\
		&= - s \frac{\alpha^2}{2} + (1-s) \Big(\frac{1}{2} + \alpha \phi(\alpha) - (1+\alpha^2) \Phi(-\alpha)\Big) - \frac{1}{2}
		= -\frac{1}{2} M_s(\alpha) \,.
	\end{align*}
\end{proof}

\begin{lemma}
	Let $s < s_{\rm max}(\delta)$ and assume that $\P(\Theta \neq 0) \leq s$. Then the following inequalities hold
	\begin{align}
		&\beta_* \tau_*(\beta_*) (\delta - M_s(\alpha_*)) \leq \delta \big( \sigma^2 + \beta_*^2\big) \,, \label{eq:ineq_tau_1}
		\\
		&- \tau_*(\beta_*) \lambda M_s'(\alpha_*) \leq \delta \sigma^2 \,,
		\label{eq:ineq_tau_2}
		\\
		&\tau_*(\beta_*) (\delta - \frac{1}{2}M_s''(\alpha_*)) \leq \beta_* \,. \label{eq:ineq_tau_3}
	\end{align}
\end{lemma}
\begin{proof}
	The inequality~\eqref{eq:ineq_tau_1} simply follows from the previous lemma and from the fact that
	$$
	\psi_{\lambda}(\beta,\tau_*(\beta)) \leq \max_{\beta \geq 0} \min_{\tau \geq \sigma} \psi_{\lambda}(\beta,\tau) \leq \frac{\sigma^2}{2} \,,
	$$
	by Lemma~\ref{lem:stupid_bound}.
	Let us prove~\eqref{eq:ineq_tau_2}. By Lemma~\ref{lem:intro_H}, we have
	$$
	\psi_{\lambda}(\beta_*,\tau_*(\beta_*)) \geq \frac{\lambda \tau_*(\beta_*)}{\delta} \E \left[ H_{\alpha_*} \! \Big(\frac{\Theta}{\tau_*(\beta_*)}\Big)\right] \,.
	$$
	Since by Lemma~\ref{lem:prop_H}, $H_{\alpha_*}$ is even, decreasing on $\R_{\geq 0}$, we have
	\begin{align*}
		\E \left[ H_{\alpha_*} \! \Big(\frac{\Theta}{\tau_*(\beta_*)}\Big)\right]
		&\geq s H_{\alpha_*}(+\infty) + (1-s) H_{\alpha_*}(0) 
		\\
		&= - s \alpha_* + 2(1-s) \big(\phi(\alpha_*) - \alpha_* \Phi(-\alpha_*)\big)
		= - \frac{1}{2} M_s'(\alpha_*) \,,
	\end{align*}
	which proves~\eqref{eq:ineq_tau_2}. To prove~\eqref{eq:ineq_tau_3} we use the optimality condition at $\beta_*$:
	\begin{equation} \label{eq:opti_beta}
		0 = \Psi_{\lambda}'(\beta_*) = \tau_*(\alpha_*) \left(1 - \frac{1}{\delta} \E \left[\Phi\Big(\frac{\Theta}{\tau_*(\alpha_*)} - \alpha_*\Big) + \Phi\Big(-\frac{\Theta}{\tau_*(\alpha_*)} - \alpha_*\Big)\right]\right) - \beta_* \,.
	\end{equation}
	The function $x \mapsto \Phi(x-\alpha_*) + \Phi(-x - \alpha_*)$ is even, increasing on $\R_{\geq 0}$. Therefore
	$$
	\E \left[\Phi\Big(\frac{\Theta}{\tau_*(\alpha_*)} - \alpha_*\Big) + \Phi\Big(-\frac{\Theta}{\tau_*(\alpha_*)} - \alpha_*\Big)\right] \leq s + 2 (1-s) \Phi(-\alpha_*)
	= \frac{1}{2} M_{s}''(\alpha_*) \,.
	$$
	Combining this inequality with~\eqref{eq:opti_beta} leads to~\eqref{eq:ineq_tau_3}.
\end{proof}


\begin{proposition}
	Let $s < s_{\rm max}(\delta)$ and assume that $\P(\Theta \neq 0) \leq s$.
	Then, there exists a constant $\tau_{\rm max} = \tau_{\rm max}(\delta,\lambda_{\rm min},\lambda_{\rm max},s,\sigma)$ such that for all $\lambda \in [\lambda_{\rm min},\lambda_{\rm max}]$,
	$$
	\sigma \leq \tau_*(\beta_*(\lambda)) \leq \tau_{\rm max} \,.
	$$
\end{proposition}
\begin{proof}
	Let $(\beta_*,\tau_*)$ be the unique optimal couple and recall $\alpha_* = \lambda / \beta_*$.
	We distinguish 3 cases:
	\\

	\textbf{Case 1}: $\alpha_* \geq \alpha_0$. In that case $\frac{1}{2} M_s''(\alpha_*) \leq M_s(\alpha_0) < \delta$. The inequality~\eqref{eq:ineq_tau_3} gives
	$$
	\tau_*(\beta_*)(\delta - \frac{1}{2} M_s''(\alpha_*)) \leq \beta_* \leq \beta_{\rm max} \,,
	$$
	which gives $\tau_*(\beta_*) \leq \frac{\beta_{\rm max}}{\delta - M_s(\alpha_0)}$.
	\\

	\textbf{Case 2}: $\alpha_* \in [(\alpha_1 + \alpha_0)/2, \alpha_0]$. In that case $\delta - M_s(\alpha_*) \geq c >0$, for some constant $c=c(\delta,s)>0$. Now, by~\eqref{eq:ineq_tau_1}
	$$
	\beta_*\tau_*(\beta_*) (\delta - M_s(\alpha_*)) \leq \delta (\sigma^2 + \beta_*^2) \,.
	$$
	Therefore,
	$$
	\tau_* \leq
	\frac{\delta}{c \beta_{\rm min}} (\sigma^2 + \beta_{\rm max}^2) \,.
	$$
	\\

	\textbf{Case 3}: $\alpha_* < (\alpha_1 + \alpha_0) /2$. In that case  $M_s'(\alpha_*) \leq - c$, for some constant $c=c(\delta,s) >0$. Consequently by~\eqref{eq:ineq_tau_2} we get
	$$
	\tau_*(\beta_*) \leq \frac{\sigma^2 \delta}{\lambda c} \,.
	$$
\end{proof}

\subsection{Dependency in \texorpdfstring{$\lambda$}{lambda}}

\begin{proposition}\label{prop:dep_lambda}
	\begin{itemize}
		\item The mapping $\lambda \mapsto \beta_*(\lambda)$ is $\cC^{\infty}$ and $2 \alpha_{\rm min}^{-1}$-Lipschitz on $\R_{>0}$.
		\item The mapping $\lambda \mapsto \tau_*(\lambda)$ is $\cC^{\infty}$ and $M$-Lipschitz on $[\lambda_{\rm min},\lambda_{\rm max}]$, for some constant $M(\Omega)>0$.
	\end{itemize}
\end{proposition}
\begin{proof}
	The first point has already been by Proposition~\ref{prop:beta_lip}.
	$\lambda \mapsto \tau_*(\lambda)$ is the composition of the mappings $\lambda \mapsto \alpha_*(\lambda)$ and $\alpha \mapsto \tau_*(\alpha)$, that are both $\cC^{\infty}$ by Lemma~\ref{lem:tau_star} and Proposition~\ref{prop:beta_lip}. Compute the derivative:
	\begin{align*}
		\frac{\partial \tau_*}{\partial \lambda}(\lambda)
		=
		\frac{\partial \alpha_*}{\partial \lambda}(\lambda)
		\frac{\partial \tau_*}{\partial \alpha}(\alpha_*(\lambda)) \,.
	\end{align*}
	Recall that $\alpha_*(\lambda) = \lambda / \beta_*(\lambda)$. Thus
	$$
	\left|\frac{\partial \alpha_*}{\partial \lambda}(\lambda)\right|
	=
	\left|\frac{1}{\beta_*(\lambda)}
	-
	\frac{\partial \beta_*}{\partial \lambda}(\lambda)
	\frac{\lambda}{\beta_*(\lambda)^2} \right|
	\leq \beta_{\rm min}^{-1} + 2 \alpha_{\rm min}^{-1} \lambda_{\rm max} \beta_{\rm  min}^{-2} \,.
	$$
	By Lemma~\ref{lem:tau_star}, we have
	$$
	\left|\frac{\partial \tau_*}{\partial \alpha}(\alpha_*(\lambda)) \right|
	\leq (\alpha_*(\lambda) + 1) \frac{\tau_*(\alpha_*)^3}{\delta \sigma^2} \,.
	$$
	Since by Theorem~\ref{th:control_beta_tau}, $\tau_{*} (\alpha_*) \leq \tau_{\rm max}(\Omega)$ and $\alpha_{*} \leq \lambda_{\rm max} / \beta_{\rm min}(\Omega)$, the derivative of $\tau_*$ with respect to $\lambda$ is bounded on $[\lambda_{\rm min},\lambda_{\rm max}]$.
\end{proof}

\section{Study of Gordon's optimization problem for \texorpdfstring{$\what{w}_{\lambda}$}{w}}\label{sec:gordon_w}

In this section we study $L_{\lambda}$ defined by~\eqref{eq:def_L}.
Define, for $w \in \R^N$ and $\beta \geq 0$
\begin{equation}\label{eq:def_l}
\ell_{\lambda}(w,\beta) = 
\left(\sqrt{\frac{\|w\|^2}{n} + \sigma^2}\ \frac{\|h\|}{\sqrt{n}}
-	\frac{1}{n} g^{\sT} w + \frac{g'\sigma}{\sqrt{n}} \right)\beta
- \frac{1}{2} \beta^2
+ \frac{\lambda}{n} |w+\theta^{\star}| 
- \frac{\lambda}{n} |\theta^{\star}|  \,.
\end{equation}
So that $L_{\lambda}(w) = \max_{\beta \geq 0} \ell_{\lambda}(w,\beta)$.
Let us define the vector $\mathsf{w}_{\lambda} \in \R^N$ by
\begin{equation}\label{eq:def_w_sf}
\mathsf{w}_{\lambda,i} = 
\eta\left(
	\theta^{\star}_i + \tau_*(\lambda) g_i, \frac{\tau_*(\lambda)\lambda}{\beta_*(\lambda)}
\right) - \theta^{\star}_i \,.
\end{equation}
The goal of this section is to prove that, with high probability, the minimizer of $L_{\lambda}$ is close to $\mathsf{w}_{\lambda}$ and that $L_{\lambda}$ is strongly convex around $\mathsf{w}_{\lambda}$.
\begin{proposition}\label{prop:w_star}
	$L_{\lambda}$ admits almost surely a unique minimizer $w_{\lambda}^*$ on $\R^N$.
\end{proposition}
\begin{proof}
	$L_{\lambda}$ is a convex function that goes to $+\infty$ at infinity, so it admits minimizers over $\R^N$. 
	\\

	\noindent\textbf{Case 1}: there exists a minimizer $w$ such that 
	$ \sqrt{\frac{\|w\|^2}{n} + \sigma^2}\ \frac{\|h\|}{\sqrt{n}} -	\frac{1}{n} g^{\sT} w + \frac{g'\sigma}{\sqrt{n}} >0$.
	\\
	In that case, there exist a neighborhood $O_w$ of $w$ such that for all $w' \in O_w$
	$$
	a(w'):= \sqrt{\frac{\|w'\|^2}{n} + \sigma^2}\ \frac{\|h\|}{\sqrt{n}} -	\frac{1}{n} g^{\sT} w' + \frac{g'\sigma}{\sqrt{n}} >0 \,.
	$$
	Thus for all $w' \in O_w$, $L_{\lambda}(w') = \frac{1}{2} a(w')^2 + \frac{\lambda}{n} |w'+\theta^{\star}| - \frac{\lambda}{n} |\theta^{\star}|$. 
	Recall that the composition of a strictly convex function and a strictly increasing function is strictly convex.
	$L_{\lambda}$ is therefore strictly convex on $O_w$ because $a$ is strictly convex and remains strictly positive on $O_w$ and because $x > 0 \mapsto x^2$ is strictly increasing.
	$w$ is thus the only minimizer of $L_{\lambda}$.
	\\

	\noindent\textbf{Case 2}: for all minimizer $w$ we have
	$ \sqrt{\frac{\|w\|^2}{n} + \sigma^2}\ \frac{\|h\|}{\sqrt{n}} -	\frac{1}{n} g^{\sT} w + \frac{g'\sigma}{\sqrt{n}} \leq 0$.
	\\
	Let $w$ be a minimizer of $L_{\lambda}$. 
	The optimality condition gives
	$$
	- \lambda^{-1}
	\left(\sqrt{\frac{\|w\|^2}{n} + \sigma^2}\ \frac{\|h\|}{\sqrt{n}}
	-	\frac{1}{n} g^{\sT} w + \frac{g'\sigma}{\sqrt{n}} \right)_{+}
	\left(\frac{w}{\sqrt{\frac{\|w\|^2}{n} + \sigma^2}}\ \frac{\|h\|}{\sqrt{n}}
	-	 g \right)
	\in \partial |\theta^{\star} + w|.
	$$
	We obtain then $0 \in \partial |\theta^{\star} + w|$ which implies $w=-\theta^{\star}$: $L_{\lambda}$ has a unique minimizer.
\end{proof}

\subsection{Local stability of Gordon's optimization}

\begin{theorem}\label{th:gordon_aux}
	There exists constants $\gamma,c,C > 0$ that only depend on $\Omega$ such that for all $\theta^{\star} \in \mathcal{D}$, all $\lambda \in [\lambda_{\rm min}, \lambda_{\rm max}]$ and all $\epsilon \in (0,1]$
	$$
	\P \Big(
		\exists w\in\R^N, \quad 
		\frac{1}{N} \| w - \mathsf{w}_{\lambda} \|^2 > \epsilon
		\quad \text{and} \quad
		L_{\lambda}(w) \leq \min\limits_{v \in \R^N} L_{\lambda}(v) + \gamma \epsilon 
	\Big)
	\leq \frac{C}{\epsilon} e^{-cn\epsilon^2} \,.
	$$
\end{theorem}

We deduce from Theorem~\ref{th:gordon_aux} that for all $\epsilon \in (0,1]$ with probability at least $1 - C \epsilon^{-1} e^{-cn\epsilon^2}$, $\frac{1}{N} \| w^*_{\lambda} - \mathsf{w}_{\lambda} \|^2 \leq \epsilon$. From this we deduce easily that with the same probability $|L_{\lambda}(w^*_{\lambda}) - L_{\lambda}(\mathsf{w}_{\lambda})| \leq M \epsilon$, for some constant $M>0$, which gives by Proposition~\ref{prop:concentration_cost_w_star}:

\begin{corollary}\label{cor:concentration_gordon}
	Define
	\begin{equation}\label{eq:def_L_0}
		L_*(\lambda) = \psi_{\lambda}(\beta_*(\lambda),\tau_*(\lambda)) \,.
	\end{equation}
	The exists constants $c,C>0$ that only depend on $\Omega$ such that 
	$$
	\P \left(\left| \min_{w \in \R^N} L_{\lambda}(w) - L_*(\lambda) \right| \geq \epsilon \right) \leq \frac{C}{\epsilon} e^{-nc\epsilon^2} \,.
	$$
\end{corollary}

\subsection{Proof of Theorem~\ref{th:gordon_aux}}

\begin{proposition}\label{prop:min_ball}
	For all $R > 0$ there exists constants $c,C>0$ that only depend on $(\Omega,R)$, such that
	for all $\epsilon \in (0,1]$,
	$$
	\forall \theta^{\star} \in \mathcal{D}, \ \forall \lambda \in [\lambda_{\rm min},\lambda_{\rm max}], \quad
	\P \Big( 
		L_{\lambda}(\mathsf{w}_{\lambda}) \leq \min_{\|w\| \leq \sqrt{n} R} L_{\lambda}(w) + \epsilon
	\Big) \geq 1 - \frac{C}{\epsilon} e^{-cn\epsilon^2}  \,.
	$$
\end{proposition}
\begin{proof}
	Notices that it suffices to proves the proposition for $\epsilon$ smaller than some constant.
	Let $\theta^{\star} \in \mathcal{D}$, $\lambda \in [\lambda_{\rm min},\lambda_{\rm max}]$.
	Let $R > 0$
	and $\epsilon \in \big(0,\min(1,\sigma^2/2)\big]$. 
	Define 
	$$
	\ell^{\circ}_{\lambda}(w,\beta) = 
	\left(\sqrt{\frac{\|w\|^2}{n} + \sigma^2} 
	-	\frac{1}{n} g^{\sT} w  \right)\beta
	- \frac{1}{2} \beta^2
	+ \frac{\lambda}{n} |w+\theta^{\star}|
	- \frac{\lambda}{n} |\theta^{\star}| \,.
	$$
	On the event 
	\begin{equation} \label{eq:event_gh}
		\left\{  \left|\frac{1}{n} \|h\|^2 - 1\right| \leq \epsilon \right\}
		\bigcap
		\left\{  \left|\frac{g'\sigma}{\sqrt{n}} \right| \leq \epsilon \right\} \,,
	\end{equation}
	which has probability at least $1- C e^{-c n\epsilon^2}$, 
	we have, for all $w \in B(0,R \sqrt{n})$ and $\beta \in [0,\beta_{\rm max}]$:
	\begin{align*}
	\left|
	\ell_{\lambda}(w,\beta) - \ell_{\lambda}^{\circ}(w,\beta)
	\right| 
	&= 
	\beta \sqrt{\frac{\|w\|^2}{n} + \sigma^2}\ \left|\frac{\|h\|}{\sqrt{n}} - 1 \right|
	+
	\beta \left|\frac{g'\sigma}{\sqrt{n}} \right|
	\\
	&\leq \beta_{\rm max} \sqrt{\sigma^2 + R^2} \left| \frac{1}{n} \| h\|^2 -1 \right| 
	+ \beta_{\rm max} \epsilon
	\leq \underbrace{\beta_{\rm max} (\sqrt{\sigma^2 + R^2}+ 1)}_{K} \epsilon \,.
	\end{align*}
	For simplicity we write $(\beta_*,\tau_*) = (\beta_*(\lambda),\tau_*(\lambda))$.
	We have on the event~\eqref{eq:event_gh}:
	\begin{align*}
		\min_{\|w\| \leq R \sqrt{n}} L_{\lambda}(w)
		&=
		\min_{\|w\| \leq R \sqrt{n}} \max_{\beta \geq 0} \ell_{\lambda}(w,\beta)
		\geq
		\min_{\|w\| \leq R \sqrt{n}} \ell_{\lambda}(w,\beta_*)
		\geq
		\min_{\|w\| \leq R \sqrt{n}} \ell_{\lambda}^{\circ}(w,\beta_*) - K \epsilon \,.
	\end{align*}
	Using the fact that for $w\in B(0,R\sqrt{n})$
	$$
	\sqrt{\frac{\|w\|^2}{n} + \sigma^2} = \min_{\sigma \leq \tau \leq \sqrt{\sigma^2 + R^2}} \left\{
		\frac{\frac{\|w\|^2}{n} + \sigma^2}{2\tau} + \frac{\tau}{2}
	\right\} \,,
	$$
	we obtain that
	\begin{align*}
		\min_{\|w\|\leq R \sqrt{n}} 
		\ell_{\lambda}^{\circ}(w,\beta_*) = 
		\min_{\sigma \leq \tau \leq \sqrt{\sigma^2 + R^2}} 
		\left\{
		\frac{\beta_*}{2} \left(  \frac{\sigma^2}{\tau} + \tau \right)- \frac{\beta_*^2}{2}
		+ \frac{1}{n} \min_{\|w\| \leq R \sqrt{n}} 
		\left\{
			\frac{\beta_*}{2\tau} \|w\|^2	- \beta_* g^{\sT} w  + \lambda |w + \theta^{\star}| - \lambda |\theta^{\star}|
		\right\}
	\right\} \,.
	\end{align*}
	For all $\tau \in [\sigma, \sqrt{\sigma^2 + R^2}]$
	the function
	$$
	g \mapsto
	\min_{\|w\| \leq R \sqrt{n}} 
	\left\{
		\frac{\beta_*}{2\tau} \|w\|^2	- \beta_* g^{\sT} w  + \lambda |w + \theta^{\star}| - \lambda |\theta^{\star}|
	\right\}
	$$
	is $\beta_{\rm max}R\sqrt{n}$-Lipschitz. 
	Therefore
	$$
	\mathsf{F}(\tau,g) = 
	\frac{\beta_*}{2} \left(  \frac{\sigma^2}{\tau} + \tau \right)- \frac{\beta_*^2}{2}
	+ \frac{1}{n} 
	\min_{\|w\| \leq R \sqrt{n}} 
	\left\{
		\frac{\beta_*}{2\tau} \|w\|^2	- \beta_* g^{\sT} w  + \lambda |w + \theta^{\star}| - \lambda |\theta^{\star}|
	\right\}
	$$
	is $\beta_{\rm max}^2 R^2 n^{-1}$-sub-Gaussian. Therefore there exists constants $C,c>0$ such that for all $\tau \in [\sigma,\sqrt{\sigma^2 + R^2}]$, we have
	$$
	\P \left( 
		\left| \mathsf{F}(\tau,g) - \E \mathsf{F}(\tau,g) \right| > \epsilon
	\right)
	\leq C e^{-cn\epsilon^2}  \,.
	$$
	$\mathsf{F}(\cdot,g)$ is almost-surely a $\beta_{\rm max}(1+\frac{R^2}{\sigma^2})$-Lipschitz function on $[\sigma,\sqrt{\sigma^2 + R^2}]$.
	Therefore, by an $\epsilon$-net argument one can find constants $C,c>0$ that only depend on $(\Omega,R)$, such that for all $\epsilon >0$ the event
	\begin{equation}\label{eq:event_uniform_psi}
		\left\{ \sup_{\tau \in [\sigma,\sqrt{\sigma^2 + R^2}]} 
			\left| \mathsf{F}(\tau,g) - \E \mathsf{F}(\tau,g) \right| \leq \epsilon
		\right\}
	\end{equation}
	has probability at least $1 - \frac{C}{\epsilon} e^{-cn\epsilon^2} $.
	On the event~\eqref{eq:event_uniform_psi} we have then $\min\limits_{\tau \in [\sigma,\sqrt{\sigma^2 + R^2}]} \mathsf{F}(\tau,g) \geq  \min\limits_{\tau \in [\sigma,\sqrt{\sigma^2 + R^2}]} \E \left[\mathsf{F}(\tau,g)\right] - \epsilon$.
	Notice that for all $\tau  >0$ we have
	\begin{align*}
		\frac{1}{n} 
		\E \left[ 
			\min_{\|w\| \leq R \sqrt{n}} 
			\left\{
				\frac{\beta_*}{2\tau} \|w\|^2	- \beta_* g^{\sT} w  + \lambda |w + \theta^{\star}| - \lambda |\theta^{\star}|
		\right\}\right]
		&\geq
		\frac{1}{n} \sum_{i=1}^N
		\E \left[ 
			\min_{w_i \in \R} 
			\left\{
				\frac{\beta_*}{2\tau} w_i^2	- \beta_* g_i w_i  + \lambda |w_i + \theta^{\star}_i| - \lambda |\theta_i^{\star}|
		\right\}\right]
		\\
		&=
		\frac{1}{\delta}
		\E \left[ 
			\min_{w \in \R} 
			\left\{
				\frac{\beta_*}{2\tau} w^2	- \beta_* Z w  + \lambda |w + \Theta| - \lambda |\Theta|
		\right\}\right],
	\end{align*}
	where the last expectation is with respect $(\Theta,Z) \sim \what{\mu}_{\theta^{\star}} \otimes \cN(0,1)$.
	Consequently on the event~\eqref{eq:event_gh} and~\eqref{eq:event_uniform_psi}, we have
	\begin{align*}
		(1 + K)\epsilon + \min_{\|w\| \leq R \sqrt{n}} L_{\lambda}(w)
		&\geq 
		\min_{\sigma \leq \tau \leq \sqrt{\sigma^2 + R^2}} 
		\psi_{\lambda}(\beta_*,\tau) 
		\geq
		\min_{\sigma \leq \tau} \psi_{\lambda}(\beta_*,\tau) 
		= \psi_{\lambda}(\beta_*,\tau_*) \,.
	\end{align*}
	By Proposition~\ref{prop:concentration_cost_w_star} we have that 
	$$
	\psi_{\lambda}(\beta_*,\tau_*) \geq L_{\lambda}(\mathsf{w}_{\lambda}) - \epsilon \,,
	$$
	with probability at least $1-C e^{-cn\epsilon^2}$.
	Then, for all $\epsilon \in (0,1)$ we have with probability at least $1 - \frac{C}{\epsilon} e^{-cn \epsilon^2}$
	$$
	\min_{\|w\| \leq R \sqrt{n}} L_{\lambda}(w) + (K + 2) \epsilon \geq L_{\lambda}(\mathsf{w}_{\lambda}) \,.
	$$
\end{proof}

\begin{lemma}\label{lem:convex_ball}
	Let $f$ be a convex function on $\R^N$. Let $w\in \R^N$ and $r>0$. Suppose that $f$ is $\gamma$-strongly convex on the ball $B(w,r)$, for some $\gamma>0$. Assume that
	$$
	f(w) \leq \min_{x \in B(w,r)} f(x) + \epsilon \,,
	$$
	for some $\epsilon \leq \frac{r^2 \gamma}{8}$.
	Then $f$ admits a unique minimizer $x^*$ over $\R^N$. We have $x^* \in B(w,r)$ and therefore
	$$
	\| x^* - w \|^2 \leq \frac{2}{\gamma} \epsilon \,.
	$$
	Moreover, for every $x \in \R^N$ such that $f(x) \leq \min f + \epsilon$ we have $\| x - w \|^2 \leq \frac{8}{\gamma} \epsilon$.
\end{lemma}
\begin{proof}
	$f$ is convex on $B(w,r)$, it admits therefore a minimizer $x^*$ on $B(w,r)$.
	By strong convexity we have
	$$
	\|x^*-w\|^2 \leq \frac{2}{\gamma} \epsilon \leq \frac{r^2}{4} \,.
	$$
	Consequently, $x^*$ is in the interior of $B(w,r)$. By strong convexity, $x^*$ is then the unique minimizer of $f$ over $\R^N$.
	By strong convexity, for any $x$ outside of $B(w,r)$ we have
	$$
	f(x) > f(x^*) + \frac{1}{2} \gamma \left(\frac{r}{2}\right)^2 \geq f(x^*) + \epsilon \,.
	$$
	Consequently, if $f(x) \leq \min f + \epsilon$ then $x \in B(w,r)$ and thus $\|x-x^*\|^2 \leq \frac{2}{\gamma} \epsilon$.
\end{proof}
\\

\begin{proof}[of Theorem~\ref{th:gordon_aux}]
	Let $t = \min(\frac{1}{16}\beta_{\rm min}, \sigma)$. By Lemma~\ref{lem:conc_w_RS} the event
	\begin{equation} \label{eq:event_beta}
		\left\{
			\left| \frac{\|\mathsf{w}_{\lambda}\|^2}{n} -  \frac{\E \|\mathsf{w}_{\lambda}\|^2}{n} \right| \leq t^2
			,
			\quad
			\frac{g^{\sT} \mathsf{w}_{\lambda}}{n} \leq \E \left[\frac{g^{\sT} \mathsf{w}_{\lambda}}{n}\right] + t
			,
			\quad
			\|g\| \leq 2 \sqrt{N},
			\quad
			\Big|\frac{\sigma g'}{\sqrt{n}} \Big| \leq \frac{\beta_{\rm min}}{4},
			\quad
			(1-\frac{\beta_{\rm min}}{8 \tau_{\rm max}}) \leq \frac{\|h\|}{\sqrt{n}} \leq 2
		\right\}
	\end{equation}
	has probability at least $1 - C e^{-cn}$, for some constants $C,c >0$. On the event~\eqref{eq:event_beta}
	\begin{align*}
		\sqrt{\frac{\|\mathsf{w}_{\lambda}\|^2}{n} + \sigma^2} 
		&\geq
		\sqrt{\frac{\E \|\mathsf{w}_{\lambda}\|^2}{n} + \sigma^2 - t^2} 
		\geq
		\tau_* - t \,.
	\end{align*}
	Therefore
	\begin{align*}
		\sqrt{\frac{\|\mathsf{w}_{\lambda}\|^2}{n} + \sigma^2} \frac{\|h\|}{\sqrt{n}}
		&\geq
		\frac{\|h\|}{\sqrt{n}}\tau_* - \frac{\|h\|}{\sqrt{n}}t
		\geq
		\tau_* - \frac{\beta_{\rm min}}{8} - 2 t \geq \tau_* - \frac{\beta_{\rm min}}{4} \,.
	\end{align*}
	Consequently, on the event~\eqref{eq:event_beta} we have
	\begin{align*}
		\sqrt{\frac{\|\mathsf{w}_{\lambda}\|^2}{n} + \sigma^2}\frac{\|h\|}{\sqrt{n}} - \frac{1}{n} g^{\sT} \mathsf{w}_{\lambda} +  \frac{\sigma g'}{\sqrt{n}}
		&\geq
		\tau_* - \frac{1}{n} \E \left[ g^{\sT} \mathsf{w}_{\lambda}\right] - \frac{3}{4} \beta_{\rm min}
		\geq \frac{1}{4} \beta_{\rm min} \,,
	\end{align*}
	because $\tau_* - \frac{1}{n} \E[g^{\sT} \mathsf{w}_{\lambda}] = \beta_* \geq \beta_{\rm min}$.
	Moreover, on the event~\eqref{eq:event_beta} the function
	$$
	f: w \mapsto \sqrt{\frac{\|w\|^2}{n} + \sigma^2}\frac{\|h\|}{\sqrt{n}} - \frac{1}{n} g^{\sT} w + \frac{g' \sigma}{\sqrt{n}}
	$$
	is $\frac{2\sqrt{N}}{n} + \frac{2}{\sqrt{n}}$-Lipschitz. We have seen above that on~\eqref{eq:event_beta}, $f(\mathsf{w}_{\lambda}) \geq \frac{1}{4} \beta_{\rm min}$.
	Thus we can find a constant $r>0$ such that on the event~\eqref{eq:event_beta} we have for all $w \in B(\mathsf{w}_{\lambda}, r \sqrt{n})$
	$$
	f(w)> \frac{1}{8} \beta_{\rm min} \,.
	$$
	By Lemma~\ref{lem:convex_sqrt}, the function $f$ is $\frac{a}{n}$-strongly convex on $B(\mathsf{w}_{\lambda},r \sqrt{n})$, for some constant $a>0$. For all $w \in B(\mathsf{w}_{\lambda},r \sqrt{n})$ we have
	$$
	L_{\lambda}(w) = \frac{1}{2} f(w)^2 + \frac{\lambda}{n} \left(|w+\theta^{\star}| - |\theta^{\star}| \right) \,.
	$$
	Compute the Hessian for $w \in B(\mathsf{w}_{\lambda},r\sqrt{n})$:
	$$
\nabla^2\Big(\frac{1}{2}f^2\Big)(w) = f(w) \nabla^2f (w) + \nabla f (w) \nabla f(w)^{\sT}
	\succeq \frac{a\beta_{\rm min}}{8n} \bbf{I}_N \,,
	$$
	which means that $L$ is $\frac{\gamma}{n}$-strongly convex on $B(\mathsf{w}_{\lambda},r\sqrt{n})$, for some constant $\gamma>0$.
	\\

	Notice that it suffices to prove Theorem~\ref{th:gordon_aux} for $\epsilon \in (0,q]$ for some constant $q >0$.
	Let $\epsilon \in (0,\frac{\gamma r^2}{8})$.
	Let now apply Proposition~\ref{prop:min_ball} with $R = \tau_{\rm max} + r$: with probability at least $1- \frac{C}{\epsilon} e^{-c n \epsilon^2}$
	\begin{equation}\label{eq:min_w_big_ball}
		\Big\{ L_{\lambda}(\mathsf{w}_{\lambda}) \leq \min_{\|w\| \leq R \sqrt{n}} L_{\lambda}(w) + \epsilon \Big\} \,.
	\end{equation}
	Notice that on the event~\eqref{eq:event_beta},
	$$
	\frac{\|\mathsf{w}_{\lambda}\|^2}{n} 
	\leq \E \left[\frac{\|\mathsf{w}_{\lambda}\|^2}{n}\right] + t^2
	\leq \E \left[\frac{\|\mathsf{w}_{\lambda}\|^2}{n}\right] + \sigma^2 = \tau_*^2 \leq \tau_{\rm max}^2 \,.
	$$
	Therefore, on~\eqref{eq:event_beta}, $B(\mathsf{w}_{\lambda},r \sqrt{n}) \subset B(0,R\sqrt{n})$. Using then~\eqref{eq:min_w_big_ball} we get
	$$
	L_{\lambda}(\mathsf{w}_{\lambda}) \leq \min_{w \in B(\mathsf{w}_{\lambda},r\sqrt{n})} L_{\lambda}(w) + \epsilon \,.
	$$
	Consequently, on the events~\eqref{eq:event_beta} and~\eqref{eq:min_w_big_ball} that have probability at least $1 - \frac{C}{\epsilon}e^{-cn\epsilon^2}$, Lemma~\ref{lem:convex_ball} gives that for all $w \in \R^N$ such that $L_{\lambda}(w) \leq \min\limits_{v \in \R^N} L_{\lambda}(v) + \epsilon$ we have $ \|\mathsf{w}_{\lambda} - w \|^2 \leq \frac{8n}{\gamma} \epsilon$.
\end{proof}

\section{Empirical distribution and risk of the Lasso}

\subsection{Proofs of local stability of the Lasso cost} \label{sec:proof_key}

\subsubsection{Application of Gordon's min-max Theorem}
\begin{proposition} \label{prop:apply_gordon}
	There exists constants $c,C >0$ that only depend on $\Omega$ such that
	for all closed set $D \subset \R^N$ and for all $\epsilon \in (0,1]$,
	$$
	\P \left(\min_{w \in D} \cC_{\lambda}(w) \leq \min_{w \in \R^N} \cC_{\lambda}(w) + \epsilon \right) 
	\leq
	2 \P \left(\min_{w \in D} L_{\lambda}(w) \leq \min_{w \in \R^N} L_{\lambda}(w) + 3 \epsilon \right)  + \frac{C}{\epsilon} e^{-cn \epsilon^2} \,.
	$$
\end{proposition}

In order to prove this, we start by showing that the optimal Lasso cost concentrates around $L_*(\lambda)$.
\begin{proposition} \label{prop:concentration_min_c}
	There exists constants $c,C >0$ that only depend on $\Omega$ such that
	for all closed set $D \subset \R^N$ and for all $\epsilon \in (0,1]$,
	$$
	\P \left(\Big|\min_{w \in \R^N} \cC_{\lambda}(w) - L_*(\lambda) \Big| \geq \epsilon \right) 
	\leq \frac{C}{\epsilon} e^{-cn\epsilon^2} \,.
	$$
\end{proposition}
\begin{proof}
	By Corollary~\ref{cor:gordon}, we have
	$$
	\P \left(\min_{w \in \R^N} \cC_{\lambda}(w) - L_*(\lambda) \geq \epsilon \right) 
	\leq
	2 \P \left(\min_{w \in \R^N} L_{\lambda}(w) - L_*(\lambda) \geq \epsilon \right) 
	\leq \frac{C}{\epsilon} e^{-cn\epsilon^2} \,,
	$$
	where the last inequality comes from Corollary~\ref{cor:concentration_gordon}. The bound of the probability of the converse inequality is proved analogously.
\end{proof}
\\

\begin{proof}[of Proposition~\ref{prop:apply_gordon}]
	Recall that $L_*(\lambda)$ is defined by~\eqref{eq:def_L_0}. Let $D \subset \R^N$ be a closed set.
	\begin{align*}
		&\P\left( \min_{w \in D} \cC_{\lambda}(w) \leq \min_{w \in \R^N} \cC_{\lambda}(w) + \epsilon \right)
		\\
		&\leq
		\P\left(
			\min_{w \in D} \cC_{\lambda}(w)
			\leq
			\min_{w \in \R^N} \cC_{\lambda}(w) + \epsilon
			\quad \text{and} \quad
			\min_{w \in \R^N} \cC_{\lambda}(w) \leq L_*(\lambda) + \epsilon
		\right)
		+
		\P\left(
			\min_{w \in \R^N} \cC_{\lambda}(w) > L_*(\lambda) + \epsilon
		\right)
		\\
		&\leq
		\P\left(
			\min_{w \in D} \cC_{\lambda}(w)
			\leq L_*(\lambda) + 2 \epsilon
		\right)
		+
		\frac{C}{\epsilon} e^{-cn\epsilon^2} \,,
	\end{align*}
	where we used Proposition~\ref{prop:concentration_min_c} above.
	We can now apply the first point of Corollary~\ref{cor:gordon} to obtain:
	\begin{equation}\label{eq:apply_gordon_1}
		\P\left(
			\min_{w \in D} \cC_{\lambda}(w)
			\leq L_*(\lambda) + 2 \epsilon
		\right)
		\leq
		2 \P\left(
			\min_{w \in D} L_{\lambda}(w)
			\leq L_*(\lambda) + 2 \epsilon
		\right) \,.
	\end{equation}
	We thus get
	\begin{align*}
		&\P\left( \min_{w \in D} \cC_{\lambda}(w) \leq \min_{w \in \R^N} \cC_{\lambda}(w) + \epsilon \right)
		\leq
		2 \P\left(
			\min_{w \in D} L_{\lambda}(w)
			\leq L_*(\lambda) + 2 \epsilon
		\right)
		+
		\frac{C}{\epsilon} e^{-cn\epsilon^2}
		\\
		&\leq
		2 \P\left(
			\min_{w \in D} L_{\lambda}(w)
			\leq \min_{w\in \R^N}L_{\lambda}(w) + 3 \epsilon
		\right)
		+
		2 \P\left(
			\min_{w \in \R^N} L_{\lambda}(w) < L_*(\lambda) - \epsilon
		\right)
		+
		\frac{C}{\epsilon} e^{-cn\epsilon^2}
		\\
		&\leq
		2 \P\left(
			\min_{w \in D} L_{\lambda}(w)
			\leq \min_{w\in \R^N}L_{\lambda}(w) + 3 \epsilon
		\right)
		+
		\frac{C}{\epsilon} e^{-cn\epsilon^2} \,,
	\end{align*}
	for some constants $c,C > 0$, because of Corollary~\ref{cor:concentration_gordon}.
\end{proof}

\subsubsection{Local stability of the empirical distribution of the Lasso estimator: proof of Theorem~\ref{th:key_law}}

For $w \in \R^N$, let $\mu_0(w)$ be the probability distribution over $\R^2$ defined by
$$
\what{\mu}_0(w) = \frac{1}{N} \sum_{i=1}^N \delta_{(w_i+\theta^{\star}_i,\theta^{\star}_i)}
$$
Theorem~\ref{th:key_law} follows from Proposition~\ref{prop:apply_gordon} and the following Lemma.

\begin{lemma}\label{lem:strong_convex_W_L_law}
	Assume that $\mathcal{D} = \cF_{p}(\xi)$ for some $\xi,p>0$.
	There exists constants $\gamma,c,C > 0$ that depend only on $\Omega$, such that for all $\epsilon \in (0,\frac{1}{2}]$ we have
	$$
	\P\left(
		\min_{w \in D_{\epsilon}} L_{\lambda}(w)
		\leq \min_{w\in \R^N}L_{\lambda}(w) + 3 \gamma \epsilon
	\right)
	\leq 
	C \epsilon^{-\max(1,a)} \exp\left(-cN \epsilon^2 \epsilon^a \log(\epsilon)^{-2} \right) \,,
	$$
	where $ D_{\epsilon} = \left\{ w \in \R^N \, \middle| \, W_2(\what{\mu}_0(w),\mu_{\lambda}^*)^2 \geq \epsilon \right\}$ and $a = \frac{1}{2} + \frac{1}{p}$.
\end{lemma}
\begin{proof}
	By Theorem~\ref{th:gordon_aux} and Proposition~\ref{prop:concentration_empirical_distribution} there exists constants $\gamma,c,C>0$ such that for all $\epsilon \in (0,\frac{1}{2}]$ the event
	\begin{equation}\label{eq:event_good}
		\Big\{ \forall w \in \R^N, \ L_{\lambda}(w) \leq \min_{v \in \R^N} L_{\lambda}(v) + 3 \gamma \epsilon  \implies \frac{1}{N} \|w-\mathsf{w}_{\lambda}\|^2 \leq \frac{\epsilon}{5} \Big\}
		\bigcap
		\Big\{ W_2\big(\mu_{\lambda}^*,\what{\mu}_0(\mathsf{w}_{\lambda})\big)^2 \leq \frac{\epsilon}{4} \Big\}
	\end{equation}
	has probability at least 
	$$
	1 - C \epsilon^{-1} \exp\big( -cn\epsilon^2\big) - C \epsilon^{-a} \exp\left(-cN \epsilon^2 \epsilon^a \log(\epsilon)^{-2} \right)
	\geq
	1 - C \epsilon^{-\max(1,a)} \exp\left(-cN \epsilon^2 \epsilon^a \log(\epsilon)^{-2} \right) \,,
	$$
	where $a=\frac{1}{2} + \frac{1}{p}$.
	On the event~\eqref{eq:event_good}, we have for all $w \in D_{\epsilon}$:
	$$
	\frac{1}{N} \| w - \mathsf{w}_{\lambda} \|^2 
	\geq W_2\big(\what{\mu}_0(w),\what{\mu}_0(\mathsf{w}_{\lambda})\big)^2
	\geq 
	\Big(
		W_2(\what{\mu}_0(w),\mu_{\lambda}^*)
		-W_2(\mu_{\lambda}^*,\what{\mu}_0(\mathsf{w}_{\lambda}))
	\Big)^2
	\geq \frac{\epsilon}{4} \,.
	$$
	This gives that on the event~\eqref{eq:event_good}, for all $w \in D_{\epsilon}$, $L_{\lambda}(w) > \min\limits_{v \in \R^N} L_{\lambda}(v) + 3 \gamma \epsilon$. The intersection of~\eqref{eq:event_good} with the event $\big\{ \min\limits_{w \in D_{\epsilon} } L_{\lambda}(w) \leq \min\limits_{w\in \R^N}L_{\lambda}(w) + 3 \gamma \epsilon \big\}$ is therefore empty: the lemma is proved.
\end{proof}

\subsubsection{Local stability of the risk of the Lasso estimator}\label{sec:proof_risk_lasso}

We prove here the analog of Theorem~\ref{th:key_law} for the risk of the Lasso estimator.
\begin{theorem}\label{th:key_risk} for the risk of the Lasso estimator.
	There exists constants $C,c,\gamma>0$ that only depend on $\Omega$ such that
	for all $\epsilon \in (0,1]$
	$$
	\sup_{\lambda \in [\lambda_{\rm min},\lambda_{\rm max}]} \
	\sup_{\theta^{\star} \in \mathcal{D}} \
	\P\! \left(\exists \theta \in \R^N\!, \ 
		\Big( \frac{1}{N}\| \theta - \theta^{\star} \|^2 - R_*(\lambda) \Big)^2 \geq \epsilon
		\quad \text{and} \quad
	\mathcal{L}_{\lambda}(\theta) \leq \min \mathcal{L}_{\lambda} + \gamma \epsilon \right) \leq \frac{C}{\epsilon} e^{-cN\epsilon^2}.
	$$
\end{theorem}

Theorem~\ref{th:key_risk} follows from Proposition~\ref{prop:apply_gordon} and the following Lemma.

\begin{lemma}\label{lem:strong_convex_W_L_risk}
	There exists constants $\gamma,c,C > 0$ that only depend on $\Omega$ such that for all $\epsilon \in (0,1]$ we have
	$$
	\P\left(
		\min_{w \in D_{\epsilon}} L_{\lambda}(w)
		\leq \min_{w\in \R^N}L_{\lambda}(w) + 3\gamma \epsilon
	\right)
	\leq \frac{C}{\epsilon} e^{-cn\epsilon^2} \,,
	$$
	where $D_{\epsilon} = \left\{ w \in \R^N \, \middle| \left(\|w\| - \sqrt{N R_*(\lambda)}\right)^2  \geq N \epsilon \right\}$.
\end{lemma}
\begin{proof}
	By Theorem~\ref{th:gordon_aux} and Lemma~\ref{lem:conc_w_RS} there exists constants $\gamma,c,C>0$ such that for all $\epsilon \in (0,1]$ the event
	\begin{equation}\label{eq:event_good0}
		\Big\{ \forall w \in \R^N, \ L_{\lambda}(w) \leq \min_{v \in \R^N} L_{\lambda}(v) + 3 \gamma \epsilon  \implies \frac{1}{N} \|w-\mathsf{w}_{\lambda}\|^2 \leq \frac{\epsilon}{5} \Big\}
		\bigcap
		\Big\{ \big( \|\mathsf{w}_{\lambda}\| - \sqrt{N R_*(\lambda)} \big)^2 \leq N \frac{\epsilon}{4} \Big\}
	\end{equation}
	has probability at least $1 - \frac{C}{\epsilon}e^{-cn\epsilon^2}$.
	On the event~\eqref{eq:event_good0}, we have for all $w \in D_{\epsilon}$:
	$$
	\frac{1}{N} \| w - \mathsf{w}_{\lambda} \|^2 \geq 
	\frac{1}{N} \big( \|w\| - \|\mathsf{w}_{\lambda}\| \big)^2
	\geq
	\frac{1}{N} \big( \sqrt{N\epsilon} - \frac{1}{2} \sqrt{N\epsilon} \big)^2
	\geq \frac{\epsilon}{4} \,.
	$$
	This gives that on the event~\eqref{eq:event_good0}, for all $w \in D_{\epsilon}$, $L_{\lambda}(w) > \min\limits_{v \in \R^N} L_{\lambda}(v) + 3 \gamma \epsilon$. The intersection of~\eqref{eq:event_good0} with the event $\big\{ \min\limits_{w \in D_{\epsilon}} L_{\lambda}(w) \leq \min\limits_{w\in \R^N}L_{\lambda}(w) + 3 \gamma \epsilon \big\}$ is therefore empty: the lemma is proved.
\end{proof}


\subsection{Uniform control over \texorpdfstring{$\lambda$}{lambda}: proofs of Theorems~\ref{th:unif_lambda_law} and~\ref{th:unif_lambda_risk}-\eqref{eq:unif_risk}} \label{sec:proof_unif}

\subsubsection{Control of the \texorpdfstring{$\ell_1$}{l1}-norm of the Lasso estimator}

\begin{proposition}\label{prop:lipschitz_lambda_l1}
	Let $\xi >0, p>0$.
	Define $K=2 \xi + \frac{2\delta\sigma^2}{\lambda_{\rm min}}$. Then
	$$
	\forall \theta^{\star} \in \cF_p(\xi), \quad \P\left( \forall \lambda \geq \lambda_{\rm min}, \quad \frac{1}{N} \big| \what{w}_{\lambda}  \big| \leq K N^{(1/p - 1)_+} \right) \geq 1 - e^{-n/2} \,.
	$$
\end{proposition}
\begin{proof}
	Since $\cF_{p'}(\xi) \subset \cF_p(\xi)$ for $p' \geq p$, it suffices to prove the Proposition for $p \in (0,1]$: we suppose now to be in that case.
	With probability at least $1 - e^{-n/2}$ we have $\|z\| \leq 2 \sqrt{n}$ and therefore
	$\min \mathcal{L}_{\lambda} \leq \mathcal{L}_{\lambda}(\theta^{\star}) \leq 2 \sigma^2 + \frac{\lambda}{n} |\theta^{\star}|$ for all $\lambda \geq 0$. One has thus with probability at least $1 - e^{-n/2}$,
	$$
	\forall \theta^{\star} \in \cF_1(\xi), \ \forall \lambda > 0, \quad
	\frac{\lambda}{n} |\what{\theta}_{\lambda}| \leq \mathcal{L}_{\lambda}(\what{\theta}_{\lambda}) \leq 2 \sigma^2 + \frac{\lambda}{n} |\theta^{\star}| \,,
	$$
	which implies that $\frac{1}{N}|\what{\theta}_{\lambda}| \leq \frac{2 \delta \sigma^2}{\lambda} + \xi N^{1/p - 1}$ since $\frac{1}{N}|\theta^{\star}| \leq \frac{1}{N} \big(\sum_{i=1}^N |\theta^{\star}_i|^p\big)^{1/p} \leq \xi N^{1/p - 1}$.
\end{proof}

\begin{proposition}\label{prop:lipschitz_lambda_s}
	Assume that $0 <\delta < 1$ and $\sigma >0$. 
	Let $s < s_{\rm max}(\delta)$.
	Then, there exists constants $c,K>0$ such that
	$$
	\forall \theta^{\star} \in \cF_0(s), \quad \P\left( 
		\forall \lambda \in [\lambda_{\rm min},\lambda_{\rm max}], \quad \frac{1}{n} \big| |\what{w}_{\lambda} + \theta^{\star} | - |\theta^{\star}| \big| \leq K
	\right) \geq 1 - 2 e^{-cn} \,.
	$$
\end{proposition}

Proposition~\ref{prop:lipschitz_lambda_s} follows from the arguments of~\cite{tropp2015convex} that we reproduce below.

\begin{lemma}
	Suppose that $\theta^{\star} \in \cF_0(s)$ for some $s<s_{\rm max}(\delta)$.
	There exists constants $c,a >0$ that only depend on $(\delta,s)$ such that with probability at least $1 - e ^{-cn}$, for all $w\in\R^N$ such that $|w + \theta^{\star}| - |\theta^{\star}| \leq 0$ we have
	$$
	\left\| X w \right\|^2 \geq a \|w\|^2 \,.
	$$
\end{lemma}
\begin{proof}
	Define 
	$$
	\mathcal{K} = \mathsf{D}\big(| \cdot |,\theta^{\star}\big) = \bigcup_{r > 0} \Big\{ u \in \R^N \  \Big| \ |\theta^{\star} + ru| \leq |\theta^{\star}| \Big\} \,,
	$$
	the descent cone of the $\ell_1$-norm at $\theta^{\star}$.
	Define
	$$
	\nu_{\rm min}(X,\mathcal{K}) = \inf \big\{ \|X x \| \, \big| \, x \in \mathcal{K}, \|x\| = 1 \big\} \,,
	$$
	Let $\omega(\mathcal{K})$ be the Gaussian width of $\mathcal{K}$:
	$$
	\omega(\mathcal{K}) = \E \left[\sup_{u \in \mathcal{K}, \|u\| = 1} \langle g , u \rangle \right] \,,
	$$
	where the expectation is taken with respect to $g \sim \cN(0,\bbf{I}_N)$.
	The following result goes back to Gordon's work,~\cite{gordon1985some},~\cite{gordon1988milman}. It can be found in for instance~\cite{tropp2015convex} (Proposition~3.3).
	\begin{proposition}
		For all $t \geq 0$,
		$$
		\P \bigg(\sqrt{n} \, \nu_{\rm min}(X,\mathcal{K}) \geq \sqrt{n-1} - \omega(\mathcal{K}) -t \bigg) \geq 1 - e^{-t^2/2} \,.
		$$
	\end{proposition}

	Recall that $M_s(\alpha) = s(1+\alpha^2) + 2 (1-s) \big((1+\alpha^2) \Phi(-\alpha) - \alpha \phi(\alpha)\big)$ is the ``critical function'' studied in Section~\ref{sec:control_sparse_balls}.
	\begin{lemma}
		For all $\alpha \geq 0$
		$$
		\omega(\mathcal{K})^2 \leq N M_s(\alpha) = N \big( s(1+\alpha^2) + 2 (1-s) \big((1+\alpha^2) \Phi(-\alpha) - \alpha \phi(\alpha)\big) \big) \,.
		$$
	\end{lemma}
	\begin{proof}
		Let $v \in \partial |\theta^{\star}|$. By convexity, we have for all $w\in \mathcal{K}$ we have $\langle w , v \rangle \leq |w+\theta^{\star}| - |\theta^{\star}| \leq 0$.
		Now for all $x \in \R^N$ and $\alpha \geq 0$
		\begin{equation}\label{eq:ineq_v}
			\|x-\alpha v\| 
			= \sup_{\|w\|=1} \langle x-\alpha v , w \rangle
			\geq \sup_{w \in \mathcal{K}, \, \|w\|=1} \big\{ \langle w, x \rangle - \langle w , \alpha v \rangle \big\}
			\geq \sup_{w \in \mathcal{K}, \, \|w\|=1} \langle w, x \rangle \,.
		\end{equation}
		Let $S_0$ denote the support of $\theta^{\star}$.
		Let $g \sim \cN(0,\bbf{I}_N)$, $\alpha \geq 0 $ and define
		$$
		v_i
		= 
		\begin{cases}
			\text{sign}(\theta^{\star}_i) & \text{if} \ i \in S_0\,, \\
			\alpha^{-1} g_i \bbf{1}(|g_i| \leq \alpha) + \text{sign}(g_i) \bbf{1}(|g_i| > \alpha) & \text{otherwise.}
		\end{cases}
		$$
		Notice that $v \in \partial |\theta^{\star}|$, therefore by~\eqref{eq:ineq_v}:
		\begin{align*}
			\omega(\mathcal{K})^2 
			&\leq \E \left[\Big(\sup_{u \in \mathcal{K}, \|u\| = 1} \langle g , u \rangle\Big)^2 \right]
			\leq 
			\E \left[ \| g - \alpha v \|^2 \right]
			\leq
			\E \left[ \sum_{i \in S_0} (g_i - \alpha \text{sign}(\theta^{\star}_i))^2 + \sum_{i \not\in S_0}\eta(g_i,\alpha)^2 \right]
			\\
			&\leq
			N s (1+\alpha^2) + 2 N (1-s) ((1+\alpha^2)\Phi(-\alpha) - \alpha \phi(\alpha)) \,.
		\end{align*}
	\end{proof}
	\\

	Since $s \leq s_{\rm max}(\delta)$, there exists (see Lemma~\ref{lem:below_s_max}) $\alpha \geq 0$ and $t\in (0,1)$ such that $M_s(\alpha) \leq \delta (1-t)^2$. Consequently $\omega(\mathcal{K}) \leq \sqrt{n}(1-t)$.
	Therefore, there exists some constants $a,c>0$ that only depends on $s$ and $\delta$ such that
	$$
	\P \left(\nu_{\rm min}(X,\mathcal{K}) \geq a \right) \geq 1 - e^{-c n} \,.
	$$
	On the above event, for all $w \in \mathcal{K}$, $\|X w\|^2 \geq a^2 \|w\|^2$, which proves the Lemma.
\end{proof}
\\

\begin{proof}[of Proposition~\ref{prop:lipschitz_lambda_s}]
	Let us work on the event
	\begin{equation} \label{eq:event_cone}
		\Big\{ \forall w \in \R^N, \quad |w+\theta^{\star}| - |\theta^{\star}| \leq 0 \implies 
			\big\| X w \big\|^2 \geq a \|w\|^2
		\Big\} \ \bigcap \
		\Big\{ \|z\| \leq 2 \sqrt{n} \Big\} \,,
	\end{equation}
	which has probability at least $1 - 3 e^{- cn/2}$. 
	Let $\lambda \in [\lambda_{\rm min},\lambda_{\rm max}]$.
	Notice that on the event~\eqref{eq:event_cone} we have $\min \cC_{\lambda} \leq \cC_{\lambda}(0) \leq 2 \sigma^2$ and therefore
	$$
	\frac{\lambda}{n} \big(|\what{w}_{\lambda} + \theta^{\star}| - |\theta^{\star}| \big) \leq 2 \sigma^2 \,.
	$$
	We distinguish two cases:
	\\

	\textbf{Case 1}: $|\what{w}_{\lambda} + \theta^{\star}| - |\theta^{\star}| \geq 0$. In that case we obtain $\frac{1}{n} \big| |\what{w}_{\lambda} + \theta^{\star} | - |\theta^{\star}| \big| \leq \frac{2\sigma^2}{\lambda_{\rm min}}$.
	\\

	\textbf{Case 2}: $|\what{w}_{\lambda} + \theta^{\star}| - |\theta^{\star}| \leq 0$. In that case
	\begin{align*}
		2 \sigma^2 \geq \cC_{\lambda}(\what{w}_{\lambda}) 
	\geq 
	\frac{1}{2n} \left\| X \what{w}_{\lambda} - \sigma z \right\|^2 - \frac{\lambda}{n}|\what{w}_{\lambda}| 
	&\geq
	\frac{1}{4n} \| X \what{w}_{\lambda} \|^2 - \frac{\sigma^2}{2n} \| z \|^2 - \frac{\lambda\sqrt{N}}{n}\|\what{w}_{\lambda}\|
	\\
	&\geq \frac{a}{2n} \|\what{w}_{\lambda}\|^2 - 2 \sigma^2 - \frac{\lambda}{\sqrt{\delta n}}\|\what{w}_{\lambda}\| \,.
	\end{align*}
	This implies that there exists a constant $C=C(s,\delta,\sigma) > 0$ such that $\frac{1}{\sqrt{n}} \|\what{w}_{\lambda}\| \leq C (1+\lambda)$.
	One conclude
	$$
	-C\delta^{-1/2}(1+\lambda_{\rm max}) \leq -\frac{1}{n} |\what{w}_{\lambda}| \leq \frac{1}{n}\big(|\what{w}_{\lambda} + \theta^{\star}| - |\theta^{\star}|\big) \leq 0 \,.
	$$
\end{proof}

\subsubsection{Lipschitz continuity of the limiting risk and empirical distribution}

\begin{proposition}\label{prop:mu_lipschitz}
	The function $\lambda \mapsto \mu_{\lambda}^*$ is $M$-Lipschitz on $[\lambda_{\rm min},\lambda_{\rm max}]$ with respect to the Wasserstein distance $W_2$, for some constant $M=M(\Omega)>0$.
\end{proposition}
\begin{proof}
	Let $\lambda_1,\lambda_2 \in [\lambda_{\rm min},\lambda_{\rm max}]$.
	\begin{align*}
		W_2(\mu_{\lambda_1}^*,\mu_{\lambda_2}^*)^2
		&\leq 
		\E \left[
			\left(
				\eta(\Theta + \tau_*(\lambda_1)Z, \alpha_*(\lambda_1) \tau_*(\lambda_1))
				-
				\eta(\Theta + \tau_*(\lambda_2)Z, \alpha_*(\lambda_2) \tau_*(\lambda_2))
			\right)^2
		\right]
		\\
		&\leq 
		2\E \left[
			\left(\tau_*(\lambda_1) Z - \tau_*(\lambda_2) Z
			\right)^2
			+ 
			(\alpha_*(\lambda_1) \tau_*(\lambda_1) - \alpha_*(\lambda_2) \tau_*(\lambda_2))^2
		\right]
		\\
		&\leq 
		2 \left(\tau_*(\lambda_1) - \tau_*(\lambda_2) \right)^2 + 
		2	(\alpha_*(\lambda_1) \tau_*(\lambda_1) - \alpha_*(\lambda_2) \tau_*(\lambda_2))^2
		\\
		&\leq 
		2 (1 + \alpha_{\rm max}^2) \left(\tau_*(\lambda_1) - \tau_*(\lambda_2) \right)^2 + 
		2 \tau_{\rm max}^2	(\alpha_*(\lambda_1) - \alpha_*(\lambda_2) )^2 \,.
	\end{align*}
	Since by Proposition~\ref{prop:dep_lambda} the functions $\lambda \mapsto \alpha_*(\lambda)$ and $\lambda \mapsto \tau_*(\lambda)$ are both $M$-Lipschitz on $[\lambda_{\rm min},\lambda_{\rm max}]$, for some constant $M=M(\Omega) >0$, we obtain:
	$$
	W_2(\mu_{\lambda_1}^*,\mu_{\lambda_2}^*)^2
	\leq 2 M^2 (1 + \alpha_{\rm max}^2 + \tau_{\rm max}^2) (\lambda_1 - \lambda_2)^2 \,,
	$$
	which proves the Lemma.
\end{proof}

\begin{proposition}\label{prop:risk_lipschitz}
The function $\lambda \mapsto R_*(\lambda)= \delta (\tau_*(\lambda)^2 - \sigma^2)$ is $M$-Lipschitz on $[\lambda_{\rm min},\lambda_{\rm max}]$, for some constant $M=M(\Omega)>0$.
\end{proposition}
\begin{proof}
	This is a consequence of Proposition~\ref{prop:dep_lambda}.
\end{proof}

\subsubsection{Proofs of Theorems~\ref{th:unif_lambda_law} and~\ref{th:unif_lambda_risk} }

\begin{lemma}\label{lem:lambda_lip}
	Assume that $\mathcal{D}$ is either $\cF_0(s)$ or $\cF_p(\xi)$ for some $s < s_{\rm max}(\delta)$ and $\xi \geq 0, p >0$. Define
	$$
	q=
	\begin{cases}
		(1/p - 1)_+ & \text{if} \quad \mathcal{D} = \cF_p(\xi)\,, \\
		0 & \text{if} \quad \mathcal{D} = \cF_0(s)\,.
	\end{cases}
	$$
	Then there exists constants $K,C,c>0$ that depend only on $\Omega$ such that for all $\theta^{\star} \in \mathcal{D}$
	\begin{equation}\label{eq:lambda_lip}
		\P \Big(
			\forall \lambda,\lambda' \in [\lambda_{\rm min},\lambda_{\rm max}],
			\quad
			\mathcal{L}_{\lambda'}(\what{\theta}_{\lambda}) \leq \min_{x \in \R^N} \mathcal{L}_{\lambda'}(x) + KN^q |\lambda - \lambda'|
		\Big)
		\geq 1 - C e^{-cn} \,.
	\end{equation}
\end{lemma}
\begin{proof}
	$K=K(\Omega)>0$ be a constant such that for all $\theta^{\star} \in \mathcal{D}$, the event
	\begin{equation}\label{eq:event_lipschitz_lambda}
		\Big\{ \forall \lambda \in [\lambda_{\rm min},\lambda_{\rm max}], \ \frac{1}{n} \big| |\what{\theta}_{\lambda}  | - |\theta^{\star}| \big| \leq K N^{q} \Big\}
	\end{equation}
	has probability at least $1 - C e^{-cn}$. Such $K$ exists by Propositions~\ref{prop:lipschitz_lambda_l1} and~\ref{prop:lipschitz_lambda_s}.
	On the event~\eqref{eq:event_lipschitz_lambda} we have for all $\lambda,\lambda' \in [\lambda_{\rm min},\lambda_{\rm max}]$:
	\begin{align*}
		\mathcal{L}_{\lambda'}(\what{\theta}_{\lambda}) 
		&= \mathcal{L}_{\lambda}(\what{\theta}_{\lambda}) + (\lambda' - \lambda) \frac{1}{n} |\what{\theta}_{\lambda} |
		\\
		&\leq \mathcal{L}_{\lambda}(\what{\theta}_{\lambda'}) + (\lambda' - \lambda) \frac{1}{n} |\what{\theta}_{\lambda}|
		= \mathcal{L}_{\lambda'}(\what{\theta}_{\lambda'}) + \frac{1}{n}(\lambda - \lambda') \big( |\what{\theta}_{\lambda'} | - |\what{\theta}_{\lambda}| \big) 
		\\
		&\leq \min_{\theta \in \R^N} \mathcal{L}_{\lambda'}(\theta) + \frac{1}{n} |\lambda' - \lambda| \Big(\big| |\what{\theta}_{\lambda'} | - |\theta^{\star}| \big| + \big| |\what{\theta}_{\lambda}  | - |\theta^{\star}| \big| \Big)
		\\
		&\leq \min_{\theta \in \R^N} \mathcal{L}_{\lambda'}(\theta) + 2 K N^q |\lambda - \lambda'| \,.
	\end{align*}
\end{proof}
\\

Theorem~\ref{th:unif_lambda_law} and Theorem~\ref{th:unif_lambda_risk}-\eqref{eq:unif_risk} are proved the same way.

\begin{proof}[of Theorem~\ref{th:unif_lambda_law}]
	Let $\gamma > 0$ as given by Theorem~\ref{th:key_law} and let $K=K(\Omega)>0$ as given by Lemma~\ref{lem:lambda_lip}.
	Let $M=M(\Omega)>0$ such that $\lambda \mapsto \mu_{\lambda}^*$ is $M$-Lipschitz with respect to the Wasserstein distance $W_2$ on $[\lambda_{\rm min},\lambda_{\rm max}]$, as given by Proposition~\ref{prop:mu_lipschitz}.
	\\

	Let $\epsilon \in (0,1]$ and define $\epsilon' = \min \left(\frac{\gamma \epsilon}{2KN^q}, \frac{\epsilon}{M+1} \right)$.
	Let $k= \big\lceil (\lambda_{\rm max} - \lambda_{\rm min}) / \epsilon' \big\rceil$. Define, for $i = 0, \dots,k$:
	$$
	\lambda_{i} = \lambda_{\rm min} + i \epsilon' \,.
	$$
	By Theorem~\ref{th:key_law}, the event
	\begin{equation}\label{eq:event_strongly_convex}
		\Big\{
			\forall i \in \{1,\dots,k\}, \
			\forall \theta \in \R^N, \quad 
			\mathcal{L}_{\lambda_i}(\theta) \leq \min_{x \in \R^N} \mathcal{L}_{\lambda_i}(x) + \gamma \epsilon
			\ \implies \ 
			W_2(\what{\mu}_{(\theta,\theta^{\star})},\mu^*_{\lambda_i})^2 \leq \epsilon 
		\Big\}
	\end{equation}
	has probability at least $1 - k C \epsilon^{-\max(1,a)} e^{-cN \epsilon^2 \epsilon^a \log(\epsilon)^{-2}} \geq 1 - CN^q \epsilon^{-\max(1,a)-1} e^{-cN \epsilon^2 \epsilon^a \log(\epsilon)^{-2}}$.
	Therefore, on the intersection of the event in~\eqref{eq:lambda_lip} and the event~\eqref{eq:event_strongly_convex} we have for all $\lambda \in [\lambda_{\rm min},\lambda_{\rm max}]$
	$$
	\mathcal{L}_{\lambda_i}(\what{\theta}_\lambda) 
	\leq \min_{x \in \R^N} \mathcal{L}_{\lambda_i}(x) + 2 K N^q | \lambda - \lambda_i |
	\leq \min_{x \in \R^N} \mathcal{L}_{\lambda_i}(x) + \gamma \epsilon \,,
	$$
	where $1 \leq i \leq k$ is such that $\lambda \in [\lambda_{i-1},\lambda_i]$.
	This implies (since we are on the event~\eqref{eq:event_strongly_convex}) that $W_2\big(\what{\mu}_{(\what{\theta}_{\lambda},\theta^{\star})}, \mu_{\lambda_i}^*\big)^2 \leq \epsilon$.
	We conclude by
	\begin{align*}
		W_2\big(\what{\mu}_{(\what{\theta}_{\lambda},\theta^{\star})}, \mu^*_{\lambda}\big)^2
		&\leq 
		2 W_2\big(\what{\mu}_{(\what{\theta}_{\lambda},\theta^{\star})}, \mu^*_{\lambda_i}\big)^2
		+
		2 W_2(\mu_{\lambda}^*, \mu^*_{\lambda_i})^2
		\leq 2 \epsilon + 2 M^2 (\lambda-\lambda_i)^2
		\leq 
		4 \epsilon \,.
	\end{align*}
	This proves the Theorem.
\end{proof}
\\

\begin{proof}[of Theorem~\ref{th:unif_lambda_risk}-\eqref{eq:unif_risk}]
	Let $\gamma > 0$ as given by Theorem~\ref{th:key_risk} and let $K=K(\Omega)>0$ as given by Lemma~\ref{lem:lambda_lip}.
	Let $M=M(\Omega)>0$ such that $\lambda \mapsto R_*(\lambda)$ is $M$-Lipschitz on $[\lambda_{\rm min},\lambda_{\rm max}]$, as given by Proposition~\ref{prop:risk_lipschitz}.
	\\

	Let $\epsilon \in (0,1]$ and define $\epsilon' = \min \left(\frac{\gamma \epsilon}{2K N^q}, \frac{\epsilon}{M+1} \right)$.
	Let $k= \big\lceil (\lambda_{\rm max} - \lambda_{\rm min}) / \epsilon' \big\rceil$. Define, for $i = 0, \dots,k$:
	$$
	\lambda_{i} = \lambda_{\rm min} + i \epsilon' \,.
	$$
	By Theorem~\ref{th:key_risk}, the event
	\begin{equation}\label{eq:event_strongly_convex_risk}
		\left\{
			\forall i \in \{1,\dots,k\}, \
			\forall \theta \in \R^N, \quad 
			\mathcal{L}_{\lambda_i}(\theta) \leq \min_{x \in \R^N} \mathcal{L}_{\lambda_i}(x) + \gamma \epsilon
			\ \implies \ 
			\Big(\frac{1}{N} \|\theta - \theta^{\star}\|^2 - R_*(\lambda_i) \Big)^2 \leq \epsilon 
		\right\}
	\end{equation}
	has probability at least $1 - k C \epsilon^{-1} e^{-cN \epsilon^2} \geq 1 - C N^q \epsilon^{-2} e^{-cN \epsilon^2 }$.
	Therefore, on the intersection of the event in~\eqref{eq:lambda_lip} and the event~\eqref{eq:event_strongly_convex_risk} we have for all $\lambda \in [\lambda_{\rm min},\lambda_{\rm max}]$
	$$
	\mathcal{L}_{\lambda_i}(\what{\theta}_\lambda) 
	\leq \min_{\theta \in \R^N} \mathcal{L}_{\lambda_i}(\theta) + 2 K N^q | \lambda - \lambda_i |
	\leq \min_{\theta \in \R^N} \mathcal{L}_{\lambda_i}(\theta) + \gamma \epsilon \,,
	$$
	where $1 \leq i \leq k$ is such that $\lambda \in [\lambda_{i-1},\lambda_i]$.
	This implies (since we are on the event~\eqref{eq:event_strongly_convex_risk}) that
	$\Big(\frac{1}{N} \|\what{\theta}_{\lambda} - \theta^{\star}\|^2 - R_*(\lambda_i) \Big)^2 \leq \epsilon$.
	We conclude by
	\begin{align*}
		\Big(\frac{1}{N} \|\what{\theta}_{\lambda} - \theta^{\star}\|^2 - R_*(\lambda) \Big)^2
		\leq
		2 \Big(\frac{1}{N} \|\what{\theta}_{\lambda} - \theta^{\star}\|^2 - R_*(\lambda_i) \Big)^2
		+
		2 \Big( R_*(\lambda_i) - R_*(\lambda) \Big)^2
		\leq 2 \epsilon + 2 M^2 (\epsilon')^2
		\leq 
		4 \epsilon \,.
	\end{align*}
	This proves~\eqref{eq:unif_risk}.
\end{proof}

\section{Study of the Lasso residual: proof of~\eqref{eq:unif_beta}-\eqref{eq:unif_prediction}} \label{sec:gordon_u}

This Section is devoted to the proof of~\eqref{eq:unif_beta}-\eqref{eq:unif_prediction} from Theorem~\ref{th:unif_lambda_risk}.
Let us define 
$$
\what{u}_{\lambda} = X \what{w}_{\lambda} - \sigma z = X \what{\theta}_{\lambda} -y \,.
$$
$\what{u}_{\lambda}$ is the unique maximizer of 
$$
u \mapsto \min_{w \in \R^N} \Big\{ u^{\sT} X w - \sigma u^{\sT} z - \frac{1}{2} \|u\|^2 + \lambda (|\theta^{\star} + w| - |\theta^{\star}|) \Big\} \,.
$$
In Section~\ref{sec:proof_norm_u} below, we prove the following Theorem:
\begin{theorem}\label{th:norm_u}
	There exists constants $c,C>0$ such that for all $\epsilon \in (0,1]$, all $\theta^{\star} \in \mathcal{D}$ and all $\lambda \in [\lambda_{\rm min},\lambda_{\rm max}]$
	$$
	\P \left(
		\Big( \frac{1}{n} \| \what{u}_{\lambda} \|^2  - \beta_*(\lambda)^2 \Big)^2 \geq \epsilon
	\right) \leq \frac{C}{\epsilon} e^{-cn\epsilon^2} \,,
	$$
	and
	$$
	\P \left(
	\Big( \frac{1}{n} \| \what{u}_{\lambda} + \sigma z \|^2  - P_*(\lambda) \Big)^2 \geq \epsilon
	\right) \leq \frac{C}{\epsilon} e^{-cn\epsilon^2} \,.
	$$
\end{theorem}

Theorem~\ref{th:unif_lambda_risk}-\eqref{eq:unif_beta}-\eqref{eq:unif_prediction} will then be deduced from Theorem~\ref{th:norm_u} in Section~\ref{sec:proof_norm_u}.

\subsection{Study of Gordon's optimization problem}

Recall that $g \sim \cN(0,I_N)$ and $h \sim \cN(0,I_n)$ are independent standard Gaussian vectors.
Define for $(w,u) \in \R^N \times \R^n$,
$$
		m_{\lambda}(w,u) = 
		-\frac{1}{n^{3/2}} \|u\| g^{\sT} w 
		+\frac{1}{n^{3/2}} \|w\| h^{\sT} u
		- \frac{\sigma}{n} u^{\sT} z - \frac{1}{2n} \|u\|^2 + 
		\frac{\lambda}{n} \big(  |w+\theta^{\star}| - |\theta^{\star}| \big)\,,
		$$
		and $U_{\lambda}(u) = \min_{w \in \R^N} m_{\lambda}(w,u)$, $\widetilde{U}_{\lambda}(u) = m_{\lambda}(\mathsf{w}_{\lambda},u)$, where $\mathsf{w}_{\lambda}$ is defined by~\eqref{eq:def_w_sf}. Obviously we have $U_{\lambda}(u) \leq \widetilde{U}_{\lambda}(u)$. We write also
		\begin{equation}
			\mathsf{u}_{\lambda} = 
			\frac{\beta_*(\lambda)}{\tau_*(\lambda)} \Big(\sqrt{\tau_*(\lambda)^2 - \sigma^2} \frac{h}{\sqrt{n}} - \frac{\sigma}{\sqrt{n}} z\Big) \,.
		\end{equation}
		\begin{lemma}\label{lem:U_tilde}
			There exists constants $C,c>0$ such that for all $\epsilon \in (0,1]$ and any $\lambda \in [\lambda_{\rm min},\lambda_{\rm max}]$ we have with probability at least $1- C e^{-cn\epsilon^2}$
			\begin{itemize}
				\item $\widetilde{U}_{\lambda}$ is $1/n$-strongly concave on $\R^n$ and admits therefore a unique maximizer $u^*_{\lambda}$ over $\R^n$.
				\item $|\max_{u \in \R^n} \widetilde{U}_{\lambda}(u) - L_*(\lambda) | \leq \epsilon$, where $L_*(\lambda)$ is defined by Corollary~\ref{cor:concentration_gordon}.
				\item $\frac{1}{n}\| u^*_{\lambda} - \mathsf{u}_{\lambda}\|^2 \leq \epsilon$.
			\end{itemize}
		\end{lemma}
		\begin{proof}
			By Lemma~\ref{lem:conc_w_RS} and Lemma~\ref{lem:opti_w_star}, $\frac{1}{n} \mathsf{w}_{\lambda}^{\sT} g$ concentrates around $\frac{s_*(\lambda)}{\delta}$ which is greater than some constant $\gamma >0$. 
			Indeed
			$$
			s_*(\lambda) = \E \left[
				\Phi\Big(\frac{\Theta}{\tau_*} - \alpha_* \Big)
				+
				\Phi\Big(-\frac{\Theta}{\tau_*} - \alpha_* \Big)
			\right]
			$$
			remains greater than some strictly positive constant while $\theta^{\star}$ vary in $\mathcal{D}$ and $\lambda$ vary in $[\lambda_{\rm min},\lambda_{\rm max}]$.
			By Lemma~\ref{lem:conc_w_RS} we have then that with probability at least $1 - Ce^{-cn}$, $\mathsf{w}_{\lambda}^{\sT} g \geq 0$ which implies that $\widetilde{U}_{\lambda}$ is $1/n$-strongly concave.
			Let us compute
			\begin{align*}
				\max_{u \in \R^n} \widetilde{U}_{\lambda}(u)
				&= \max_{\beta \geq 0}
				\Big\{
					\Big(
					\Big\| \frac{1}{n} \|\mathsf{w}_{\lambda}\| h - \frac{\sigma}{\sqrt{n}} z \Big\|
					-\frac{1}{n} g^{\sT} \mathsf{w}_{\lambda} 
				\Big) \beta
					- \frac{1}{2} \beta^2 + 
				\frac{\lambda}{n} \big(  |\mathsf{w}_{\lambda}+\theta^{\star}| - |\theta^{\star}| \big) \Big\}
				\\
				&= 
				\frac{1}{2}
				\left(
					\Big\| \frac{1}{n} \|\mathsf{w}_{\lambda}\| h - \frac{\sigma}{\sqrt{n}} z \Big\|
					-\frac{1}{n} g^{\sT} \mathsf{w}_{\lambda} 
				\right)_{\!\!+}^2
				+\frac{\lambda}{n} \big(  |\mathsf{w}_{\lambda}+\theta^{\star}| - |\theta^{\star}| \big) \,.
			\end{align*}
			By the concentration properties of $\mathsf{w}_{\lambda}$ (see Section~\ref{sec:concentration_w}), have that with probability at least $1- C e^{-cn\epsilon^2}$,
			$|\max_{u \in \R^n} \widetilde{U}_{\lambda}(u) - L_*(\lambda) | \leq \epsilon$. One verify analogously that $\widetilde{U}_{\lambda}(\mathsf{u}_{\lambda}) \geq L_*(\lambda) - \epsilon$ with the same probability, which implies the third point by strong concavity.
		\end{proof}

		\subsection{Proof of Theorem~\ref{th:norm_u}} \label{sec:proof_norm_u}
		Let us only prove the second point since the first one follows from the same arguments. 
		Let $\epsilon \in (0,1]$ and define
		$$
		D_{\epsilon} = \Big\{ u \in \R^n \, \Big| \, \big| \frac{1}{\sqrt{n}} \|u + \sigma z \|  - \sqrt{P_*(\lambda)} \big| \geq 6\epsilon^{1/2} \Big\} \,.
		$$
Let us define for $(w,u) \in \R^N \times \R^n$:
\begin{equation}\label{eq:def_c}
		c_{\lambda}(w,u) = \frac{1}{n} u^{\sT} Xw - \frac{\sigma}{n} u^{\sT}z - \frac{1}{2n} \|u\|^2 + \frac{\lambda}{n}\big( |w + \theta^{\star}| - |\theta^{\star}| \big)
	\end{equation}
	\begin{lemma}\label{lem:minmax_c}
		We have almost surely
		$$
		\min_{w \in \R^N} \max_{u \in \R^n} c_{\lambda}(w,u)
		=
		\max_{u \in \R^n} \min_{w \in \R^N} c_{\lambda}(w,u)
		= c_{\lambda}(\what{w}_{\lambda},\what{u}_{\lambda}) \,.
		$$
	\end{lemma}
	\begin{proof}
		By definition of $\what{w}_{\lambda}$ and $\what{u}_{\lambda}$ we have
		$$
c_{\lambda}(\what{w}_{\lambda},\what{u}_{\lambda})
=
		\min_{w \in \R^N} \max_{u \in \R^n} c_{\lambda}(w,u)
		\geq
		\max_{u \in \R^n} \min_{w \in \R^N} c_{\lambda}(w,u) \,.
		$$
		Let us prove the converse inequality. 
		The optimality condition of $\what{w}_{\lambda}$ gives that there exists $v \in \partial |\theta^{\star} + \what{w}_{\lambda}|$ such that
		$$
		X^{\sT} \what{u}_{\lambda} + \lambda v =
		X^{\sT} (X \what{w}_{\lambda} - \sigma z) + \lambda v = 0 \,.
		$$
		The function $w \mapsto c_{\lambda}(w,\what{u}_{\lambda})$ is convex and 
		$$
		\frac{1}{n} X^{\sT} \what{u}_{\lambda} + \frac{\lambda}{n}v = 0
		$$
	is a subgradient at $\what{w}_{\lambda}$. Therefore $\min_{w \in \R^N} c_{\lambda}(w,\what{u}_{\lambda}) = c_{\lambda}(\what{w}_{\lambda},\what{u}_{\lambda})$, which proves the lemma.
	\end{proof}
	\\

	We compute now
		\begin{align*}
		\P\big(\what{u}_{\lambda} \in D_{\epsilon} \big) 
		&=
		\P\Big(\max_{u \in D_{\epsilon}} \min_{w\in \R^N} c_{\lambda}(w,u) \geq \max_{u \in \R^n} \min_{w\in \R^N} c_{\lambda}(w,u) \Big) 
		\\
		&\leq 
		\P\Big(\max_{u \in D_{\epsilon}} \min_{w\in \R^N} c_{\lambda}(w,u) \geq L_*(\lambda) - \epsilon \Big) + 
			\P \Big(\max_{u \in \R^n} \min_{w\in \R^N} c_{\lambda}(w,u) \leq L_*(\lambda) - \epsilon \Big)  \,.
		\end{align*}
		By Lemma~\ref{lem:minmax_c} and Proposition~\ref{prop:concentration_min_c} we can bound
		$$
			\P \Big(\max_{u \in \R^n} \min_{w\in \R^N} c_{\lambda}(w,u) \leq L_*(\lambda) - \epsilon \Big) 
			=
			\P \Big(\min_{w\in \R^N} \cC_{\lambda}(w) \leq L_*(\lambda) - \epsilon \Big) 
			\leq \frac{C}{\epsilon} e^{-cn\epsilon^2} \,.
			$$
			Now by the same reasoning than Corollary~\ref{cor:gordon} (we omit here the details for the sake of brevity) we have
			$$
		\P\Big(\max_{u \in D_{\epsilon}} \min_{w\in \R^N} c_{\lambda}(w,u) \geq L_*(\lambda) - \epsilon \Big) 
		\leq
		2\P\Big(\max_{u \in D_{\epsilon}} \min_{w\in \R^N} m_{\lambda}(w,u) \geq L_*(\lambda) - \epsilon \Big) 
		=
		2\P\Big(\max_{u \in D_{\epsilon}} U_{\lambda}(u) \geq L_*(\lambda) - \epsilon \Big)  \,.
		$$
		Since $U_{\lambda} \leq \widetilde{U}_{\lambda}$ we obtain
			$$
		\P\Big(\max_{u \in D_{\epsilon}} \min_{w\in \R^N} c_{\lambda}(w,u) \geq L_*(\lambda) - \epsilon \Big) 
		\leq
		2\P\Big(\max_{u \in D_{\epsilon}} \widetilde{U}_{\lambda}(u) \geq L_*(\lambda) - \epsilon \Big) \,.
		$$
		Let $E$ be the event of Lemma~\ref{lem:U_tilde} above and let us work on the event
		\begin{equation}
			E \
			\bigcap \
			\Big\{  \big| \frac{1}{\sqrt{n}} \| \mathsf{u}_{\lambda} + \sigma z \| - \sqrt{P_*(\lambda)} \big| \leq \epsilon^{1/2} \Big\} \,,
		\end{equation}
		which has probability at least $1-\frac{C}{\epsilon} e^{-cn\epsilon^2}$ (the fact that the second event in the intersection has this probability follows from standard concentration arguments as in Section~\ref{sec:concentration_w}).
		Let now $u \in D_{\epsilon}$, by the definition of $D_{\epsilon}$ and the event above we have $\frac{1}{\sqrt{n}} \|u - \mathsf{u}_{\lambda} \| \geq 5\epsilon^{1/2}$ and thus $\frac{1}{\sqrt{n}} \|u - u^*_{\lambda} \| \geq 4\epsilon^{1/2}$. By $1/n$-strong concavity of $\widetilde{U}_{\lambda}$ we get
		$$
		\widetilde{U}_{\lambda}(u) \leq \max_{u' \in \R^n} \widetilde{U}_{\lambda}(u') - 8 \epsilon \leq L_*(\lambda) -7 \epsilon \,.
		$$
		Consequently $\displaystyle \P\Big(\max_{u \in D_{\epsilon}} \widetilde{U}_{\lambda}(u) \geq L_*(\lambda) - \epsilon \Big) \leq \frac{C}{\epsilon} e^{-cn\epsilon^2}$, which proves the result.

		\subsection{Uniform control over \texorpdfstring{$\lambda$}{lambda}: proof of Theorem~\ref{th:unif_lambda_risk}-\eqref{eq:unif_beta}-\eqref{eq:unif_prediction}} \label{sec:proof_unif_u}

	Let $\mathcal{D}$ be either $\cF_0(s)$ for some $s < s_{\rm max}(\delta)$ or $\cF_p(\xi)$ for some $\xi \geq 0$, $p >0$. 
	Let $q=0$ if $\mathcal{D} = \cF_0(s)$ and $q=(1/p-1)_+$ if $\mathcal{D} = \cF_p(\xi)$.
	By Propositions~\ref{prop:lipschitz_lambda_l1} and~\ref{prop:lipschitz_lambda_s} there exists a constant $K= K(\Omega)$ such that the event
		\begin{equation}\label{eq:event_control_l1}
			\Big\{ \forall \lambda \in [\lambda_{\rm min},\lambda_{\rm max}], \ 
			\frac{1}{n} \big| |\what{w}_{\lambda} + \theta^{\star}| - |\theta^{\star}| \big| \leq K N^q \Big\}
		\end{equation}
		has probability at least $1 - C e^{-cn}$. Let us fix this constant $K$ and let us write
		$$
		D_K = \Big\{ w \in \R^N \, \Big| \, 
			\frac{1}{n} \big| |w + \theta^{\star}| - |\theta^{\star}| \big| \leq K N^q \Big\} \,.
		$$
		We define also
		$$
		\mathcal{U}_{\lambda}(u) = \min_{w \in D_K} \Big\{
		\frac{1}{n} u^{\sT} X w - \frac{\sigma}{n} u^{\sT} z - \frac{1}{2n}\|u\|^2 + \frac{\lambda}{n} \big(|w + \theta^{\star}| - |\theta^{\star}|\big) \Big\} \,.
		$$
		\begin{lemma}
			The function $\mathcal{U}_{\lambda}$ is $1/n$-strongly concave. On the event~\eqref{eq:event_control_l1}, $\what{u}_{\lambda}$ is the (unique) maximizer of $\mathcal{U}_{\lambda}$.
		\end{lemma}
		\begin{proof}
			Let us work on the event~\eqref{eq:event_control_l1} and let $\lambda \in [\lambda_{\rm min},\lambda_{\rm max}]$. We have, by permutation of max and min:
			$$
			\max_{u \in \R^n} \mathcal{U}_{\lambda}(u) \leq \min_{w \in D_K} \mathcal{C}_{\lambda}(w) = \cC_{\lambda}(\what{w}_{\lambda}) \,,
			$$
			because on the event~\eqref{eq:event_control_l1}, $\what{w}_{\lambda}$ (the minimizer of $\cC_{\lambda}$) is in $D_K$. By the optimality condition of $\what{w}_{\lambda}$, one verify easily that $\mathcal{U}_{\lambda}(\what{u}_{\lambda}) = \cC_{\lambda}(\what{w}_{\lambda})$ which proves the lemma.
		\end{proof}
		\\

Theorem~\ref{th:unif_lambda_risk}-\eqref{eq:unif_beta}-\eqref{eq:unif_prediction} follow then easily from Theorem~\ref{th:norm_u} (by an $\epsilon$-net argument as in the proof of Theorems~\ref{th:unif_lambda_law} and ~\ref{th:unif_lambda_risk}-\eqref{eq:unif_risk}, see Section~\ref{sec:proof_unif}) and the following Proposition:
		\begin{proposition}\label{prop:lipschitz_u_hat}
	Let $q=0$ if $\mathcal{D} = \cF_0(s)$ and $q=(1/p-1)_+$ if $\mathcal{D} = \cF_p(\xi)$.
			There exists constants $C,c, \kappa >0$ such that for all $\theta^{\star} \in \mathcal{D}$ the following event
			$$
			\Big\{
				\forall \lambda,\lambda' \in [\lambda_{\rm min},\lambda_{\rm max}], \
				\frac{1}{n} \| \what{u}_{\lambda} - \what{u}_{\lambda'} \|^2 \leq \kappa N^q |\lambda - \lambda'|
			\Big\}
			$$
			has probability at least $1- Ce^{-cn}$.
		\end{proposition}
		\begin{proof}
			Let us work on the event~\eqref{eq:event_control_l1}, which has probability at least $1- Ce^{-cn}$. Let $\lambda,\lambda' \in [\lambda_{\rm min},\lambda_{\rm max}]$. We have
			$$
			\sup_{u \in \R^n} \big| \mathcal{U}_{\lambda}(u) - \mathcal{U}_{\lambda'}(u) \big| \leq \sup_{w \in D_K} \Big| \frac{\lambda - \lambda'}{n} (|w+\theta^{\star}|-|\theta^{\star}| ) \Big|
			\leq K N^q |\lambda - \lambda'| \,.
			$$
			Therefore
			\begin{align*}
				\mathcal{U}_{\lambda'}(\what{u}_{\lambda})
				\geq
				\mathcal{U}_{\lambda}(\what{u}_{\lambda}) - N^qK|\lambda-\lambda'|
				\geq
				\mathcal{U}_{\lambda}(\what{u}_{\lambda'}) - N^qK|\lambda-\lambda'|
				\geq
				\mathcal{U}_{\lambda'}(\what{u}_{\lambda'}) - 2KN^q|\lambda-\lambda'| \,,
			\end{align*}
			which gives that $\frac{1}{n}\|\what{u}_{\lambda} - \what{u}_{\lambda'} \|^2 \leq 4 KN^q |\lambda -\lambda'|$ by $1/n$-strong concavity.
		\end{proof}

		\section{Study of the subgradient \texorpdfstring{$\what{v}_{\lambda}$}{v}}\label{sec:gordon_v}

The goal of this section is to analyze the vector
$$
\what{v}_{\lambda} = \frac{1}{\lambda} X^{\sT} (y - X \what{\theta}_{\lambda}) \,,
$$
which is a subgradient of the $\ell_1$-norm at $\what{\theta}_{\lambda}$. Let us define
$$
B_{\infty}(0,1) = \Big\{ v \in \R^N \, \Big| \, \|v\|_{\infty} \leq 1 \Big\} \,.
$$
\subsection{Main results}

Let $B = \big\{w\in\R^N \, | \, |w| \leq 2 |\theta^{\star}| + 5 \sigma^2 \lambda_{\rm min}^{-1} n +K \big\}$, where $K >0$ is some constant (depending only on $\Omega$) that will be fixed later in the analysis (in fact $K$ is the constant given by Lemma~\ref{lem:kappa}).
Notice that $\what{w}_{\lambda} \in B$, with probability at least $1 - e^{-n/2}$.
Define 
$$
\mathcal{V}_{\lambda}(v)= \min_{w \in B}  \left\{
\frac{1}{2n} \| X w - \sigma  z\|^2 + \frac{\lambda}{n} v^{\sT}(\theta^{\star}+w) - \frac{\lambda}{n}|\theta^{\star}| \right\} \,.
$$
\begin{lemma}\label{lem:v_is_opt}
	With probability at least $1-e^{-n/2}$ we have for all $\lambda \geq \lambda_{\rm min}$
	$$
	\min_{w \in \R^N} \cC_{\lambda}(w)	=
	\max_{\|v\|_{\infty}\leq 1} \mathcal{V}_{\lambda}(v)
	$$
	and $\what{v}_{\lambda} = - \lambda^{-1} X^{\sT}(X\what{w}_{\lambda}- \sigma z)$ is a maximizer of $\mathcal{V}_{\lambda}$.
\end{lemma}
\begin{proof}
	Let us work on the event $\{ \|z\| \leq 2 \sqrt{n} \}$ which has probability at least $1-e^{-n/2}$. On this event we have $\what{w}_{\lambda} \in B$ and therefore
	$$
	\min_{w \in \R^N} \cC_{\lambda}(w) 
	=\min_{w \in B} \cC_{\lambda}(w) = \max_{\|v\|_{\infty} \leq 1} \mathcal{V}_{\lambda}(v)\,,
	$$
	where the permutation of the min-max is authorized by Proposition~\ref{prop:rockafellar}.
	The optimality condition of $\what{w}_{\lambda}$ gives that
	$$
	\what{v}_{\lambda}= - \lambda^{-1} X^{\sT}(X\what{w}_{\lambda}- \sigma z) \in \partial |\theta^{\star} + \what{w}_{\lambda}| \,.
	$$
	Therefore $\what{v}_{\lambda}^{\sT} (\what{w}_{\lambda} + \theta^{\star})= |\what{w}_{\lambda} + \theta^{\star}|$. Using the optimality condition again we obtain
	\begin{align*}
		\mathcal{V}_{\lambda}(\what{v}_{\lambda}) = \min_{w \in B}  \left\{
		\frac{1}{2} \| X w - \sigma  z\|^2 + \lambda \what{v}_{\lambda}^{\sT}(\theta^{\star}+w)\right\}
		&=\frac{1}{2} \| X \what{w}_{\lambda} - \sigma  z\|^2 + \lambda \what{v}_{\lambda}^{\sT}(\theta^{\star}+\what{w}_{\lambda})
		\\
		&=\frac{1}{2} \| X \what{w}_{\lambda} - \sigma  z\|^2 + \lambda |\theta^{\star}+\what{w}_{\lambda}| \,.
	\end{align*}
	Therefore $\what{v}_{\lambda}$ achieves the optimal value.
\end{proof}

\subsubsection{The empirical law of the subgradient}
Let $\nu^*_{\lambda}$ be the law of the couple
\begin{equation}\label{eq:def_nu_star}
\left( 
	- \frac{1}{\alpha_*(\lambda)\tau_*(\lambda)} \Big(\eta\big(\Theta + \tau_*(\lambda)Z, \alpha_*(\lambda) \tau_*(\lambda) \big) - \Theta - \tau_*(\lambda) Z \Big), \ \Theta
\right) \,,
\end{equation}
where $(\Theta,Z) \sim \what{\mu}_{\theta^{\star}} \otimes \cN(0,1)$. For $v \in \R^N$ we define
$$
\what{\mu}_{(v,\theta^{\star})} = \frac{1}{N} \sum_{i=1}^N \delta_{(v_i,\theta^{\star}_i)} \,.
$$

\begin{theorem}\label{th:law_v}
	Assume that $\mathcal{D} = \cF_{p}(\xi)$ for some $\xi,p>0$.
	There exists constants $C,c>0$ that only depend on $\Omega$ such that
	for all $\lambda \in [\lambda_{\rm min},\lambda_{\rm max}]$ and all $\epsilon \in (0,\frac{1}{2}]$,
	\begin{align*}
		\sup_{\theta^{\star} \in \mathcal{D} } \P \Big( W_2(\what{\mu}_{(\what{v}_{\lambda},\theta^{\star})},\nu^*_{\lambda})^2 \geq \epsilon  \Big) 
	\leq C \epsilon^{-\max(1,a)} \exp\left(-cN \epsilon^2 \epsilon^a \log(\epsilon)^{-2} \right) \,,
	\end{align*}
	where $a = \frac{1}{2} + \frac{1}{p}$.
\end{theorem}

Theorem~\ref{th:law_v} is proved in Section~\ref{sec:proof_th_law_v}.

\begin{theorem}\label{th:unif_law_v}
	Let $\mathcal{D}$ be $\cF_{p}(\xi)$ for some $\xi >0$ and $p > 0$.
	For all $\epsilon \in (0,\frac{1}{2}]$,
	$$
	\sup_{\theta^{\star} \in \mathcal{D} } \P \Big(\sup_{\lambda \in [\lambda_{\rm min},\lambda_{\rm max}]} W_2(\what{\mu}_{(\what{v}_{\lambda},\theta^{\star})},\nu^*_{\lambda})^2 \geq \epsilon  \Big) 
	\leq C \epsilon^{-\max(1,a)-1} N^{(1/p-1)_+}\exp\left(-cN \epsilon^2 \epsilon^a \log(\epsilon)^{-2} \right) \,,
	$$
	where $a = \frac{1}{2} + \frac{1}{p}$.
\end{theorem}
Theorem~\ref{th:unif_law_v} is deduced from Theorem~\ref{th:law_v} in Section~\ref{sec:proof_unif_v}.

\subsubsection{The norm of the subgradient}

Let us define
\begin{equation}\label{eq:def_kappa}
\kappa_*(\lambda) = \frac{\beta_*(\lambda)^2}{\lambda^2} \Big(1 + \delta - 2 s_*(\lambda) - \delta \frac{\sigma^2}{\tau_*(\lambda)^2}\Big) \,.
\end{equation}

\begin{theorem}\label{th:norm_v}
	There exists a constant $C,c>0$ such that for all $\lambda \in [\lambda_{\rm min},\lambda_{\rm max}]$ and all $\epsilon \in (0,1]$,
	$$
	\sup_{\theta^{\star} \in \mathcal{D} } \P \Big( \Big(\frac{1}{N}\|\what{v}_{\lambda}\|^2 - \kappa_*(\lambda)\Big)^2 \geq \epsilon \Big) \leq \frac{C}{\epsilon} e^{-cn\epsilon^2} \,.
	$$
\end{theorem}

Theorem~\ref{th:norm_v} is proved in Section~\ref{sec:proof_th_norm_v}. We deduce as before:

\begin{theorem}\label{th:unif_norm_v}
	Let $\mathcal{D}$ be either $\cF_0(s)$ for $s<s_{\rm max}(\delta)$ or $\cF_{p}(\xi)$ for some $\xi >0$ and $p > 0$.
	There exists constants $C,c>0$ such that for all $\epsilon \in (0,1]$,
	$$
	\sup_{\theta^{\star} \in \mathcal{D} } \P \Big(\sup_{\lambda \in [\lambda_{\rm min},\lambda_{\rm max}]} 
		\Big(\frac{1}{N}\|\what{v}_{\lambda}\|^2 - \kappa_*(\lambda)\Big)^2 \geq \epsilon
	\Big) \leq \frac{C}{\epsilon^2} N^q e^{-cn\epsilon^2 } \,,
	$$
	where $q = 0$ if\, $\mathcal{D} = \cF_0(s)$ and $q=(1/p - 1)_+$ if\, $\mathcal{D} = \cF_p(\xi)$.
\end{theorem}
Theorem~\ref{th:unif_norm_v} is deduced from Theorem~\ref{th:norm_v} in Section~\ref{sec:proof_unif_v}.

\subsubsection{Upper bound on the sparsity of the Lasso estimator}
\label{sec:UpperBoundSparsity}

Studying $\what{v}_{\lambda}$ allows to get an upper-bound on the $\ell_0$ norm of $\what{\theta}_{\lambda}$. Indeed if $\what{\theta}_{\lambda,i} \neq 0$ then $|\what{v}_{\lambda,i}| = 1$: therefore $\|\what{\theta}_{\lambda}\|_0 \leq \# \big\{ i \, \big| \, |\what{v}_{\lambda,i}| =1 \big\}$. For this reason, the following results will be used to prove Theorem~\ref{th:sparsity_uniform} in Appendix~\ref{sec:proof_sparsity}.

\begin{theorem}\label{th:key_1_v}
	There exists constants $C,c>0$ such that for all $\lambda \in [\lambda_{\rm min},\lambda_{\rm max}]$ and all $\epsilon \in (0,1]$,
	$$
	\sup_{\theta^{\star} \in \mathcal{D} } 
	\P \Big(
		\frac{1}{N}\# \big\{ i \, \big| \, |\what{v}_{\lambda,i}| \geq 1 - \epsilon \big\} \geq s_*(\lambda) + 2 (1 + \alpha_{\rm max}) \epsilon
	\Big) \leq \frac{C}{\epsilon^3} e^{-cn\epsilon^6} \,.
	$$
\end{theorem}

Theorem~\ref{th:key_1_v} is proved in Section~\ref{sec:proof_th_key_1_v}.

	\begin{theorem}\label{th:unif_1_v}
	Let $\mathcal{D}$ be either $\cF_0(s)$ for $s<s_{\rm max}(\delta)$ or $\cF_{p}(\xi)$ for some $\xi >0$ and $p > 0$.
	We have for all $\epsilon \in (0,1]$,
	$$ \sup_{\theta^{\star} \in \mathcal{D} } \P \Big(\exists \lambda \in [\lambda_{\rm min},\lambda_{\rm max}],\quad  \frac{1}{N} \# \big\{ i \, \big| \, |\what{v}_{\lambda,i}| =1 \big\}  \geq s_*(\lambda) + \epsilon  \Big) \leq \frac{C}{\epsilon^6} N^q e^{-cn\epsilon^6} \,,
	$$
	where $q = 0$ if\, $\mathcal{D} = \cF_0(s)$ and $q=(1/p - 1)_+$ if\, $\mathcal{D} = \cF_p(\xi)$.
\end{theorem}

Theorem~\ref{th:unif_1_v} is deduced from Theorem~\ref{th:key_1_v} in Section~\ref{sec:proof_unif_v}.

\subsection{Gordon's strategy for the subgradient}

\subsubsection{Application of Gordon's Theorem}

Let $g \sim \cN(0,\bbf{I}_N)$ and $h \sim \cN(0,\bbf{I}_n)$ be independent standard Gaussian vectors. We define:
\begin{align*}
	V_{\lambda}(v) 
	&= 
	\min_{w \in B}
	\left\{
		\frac{1}{2} \left(
			\sqrt{ \frac{1}{n}\|w\|^2  + \sigma^2 } \frac{\|h\|}{\sqrt{n}}
		-	\frac{1}{n} g^{\sT} w  + \frac{g'\sigma}{\sqrt{n}} \right)_{\!+}^2
		+ \frac{\lambda}{n} v^{\sT}(w+\theta^{\star})
		- \frac{\lambda}{n} |\theta^{\star}|
	\right\} \,.
\end{align*}

The following Proposition is the analog of Corollary~\ref{cor:gordon}.
\begin{proposition}
	Let $D \subset \big\{ v \in \R^N \, \big| \, \|v\|_{\infty} \leq 1 \big\}$ be a closed set.
	\begin{itemize}
		\item We have for all $t \in \R$
			$$
			\P\Big(\max_{v \in D} \mathcal{V}_{\lambda}(v) \geq t\Big) \leq 2 \P \Big(\max_{v \in D} V_{\lambda}(v) \geq t\Big) \,.
			$$
		\item If $D$ is convex, then we have for all $t \in \R$
			$$
			\P\Big(\max_{v \in D} \mathcal{V}_{\lambda}(v) \leq t\Big) \leq 2 \P \Big(\max_{v \in D} V_{\lambda}(v) \leq t\Big) \,.
			$$
	\end{itemize}
\end{proposition}
\begin{proof}
	Let $v \in \R^N$. By Proposition~\ref{prop:rockafellar} one can permute the min-max and obtain:
	\begin{align*}
		\mathcal{V}_{\lambda}(v)
		&=
		\min_{w \in B} \max_{u \in \R^n} \left\{
		\frac{1}{n}u^{\sT} X w - \frac{\sigma}{n}  u^{\sT}z - \frac{1}{2n}\|u\|^2 + \frac{\lambda}{n} v^{\sT}(\theta^{\star}+w) - \frac{\lambda}{n} |\theta^{\star}| \right\}
		\\
		&=
		\max_{u \in \R^n}\min_{w \in B}  
		\left\{
		\frac{1}{n}u^{\sT} X w - \frac{\sigma}{n}  u^{\sT}z - \frac{1}{2n}\|u\|^2 + \frac{\lambda}{n} v^{\sT}(\theta^{\star}+w) - \frac{\lambda}{n} |\theta^{\star}| \right\} \,.
	\end{align*}
	Let $D \subset \big\{ v \in \R^N \, \big| \, \|v\|_{\infty} \leq 1 \big\}$ be a closed set. 
	We can the apply Gordon's Theorem (Corollary~\ref{cor:gordon0}) in order to compare
	\begin{equation}\label{eq:gauss_p_v}
	\max_{(v,u)\in D \times \R^n} \min_{w \in B}
	\left\{
	\frac{1}{n}u^{\sT} X w - \frac{\sigma}{n}  u^{\sT}z - \frac{1}{2n}\|u\|^2 + \frac{\lambda}{n} v^{\sT}(\theta^{\star}+w) - \frac{\lambda}{n} |\theta^{\star}| \right\} \,,
	\end{equation}
	with
	\begin{align}
		&\max_{(v,u)\in D \times \R^n} \min_{w \in B}
		\left\{
		\sqrt{\frac{\|w\|^2}{n} + \sigma^2} \frac{h^{\sT}u}{n} - \frac{1}{n^{3/2}}\|u\| g^{\sT} w  + \frac{g'\sigma}{\sqrt{n}}  - \frac{1}{2n}\|u\|^2 + \frac{\lambda}{n} v^{\sT}(\theta^{\star}+w) - \frac{\lambda}{n} |\theta^{\star}| \right\} \label{eq:gauss_p_vg}
		\\
		&\quad =\max_{v \in D} \min_{w \in B} \max_{u \in \R^n}
		\left\{
		\sqrt{\frac{\|w\|^2}{n} + \sigma^2} \frac{h^{\sT}u}{n} - \frac{1}{n^{3/2}}\|u\| g^{\sT} w  + \frac{g'\sigma}{\sqrt{n}}  - \frac{1}{2n}\|u\|^2 + \frac{\lambda}{n} v^{\sT}(\theta^{\star}+w) - \frac{\lambda}{n} |\theta^{\star}| \right\} \,,
		\nonumber
	\end{align}
	which is equal to $\max_{v \in D} V_{\lambda}(v)$.
	Note that the maximums in~\eqref{eq:gauss_p_v} and~\eqref{eq:gauss_p_vg} are not defined on compact sets (since $D \times \R^n$ is not bounded). One has therefore to follow the same procedure than for Corollary~\ref{cor:gordon}, and show that there exists a compact set $\mathcal{K} \subset \R^n$ such that with high probability, the maximum over $u \in \R^n$ is achieved in $\mathcal{K}$. For the sake of brevity we do not provide a complete execution of this argument and refer to the proof of Corollary~\ref{cor:gordon}.
\end{proof}

\subsubsection{Study of Gordon's optimization problem}

In this section we study the optimization problem $\max_{\|v\|_{\infty} \leq 1} V_{\lambda}(v)$. 
Let us define
$$
\mathsf{v}_{\lambda} = -\alpha_*(\lambda)^{-1} \tau_*(\lambda)^{-1} \Big( \eta\big(\theta^{\star} + \tau_*(\lambda) g,\alpha_*(\lambda) \tau_*(\lambda)\big) - \theta^{\star} - \tau_*(\lambda) g\Big) \,.
$$
The goal of this section is to prove:
\begin{theorem}\label{th:gordon_aux_v}
	There exists constants $\gamma,c,C > 0$ that only depend on $\Omega$ such that for all $\theta^{\star} \in \mathcal{D}$, all $\lambda \in [\lambda_{\rm min}, \lambda_{\rm max}]$ and
	all $\epsilon \in (0,1]$
	$$
	\P \Big(
		\exists v\in B_{\infty}(0,1), \quad 
		\frac{1}{N} \| v - \mathsf{v}_{\lambda} \|^2 \geq \epsilon
		\quad \text{and} \quad
		V_{\lambda}(v) \geq \max\limits_{v \in \R^N} V_{\lambda}(v) - \gamma \epsilon 
	\Big)
	\leq \frac{C}{\epsilon} e^{-cn\epsilon^2} \,.
	$$
\end{theorem}

Recall that $w^*_{\lambda}$ is by Lemma~\ref{prop:w_star} the unique minimizer of $L_{\lambda}$ over $\R^N$.
\begin{lemma}\label{lem:opti_v_star}
	With probability at least $1 - 2 e^{-n/2}$ we have
	$$
	\min_{w \in \R^N} L_{\lambda}(w)
	=
	\max_{\|v\|_{\infty}\leq 1} V_{\lambda}(v) 
	$$
	and the vector
	\begin{equation}\label{eq:def_v_star}
	v_{\lambda}^* = -\lambda^{-1} \left(\sqrt{\frac{\|w_{\lambda}^*\|^2}{n} + \sigma^2}\ \frac{\|h\|}{\sqrt{n}}
	-	\frac{1}{n} g^{\sT} w_{\lambda}^*+ \frac{g'\sigma}{\sqrt{n}} \right)_{\!\!+}
	\left(\frac{\|h\|}{\sqrt{n}} \frac{w_{\lambda}^*}{\sqrt{\|w_{\lambda}^*\|^2/n + \sigma^2}} - g \right)
\end{equation}
	verifies $\|v^*_{\lambda}\|_{\infty} \leq 1$ and is a maximizer of $V_{\lambda}$.
\end{lemma}
\begin{proof}
	By Proposition~\ref{prop:rockafellar}, one can switch the min-max:
	\begin{align}
		\max_{\|v\|_{\infty} \leq 1} V_{\lambda}(v) 
		&= \min_{w \in B} \max_{\|v\|_{\infty} \leq 1}
		\left\{
			\frac{1}{2} \left(
				\sqrt{\frac{1}{n}\|w\|^2 + \sigma^2} \frac{\|h\|}{\sqrt{n}}
			-	\frac{1}{n} g^{\sT} w  + \frac{g'\sigma}{\sqrt{n}} \right)_{\!+}^2
			+ \frac{\lambda}{n} v^{\sT}(w+\theta^{\star})
			- \frac{\lambda}{n} |\theta^{\star}|
		\right\}
		\\
		&= \min_{w \in B} L_{\lambda}(w) \,.
	\end{align}
	Let us work on the event $\big\{\|h\| \leq 2 \sqrt{n}\big\} \cap \{ g' \leq \sqrt{n} \}$ which has probability at least $1- 2e^{-n/2}$.
	We have $\frac{\lambda}{n} (|w^*_{\lambda} + \theta^{\star}| - |\theta^{\star}|) \leq L_{\lambda}(w^*_{\lambda}) \leq L_{\lambda}(0) \leq 5 \sigma^2$. This gives $w^*_{\lambda} \in B$ and thus 
	$\max_{\|v\|_{\infty} \leq 1} V_{\lambda}(v) = \min_{w \in B} L_{\lambda}(w) = \min_{w \in \R^N} L_{\lambda}(w)$.
	\\

	The optimality condition of $w_{\lambda}^*$ gives that
	$$
	v^*_{\lambda} = -\lambda^{-1} \left(\sqrt{\frac{\|w_{\lambda}^*\|^2}{n} + \sigma^2}\ \frac{\|h\|}{\sqrt{n}}
	-	\frac{1}{n} g^{\sT} w^* + \frac{g'\sigma}{\sqrt{n}} \right)_{\!\!+}
	\left(\frac{\|h\|}{\sqrt{n}} \frac{w_{\lambda}^*}{\sqrt{\|w_{\lambda}^*\|^2/n + \sigma^2}} - g \right)
	\in \partial |\theta^{\star} + w_{\lambda}^*| \,.
	$$
	Therefore $v_{\lambda}^{*\sT} (w_{\lambda}^* + \theta^{\star})= |w_{\lambda}^* + \theta^{\star}|$. Using the optimality condition again we obtain
	\begin{align*}
		V_{\lambda}(v^{*}_{\lambda})
		&= \min_{w \in B} 
		\left\{
			\frac{1}{2} \left(
				\sqrt{\frac{1}{n}\|w\|^2 + \sigma^2} \frac{\|h\|}{\sqrt{n}}
			-	\frac{1}{n} g^{\sT} w + \frac{g'\sigma}{\sqrt{n}} \right)_{\!+}^2
			+ \frac{\lambda}{n} v_{\lambda}^{*\sT}(w+\theta^{\star})
			- \frac{\lambda}{n} |\theta^{\star}|
		\right\}
		\\
		&= 
		\frac{1}{2} \left(
			\sqrt{\frac{1}{n}\|w_{\lambda}^*\|^2 + \sigma^2} \frac{\|h\|}{\sqrt{n}}
		-	\frac{1}{n} g^{\sT} w_{\lambda}^*  + \frac{g'\sigma}{\sqrt{n}} \right)_{\!+}^2
		+ \frac{\lambda}{n} |w_{\lambda}^*+\theta^{\star}|
		- \frac{\lambda}{n} |\theta^{\star}|
=\min_{w \in \R^N} L_{\lambda}(w)
		= \max_{\|v\|_{\infty} \leq 1} V_{\lambda}(v) \,.
	\end{align*}
	Therefore $v_{\lambda}^*$ achieves the optimal value.
\end{proof}
\\

\begin{proposition}
	For all $\theta^{\star} \in \mathcal{D}$ and all $\lambda \in [\lambda_{\rm min},\lambda_{\rm max}]$ we have for all $\epsilon \in (0,1]$
	$$
	\P \Big(\frac{1}{N} \| v_{\lambda}^* - \mathsf{v}_{\lambda}\|^2 \geq \epsilon \Big) \leq \frac{C}{\epsilon} e^{-cn\epsilon} \,.
	$$
\end{proposition}
\begin{proof}
	By Theorem~\ref{th:gordon_aux} we have for all $\epsilon \in (0,1]$
	$$
	\P \Big(\frac{1}{N} \| w_{\lambda}^* - \mathsf{w}_{\lambda}\|^2 \geq \epsilon \Big) \leq \frac{C}{\epsilon} e^{-cn\epsilon} \,,
	$$
	so we deduce the result from the expression~\eqref{eq:def_v_star} of $v^*_{\lambda}$ and the concentration properties of $\mathsf{w}_{\lambda}$ (see Section~\ref{sec:concentration_w}).
\end{proof}
\\

By the same arguments used for proving Lemma~\ref{lem:opti_v_star} it is not difficult to prove:
\begin{lemma}\label{lem:b_star}
	The function
	$$
	\beta \geq 0 \mapsto \min_{w \in B} \ell_{\lambda}(w,\beta)
	$$
	(recall that $\ell_{\lambda}$ is defined by Equation~\ref{eq:def_l}) admits a unique maximizer $b^*_{\lambda}$ over $\R_{\geq 0}$ and
	$$
	b^*_{\lambda} = 
	\left(\sqrt{\frac{\|w_{\lambda}^*\|^2}{n} + \sigma^2}\ \frac{\|h\|}{\sqrt{n}}
	-	\frac{1}{n} g^{\sT} w_{\lambda}^* + \frac{g'\sigma}{\sqrt{n}} \right)_{\!\!+} \,.
	$$
	Moreover, for all $\epsilon \in (0,1]$ we have $\P\Big(|b^*_{\lambda} - \beta_*(\lambda)| > \epsilon \Big) \leq \frac{C}{\epsilon} e^{-cn\epsilon}$.
\end{lemma}

\subsubsection{Proof of Theorem~\ref{th:gordon_aux_v}}
	Let $v \in \R^N$ such that $\|v\|_{\infty} \leq 1$.
	We have by Proposition~\ref{prop:rockafellar}
	\begin{align*}
		V_{\lambda}(v)
		&=
		\min_{w \in B}
		\max_{\beta \geq 0}
		\left\{
			\beta \left(
				\sqrt{\frac{1}{n}\|w\|^2  + \sigma^2 } \frac{\|h\|}{\sqrt{n}}
			-	\frac{1}{n} g^{\sT} w  + \frac{g'\sigma}{\sqrt{n}} \right) - \frac{\beta^2}{2}
			+ \frac{\lambda}{n} v^{\sT}(w+\theta^{\star})
			- \frac{\lambda}{n} |\theta^{\star}|
		\right\}
		\\
		&=
		\max_{\beta \geq 0}
		\min_{w \in B}
		\left\{
			\beta \left(
				\sqrt{\frac{1}{n}\|w\|^2  + \sigma^2 } \frac{\|h\|}{\sqrt{n}}
			-	\frac{1}{n} g^{\sT} w  + \frac{g'\sigma}{\sqrt{n}} \right) - \frac{\beta^2}{2}
			+ \frac{\lambda}{n} v^{\sT}(w+\theta^{\star})
			- \frac{\lambda}{n} |\theta^{\star}|
		\right\}
		\\
		&=
		\max_{\beta \geq 0}
		\min_{0 \leq r \leq R}
		\left\{
			\beta\sqrt{r^2  + \sigma^2 } \frac{\|h\|}{\sqrt{n}}
			-\frac{r}{\sqrt{n}} \|\beta g - \lambda v\|
			+ \beta\frac{g'\sigma}{\sqrt{n}}  
			- \frac{\beta^2}{2}
			+ \frac{\lambda}{n} v^{\sT}\theta^{\star}
			- \frac{\lambda}{n} |\theta^{\star}|
		\right\} \,,
	\end{align*}
because the minimization over the direction of $w$ is easy to perform.
	Let us define for $\kappa >0$
	$$
	D_{\kappa} = \Big\{ v \in B_{\infty}(0,1), \quad \frac{1}{N} \|v - v_{\lambda}^*\|^2 \leq \kappa^2 \quad \text{and} \quad V_{\lambda}(v) \geq \max_{\|v'\|_{\infty} \leq 1} V_{\lambda}(v') - \frac{1}{8}\kappa^2 \Big\} \,.
	$$
	By concavity of $V_{\lambda}$, $D_{\kappa}$ is convex. 
	\begin{proposition}\label{prop:v_tilde}
		There exists a constant $\kappa > 0$ such that with probability at least $1 - C e^{-cn}$ we have $\forall v \in D_{\kappa}, V_{\lambda}(v) = \widetilde{V}_{\lambda}(v)$ where
		$$
		\widetilde{V}_{\lambda}: v \mapsto
		\min_{w \in \R^N}
		\left\{
			\max_{\beta \in [\beta_*-\kappa,\beta_*+\kappa]}
			\left\{
				\beta \left(
					\sqrt{ \frac{1}{n}\|w\|^2  + \sigma^2 } \frac{\|h\|}{\sqrt{n}}
				-	\frac{1}{n} g^{\sT} w +  \frac{g'\sigma}{\sqrt{n}} \right) - \frac{\beta^2}{2}
			\right\}
			+ \frac{\lambda}{n} v^{\sT}(w+\theta^{\star})
			- \frac{\lambda}{n} |\theta^{\star}|
		\right\} \,.
		$$
	\end{proposition}

	In order to prove Proposition~\ref{prop:v_tilde}, we start with a Lemma:
		\begin{lemma}\label{lem:min_f}
		For all $v \in D_{\kappa}$, the function
		$$
		f_v:\beta \mapsto	
		\min_{w \in B}
		\left\{
			\beta \left(
				\sqrt{ \frac{1}{n}\|w\|^2  + \sigma^2 } \frac{\|h\|}{\sqrt{n}}
			-	\frac{1}{n} g^{\sT} w  + \frac{g'\sigma}{\sqrt{n}} \right) - \frac{\beta^2}{2}
			+ \frac{\lambda}{n} v^{\sT}(w+\theta^{\star})
			- \frac{\lambda}{n} |\theta^{\star}|
		\right\}
		$$
		admits a unique maximizer $b_{\lambda}(v)$ on $[0,+\infty)$ and one has $|b_{\lambda}(v) - b_{\lambda}^*| \leq \kappa/2$.
	\end{lemma}
	\begin{proof}
		Let $v \in D_{\kappa}$. $f_v$ is $1$-strongly concave so it admits a unique maximizer $b_{\lambda}(v)$ on $\R_{\geq 0}$. We have
		$$
		f_v(b_{\lambda}(v)) = \max_{\beta\geq 0} f_v(\beta) = V_{\lambda}(v) \geq \max_{\|v'\|_{\infty} \leq 1} V_{\lambda}(v') - \frac{1}{8} \kappa^2 \,,
		$$
		because $v \in D_{\kappa}$. Notice now that $f_v(b_{\lambda}(v)) \leq \min_{w \in B} \ell_{\lambda}(w,b_{\lambda}(v))$ because $v^{\sT}(w+\theta^{\star}) \leq |w + \theta^{\star}|$.
		Permuting the min-max (using Proposition~\ref{prop:rockafellar}), we have $\displaystyle \max_{\|v'\|_{\infty} \leq 1} V_{\lambda}(v') = \max_{\beta \geq 0} \min_{w \in B} \ell_{\lambda}(w,\beta)$, where we recall that $\ell_{\lambda}$ is defined by~\eqref{eq:def_l}. We get
		$$
		\min_{w \in B} \ell_{\lambda}(w,b_{\lambda}(v)) \geq \max_{\beta \geq 0}\min_{w \in B} \ell_{\lambda}(w,\beta) - \frac{1}{8} \kappa^2 \,.
		$$
		The function $\beta \mapsto \min_{w \in B}\ell_{\lambda}(w,\beta)$ is $1$-strongly concave and maximized (by Lemma~\ref{lem:b_star} above) at $b^*_{\lambda}$, hence $|b_{\lambda}(v) - b_{\lambda}^*| \leq \kappa/2$.
	\end{proof}
	\\

	\begin{lemma}\label{lem:kappa}
		There exist constants $K,\kappa > 0$ such that with probability at least $1 - Ce^{-cn}$ the following happens.
		For all $\beta \geq 0$, $v \in \R^N$ such that $|\beta - b_{\lambda}^*| \leq 2 \kappa$ and $\|v - v_{\lambda}^*\| \leq \sqrt{N} \kappa$ the minimum over $\R^N$ of
		$$
		w \mapsto
		\beta \left(
			\sqrt{ \frac{1}{n}\|w\|^2  + \sigma^2 } \frac{\|h\|}{\sqrt{n}}
		-	\frac{1}{n} g^{\sT} w +\frac{g'\sigma}{\sqrt{n}} \right) - \frac{\beta^2}{2}
		+ \frac{\lambda}{n} v^{\sT}(w+\theta^{\star})
		- \frac{\lambda}{n} |\theta^{\star}|
		$$
		is achieved on $B(0,\sqrt{N}K)$.
	\end{lemma}
	\begin{proof}
		The minimization with respect to the direction of $w$ is easy to perform: $w$ has to be a non-negative multiple of $\beta g -\lambda v$. It remains thus to minimizes with respect to the norm of $w$. We have to show that under the conditions of the lemma, the minimum of
		\begin{equation}\label{eq:min_r}
		r \geq 0 \mapsto \beta\sqrt{r^2  + \sigma^2 } 
		-\frac{r}{\|h\|} \|\beta g - \lambda v\|
	\end{equation}
		is achieved for $r$ smaller than some constant.
		By Theorem~\ref{th:gordon_aux} and Lemma~\ref{lem:b_star}
		there exists a constant $R>0$ (for instance $R = \tau_{\rm max} + 1$) such that
		the event
		\begin{equation}\label{eq:event_ras}
			\Big\{ \|w_{\lambda}^*\|^2 \leq n R^2 \Big\} \bigcap 
			\Big\{ b^*_{\lambda} \geq \beta_{\rm min} / 2 \Big\} \bigcap
			\Big\{ \|h\| \geq \sqrt{n} / 2 \Big\}
		\end{equation}
		has probability at least $1 - Ce^{-cn}$.
		Let us define the constants $a = \frac{R}{\sqrt{R^2 + \sigma^2}} < 1$ and
		$$
		\kappa = \min \left( \frac{\sqrt{\delta} \beta_{\rm min}^2 (1-a)}{256 \lambda_{\rm max}}  , \quad \beta_{\rm min}/8, \quad
		\frac{(1-a)\sqrt{\delta} \beta_{\rm min}}{16\lambda_{\rm max}} \right).
		$$
		Let us now work on the event~\eqref{eq:event_ras}. Let $v \in \R^N$ and $\beta \geq 0$ such that $\|v-v^*_{\lambda}\| \leq \sqrt{N}\kappa$ and $|\beta - b^*_{\lambda}|\leq 2 \kappa$. We have
		\begin{align*}
			\| g - \frac{\lambda}{\beta} v\|
			\leq
			\| g - \frac{\lambda}{b^*_{\lambda}} v_{\lambda}^*\|
			+
			\| \frac{\lambda}{\beta} v - \frac{\lambda}{b^*_{\lambda}} v_{\lambda}^*\| \,.
		\end{align*}
		Compute
		$$
		\| g - \frac{\lambda}{b^*_{\lambda}} v_{\lambda}^*\|
		= \|h\| \frac{\|w_{\lambda}^*\| / \sqrt{n}}{\sqrt{\frac{\|w_{\lambda}^*\|^2}{n} + \sigma^2}}
		\leq \|h\| a \,,
		$$
		with probability at least $1-Ce^{-cn}$.
		Now
		\begin{align*}
			\| \frac{1}{\beta} v - \frac{1}{b^*_{\lambda}} v_{\lambda}^*\|
			&\leq
			\frac{1}{\beta^*}\| v -  v_{\lambda}^*\|
			+
			\| \frac{1}{\beta} v - \frac{1}{b^*_{\lambda}} v\|
			\leq
			\frac{2}{\beta_{\rm min}}\| v -  v_{\lambda}^*\|
			+
			\frac{\sqrt{N}}{\min(\beta,b^*_{\lambda})^2}|\beta - b^*_{\lambda}|
			\\
			&\leq
			\frac{2}{\beta_{\rm min}}\| v -  v^*_{\lambda}\|
			+
			\frac{16\sqrt{N}}{\beta_{\rm  min}^2}	|\beta - b^*_{\lambda}| \leq \frac{1 - a}{4} \sqrt{n} \,.
		\end{align*}
		Putting all together:
		$$
		\frac{1}{\|h\|}\| g - \frac{\lambda}{\beta} v\| \leq a + \frac{1-a}{2} = \frac{1+a}{2} < 1 \,.
		$$
		This gives that the minimum of~\eqref{eq:min_r} is achieved for $r \leq \frac{\sigma}{\sqrt{1 - ((1+a)/2)^2}}$. One can thus chose $K = \frac{\delta \sigma}{\sqrt{1 - ((1+a)/2)^2}}$.
	\end{proof}
	\\

	\begin{proof}[of Proposition~\ref{prop:v_tilde}]
	Let us now fix a constant $\kappa \in (0,\beta_{\rm min}/2)$ that verify the statement of Lemma~\ref{lem:kappa}. 
	Let us work on the intersection of the event $\{ |b^*_{\lambda} - \beta_*(\lambda)|\leq \kappa /2 \}$ with the event of Lemma~\ref{lem:kappa}. This intersection has by Lemma~\ref{lem:b_star} and Lemma~\ref{lem:kappa} probability at least $1-Ce^{-cn}$.

	Let $v \in D_{\kappa}$. By Lemma~\ref{lem:min_f} the unique maximizer $b_{\lambda}(v)$ of $f_v$ verify $|b_{\lambda}(v) - b_{\lambda}^*| \leq \kappa/2$ and therefore $|b_{\lambda}(v) - \beta_*| \leq \kappa$. Consequently
	$$
		V_{\lambda}(v)
		=
		\max_{\beta \in [\beta_*-\kappa,\beta_*+\kappa]}
		\min_{w \in B}
		\left\{
			\beta \left(
				\sqrt{ \frac{1}{n}\|w\|^2  + \sigma^2 } \frac{\|h\|}{\sqrt{n}}
			-	\frac{1}{n} g^{\sT} w  + \frac{g'\sigma}{\sqrt{n}} \right) - \frac{\beta^2}{2}
			+ \frac{\lambda}{n} v^{\sT}(w+\theta^{\star})
			- \frac{\lambda}{n} |\theta^{\star}|
		\right\} \,.
	$$
	Now, for $\beta \in [\beta_*-\kappa,\beta_* +\kappa]$, we have $|\beta - b_{\lambda}^*| \leq 2 \kappa$. Since we are working on the event of Lemma~\ref{lem:kappa}, we obtain
	\begin{align*}
		V_{\lambda}(v)
		&=
		\max_{\beta \in [\beta_*-\kappa,\beta_*+\kappa]}
		\min_{w \in \R^N}
		\left\{
			\beta \left(
				\sqrt{\frac{1}{n}\|w\|^2  + \sigma^2 } \frac{\|h\|}{\sqrt{n}}
			-	\frac{1}{n} g^{\sT} w  + \frac{g'\sigma}{\sqrt{n}} \right) - \frac{\beta^2}{2}
			+ \frac{\lambda}{n} v^{\sT}(w+\theta^{\star})
			- \frac{\lambda}{n} |\theta^{\star}|
		\right\} \,,
	\end{align*}
	and Proposition~\ref{prop:v_tilde} follows from the permutation of the $\min-\max$ using Proposition~\ref{prop:rockafellar}.
	\end{proof}
	\\

	\begin{lemma}\label{lem:v_tilde_concave}
		There exists a constant $C,c,\gamma >0$ such that
		$\widetilde{V}_{\lambda}$ is $\gamma/N$-strongly concave, with probability at least $1- Ce^{-cn}$.
	\end{lemma}
	\begin{proof}
		The function
		$$
		f^*:
		v\in \R^N \mapsto 
		\max_{w \in \R^N}
		\left\{
			v^{\sT}w - 
			\frac{n}{\lambda}
			\max_{\beta \in [\beta_*-\kappa,\beta_*+\kappa]}
			\left\{
				\beta \left(
					\sqrt{ \frac{1}{n}\|w\|^2  + \sigma^2 } \frac{\|h\|}{\sqrt{n}}
				-	\frac{1}{n} g^{\sT} w  + \frac{g'\sigma}{\sqrt{n}} \right) - \frac{\beta^2}{2}
			\right\}
		\right\}
		$$
		is the convex conjugate of the convex function 	
	\begin{align*}
		f:w\in \R^N \mapsto
		&\frac{n}{\lambda}
		\max_{\beta \in [\beta_*-\kappa,\beta_*+\kappa]}
		\left\{
			\beta \left(
				\sqrt{ \frac{1}{n}\|w\|^2  + \sigma^2 } \frac{\|h\|}{\sqrt{n}}
			-	\frac{1}{n} g^{\sT} w  + \frac{g'\sigma}{\sqrt{n}} \right) - \frac{\beta^2}{2}
		\right\}
		\\
		&= \frac{n}{\lambda} \varphi\left(
			\sqrt{ \frac{1}{n}\|w\|^2  + \sigma^2 } \frac{\|h\|}{\sqrt{n}}
		-	\frac{1}{n} g^{\sT} w  + \frac{g'\sigma}{\sqrt{n}} \right) \,,
	\end{align*}
		where $\varphi$ is the $\cC^1$ function
		$$
		\varphi(x) =
		\begin{cases}
			\frac{1}{2} x^2 & \text{if} \quad x \in [\beta_*-\kappa,\beta_*+\kappa]\,, \\
			(\beta_* - \kappa) x - \frac{1}{2} (\beta_* -\kappa)^2 & \text{if} \quad x \leq \beta_*-\kappa\,, \\
			(\beta_* + \kappa) x - \frac{1}{2} (\beta_* +\kappa)^2 & \text{if} \quad x \geq \beta_*+\kappa\,. \\
		\end{cases}
		$$
		$f$ is a proper closed convex function (because $f$ is convex and its domain is $\R^N$), therefore its convex conjugate $f^*$ is also a proper closed convex function. The Fenchel-Moreau Theorem gives then that $f=f^{**}$.
		Let us compute the gradient of $f$ for $w \in \R^N$
		\begin{align*}
			\nabla f(w) 
			&= \frac{n}{\lambda}
			\varphi'\left(
				\sqrt{ \frac{1}{n}\|w\|^2  + \sigma^2 } \frac{\|h\|}{\sqrt{n}}
			-	\frac{1}{n} g^{\sT} w + \frac{g'\sigma}{\sqrt{n}} \right)
			\left(
				\frac{\|h\|}{\sqrt{n}} \frac{w/n}{\sqrt{\frac{\|w\|^2}{n} + \sigma^2}} - \frac{g}{n}
			\right) \,.
		\end{align*}
		It is not difficult to verify that there exists a constant $L$ such that $\nabla f$ is $L$-Lipschitz on $\R^N$, with probability at least $1 - Ce^{-cn}$.
		$f = f^{**}$ is therefore $1/L$-strongly smooth (see Definition~\ref{def:smooth}).
		By Proposition~\ref{prop:duality} this gives that $f^*$ is $1/L$ strongly convex. One deduces then that $\widetilde{V}_{\lambda}$ is $\gamma / N$- strongly concave with $\gamma = \lambda / (L\delta)$.
	\end{proof}
	\\

	Let $0< \gamma < 1/2$ be a constant that verify the statement of Lemma~\ref{lem:v_tilde_concave} and let $\kappa > 0$ be a constant given by Proposition~\ref{prop:v_tilde}.
	Notice that it suffices to prove Theorem~\ref{th:gordon_aux_v} for $\epsilon$ small enough and let $\epsilon \in (0,\kappa^2)$. 
	\begin{align}
		&\P \Big(
		\exists v\in B_{\infty}(0,1), \quad 
		\frac{1}{N} \| v - \mathsf{v}_{\lambda} \|^2 > \epsilon
		\quad \text{and} \quad
		V_{\lambda}(v) \geq \max\limits_{\|v'\|_{\infty} \leq 1} V_{\lambda}(v') - \frac{1}{4}\gamma \epsilon 
	\Big)\nonumber
	\\
	&\qquad \leq 
	\P \Big(
		\exists v\in B_{\infty}(0,1), \quad 
		\frac{1}{N} \| v - v_{\lambda}^* \|^2 > \frac{\epsilon}{2}
		\quad \text{and} \quad
		V_{\lambda}(v) \geq \max\limits_{\|v'\|_{\infty} \leq 1} V_{\lambda}(v') - \frac{1}{4}\gamma \epsilon 
	\Big)
	+ \frac{C}{\epsilon} e^{-cn\epsilon^2} \nonumber
	\\
	&\qquad \leq 
	\P \Big(
		\exists v\in D_{\kappa}, \quad 
		\frac{1}{N} \| v - v_{\lambda}^* \|^2 > \frac{\epsilon}{2}
		\quad \text{and} \quad
		V_{\lambda}(v) \geq \max\limits_{v' \in D_{\kappa}} V_{\lambda}(v') - \frac{1}{4}\gamma \epsilon 
	\Big)
	+ \frac{C}{\epsilon} e^{-cn\epsilon^2} \,, \label{eq:last_proba}
\end{align}
because, if there exists $v \in B_{\infty}(0,1)$ such that
$\frac{1}{N} \| v - v_{\lambda}^* \|^2 > \frac{\epsilon}{2}$ and $V_{\lambda}(v) \geq \max\limits_{\|v'\|_{\infty} \leq 1} V_{\lambda}(v') - \frac{1}{4}\gamma \epsilon$, we can construct $\widetilde{v} \in D_{\kappa}$ that verifies the same conditions. Indeed:
\begin{itemize}
	\item if $\frac{1}{N} \| v - v_{\lambda}^* \|^2 \leq \kappa^2$, one simply take $\widetilde{v} = v$.
	\item otherwise, $\widetilde{v} = v^*_{\lambda} + \kappa (v - v^*_{\lambda}) / \|v - v^*_{\lambda} \|$ is in $D_{\kappa}$ and by concavity $V_{\lambda}(\widetilde{v}) \geq V_{\lambda}(v)$.
\end{itemize}
Since with probability at least $1 - Ce^{-cn}$ we have $V_{\lambda}(v) = \widetilde{V}_{\lambda}(v)$ for all $v \in D_{\kappa}$ and $\widetilde{V}_{\lambda}$ is $\gamma / N$-strongly concave, the probability in~\eqref{eq:last_proba} above is less that $C e^{-cn}$.

\subsection{Proofs of the main results about the subgradient}

Let us start with the analog of Proposition~\ref{prop:apply_gordon} for the costs functions $\mathcal{V}_{\lambda}$ and $V_{\lambda}$:
\begin{proposition} \label{prop:apply_gordon_v}
	There exists constants $c,C >0$ that only depend on $\Omega$ such that
	for all closed set $D \subset \R^N$ and for all $\epsilon \in (0,1]$,
	$$
	\P \left(\max_{v \in D} \mathcal{V}_{\lambda}(v) \geq \max_{\|v\|_{\infty} \leq 1} \mathcal{V}_{\lambda}(v) - \epsilon \right) 
	\leq
	2 \P \left(\max_{v \in D} V_{\lambda}(v) \geq \max_{\|v\|_{\infty} \leq 1} V_{\lambda}(v) - 3 \epsilon \right)  + \frac{C}{\epsilon} e^{-cn \epsilon^2} \,.
	$$
\end{proposition}

The proof of Proposition~\ref{prop:apply_gordon_v} is omitted for the sake of brevity, and because it follows from the exact same arguments than Proposition~\ref{prop:apply_gordon}.

\subsubsection{The norm of \texorpdfstring{$\what{v}_{\lambda}$}{v}: proof of Theorem~\ref{th:norm_v}}\label{sec:proof_th_norm_v}

\begin{lemma}\label{lem:strong_concave_v_risk}
	There exists constants $\gamma,c,C > 0$ that depend only on $\Omega$, such that for all $\epsilon \in (0,1]$ we have
	$$
	\P\left(
		\max_{v \in D_{\epsilon}} V_{\lambda}(v)
		\geq \max_{\|v\|_{\infty} \leq 1}V_{\lambda}(v) - 3 \gamma \epsilon
	\right)
	\leq 
	\frac{C}{\epsilon} e^{-cn\epsilon^2} \,,
	$$
	where $ D_{\epsilon} = \left\{ v \in B_{\infty}(0,1) \, \middle| \, \big( \|v\| - \sqrt{N \kappa_*(\lambda)}\big)^2 \geq \epsilon \right\}$ and $\kappa_*(\lambda)$ is defined by~\eqref{eq:def_kappa}.
\end{lemma}
\begin{proof}
	Similarly to Proposition~\ref{prop:concentration_cost_w_star} it is not difficult to prove that for all $\epsilon \in (0,1]$,
	$$
	\P \Big( \Big| \frac{1}{N} \|\mathsf{v}_{\lambda}\|^2 - \kappa_*(\lambda)^2 \Big| > \epsilon \Big) \leq C e^{-cN\epsilon^2} \,,
	$$
	for some constants $c,C >0$.
	By Theorem~\ref{th:gordon_aux_v} there exists constants $\gamma,c,C>0$ such that for all $\epsilon \in (0,1]$ the event
	\begin{equation}\label{eq:event_good_v_r}
		\Big\{ \forall v \in B_{\infty}(0,1), \ V_{\lambda}(v) \geq \max_{\|v\|_{\infty} \leq 1} V_{\lambda}(v) - 3 \gamma \epsilon  \implies \frac{1}{N} \|v-\mathsf{v}_{\lambda}\|^2 \leq \frac{\epsilon}{5} \Big\}
		\bigcap
		\Big\{ \big(\|\mathsf{v}_{\lambda}\| - \sqrt{N\kappa_*(\lambda)}\big)^2 \leq N \frac{\epsilon}{4} \Big\}
	\end{equation}
	has probability at least $\frac{C}{\epsilon} e^{-cn\epsilon^2}$.
	On the event~\eqref{eq:event_good_v_r}, we have for all $v \in D_{\epsilon}$:
	$$
	\frac{1}{N} \| v - \mathsf{v}_{\lambda} \|^2 \geq 
	\frac{1}{N} \big( \|v\| - \|\mathsf{v}_{\lambda}\| \big)^2
	\geq
	\frac{1}{N} \big( \sqrt{N\epsilon} - \frac{1}{2} \sqrt{N\epsilon} \big)^2
	\geq \frac{\epsilon}{4} \,.
	$$
	This gives that on the event~\eqref{eq:event_good_v_r}, for all $v \in D_{\epsilon}$, $V_{\lambda}(v) < \max\limits_{\|v'\|_{\infty} \leq 1} V_{\lambda}(v') - 3 \gamma \epsilon$. The intersection of~\eqref{eq:event_good_v_r} with the event $\big\{ \max\limits_{v \in D_{\epsilon} } V_{\lambda}(v) \geq \max\limits_{\|v\|_{\infty} \geq 1}V_{\lambda}(v) - 3 \gamma \epsilon \big\}$ is therefore empty: the lemma is proved.
\end{proof}
\\

\begin{proof}[of Theorem~\ref{th:norm_v}]
Let $\gamma >0$ be a constant that verify the statement of Lemma~\ref{lem:strong_concave_v_risk}.
Let $\epsilon \in (0,1]$ and define 
$$ D_{\epsilon} = \left\{ v \in B_{\infty}(0,1) \, \middle| \, \big(\|v\| - \sqrt{N \kappa_*(\lambda)}\big)^2 \geq N \epsilon \right\} \,.
$$
$D_{\epsilon}$ is a closed set.
\begin{align*}
	\P \Big(\exists v \in B_{\infty}(0,1)&, \ \big|\frac{1}{N}\|v\|^2 - \kappa_*(\lambda)\big| \geq \epsilon \ \text{and} \ \mathcal{V}_{\lambda}(v) \geq \max \mathcal{V}_{\lambda} - \gamma \epsilon \Big)
	= \P\Big(\max_{v \in D_{\epsilon}} \mathcal{V}_{\lambda}(v) \geq \max_{\|v\|_{\infty} \leq 1}\mathcal{V}_{\lambda}(v) - \gamma \epsilon \Big)
	\\
	&\leq 2 \P\Big(\max_{v \in D_{\epsilon}} V_{\lambda}(v) \geq \max_{\|v\|_{\infty} \leq 1}V_{\lambda}(v) - 3 \gamma \epsilon \Big) + C e^{-cn\epsilon}
	\leq \frac{C}{\epsilon} e^{-cn\epsilon^2} \,,
\end{align*}
where we used successively Proposition~\ref{prop:apply_gordon_v} and Lemma~\ref{lem:strong_concave_v_risk}.
\end{proof}

\subsubsection{The empirical law of \texorpdfstring{$\what{v}_{\lambda}$}{v}: proof of Theorem~\ref{th:law_v}}\label{sec:proof_th_law_v}

Theorem~\ref{th:law_v} follows now from Proposition~\ref{prop:apply_gordon_v} and the following Lemma.

\begin{lemma}\label{lem:strong_concave_v_law}
	There exists constants $\gamma,c,C > 0$ that depend only on $\Omega$, such that for all $\epsilon \in (0,\frac{1}{2}]$ we have
	$$
	\P\left(
		\max_{v \in D_{\epsilon}} V_{\lambda}(v)
		\geq \max_{\|v\|_{\infty} \leq 1}V_{\lambda}(v) - 3 \gamma \epsilon
	\right)
	\leq 
	\frac{C}{\epsilon} e^{-cn\epsilon^2 \log(\epsilon)^{-2}} \,,
	$$
	where $ D_{\epsilon} = \left\{ v \in B_{\infty}(0,1) \, \middle| \, W_2(\what{\mu}_{(v,\theta^{\star})},\nu_{\lambda}^*)^2 \geq \epsilon \right\}$.
\end{lemma}
\begin{proof}
	By Theorem~\ref{th:gordon_aux_v} and Proposition~\ref{prop:concentration_empirical_distribution} there exists constants $\gamma,c,C>0$ such that for all $\epsilon \in (0,\frac{1}{2}]$ the event
	\begin{equation}\label{eq:event_good_v}
		\Big\{ \forall v \in B_{\infty}(0,1), \ V_{\lambda}(v) \geq \max_{\|v\|_{\infty} \leq 1} V_{\lambda}(v) - 3 \gamma \epsilon  \implies \frac{1}{N} \|v-\mathsf{v}_{\lambda}\|^2 \leq \frac{\epsilon}{5} \Big\}
		\bigcap
		\Big\{ W_2\big(\nu_{\lambda}^*,\what{\mu}_{(\mathsf{v}_{\lambda},\theta^{\star})}\big)^2 \leq \frac{\epsilon}{4} \Big\}
	\end{equation}
	has probability at least
	$$
	1 - C \epsilon^{-1} \exp\big( -cn\epsilon^2\big) - C \epsilon^{-a} \exp\left(-cN \epsilon^2 \epsilon^a \log(\epsilon)^{-2} \right)
	\geq
	1 - C \epsilon^{-\max(1,a)} \exp\left(-cN \epsilon^2 \epsilon^a \log(\epsilon)^{-2} \right).
	$$
	On the event~\eqref{eq:event_good_v}, we have for all $v \in D_{\epsilon}$:
	$$
	\frac{1}{N} \| v - \mathsf{v}_{\lambda} \|^2 
	\geq W_2\big(\what{\mu}_{(v,\theta^{\star})},\what{\mu}_{(\mathsf{v}_{\lambda},\theta^{\star})}\big)^2
	\geq 
	\Big(
		W_2(\what{\mu}_{(v,\theta^{\star})},\nu_{\lambda}^*)
		-W_2(\nu_{\lambda}^*,\what{\mu}_{(\mathsf{v}_{\lambda},\theta^{\star})})
	\Big)^2
	\geq \frac{\epsilon}{4} \,.
	$$
	This gives that on the event~\eqref{eq:event_good_v}, for all $v \in D_{\epsilon}$, $V_{\lambda}(v) < \max\limits_{\|v'\|_{\infty} \leq 1} V_{\lambda}(v') - 3 \gamma \epsilon$. The intersection of~\eqref{eq:event_good_v} with the event $\big\{ \max\limits_{v \in D_{\epsilon} } V_{\lambda}(v) \geq \max\limits_{\|v\|_{\infty} \geq 1}V_{\lambda}(v) - 3 \gamma \epsilon \big\}$ is therefore empty: the lemma is proved.
\end{proof}
\\

\begin{proof}[of Theorem~\ref{th:law_v}]
Let $\gamma >0$ be a constant that verify the statement of Lemma~\ref{lem:strong_concave_v_law}.
Let $\epsilon \in (0,\frac{1}{2}]$ and define 
$$ D_{\epsilon} = \left\{ v \in B_{\infty}(0,1) \, \middle| \,  W_2(\what{\mu}_v,\nu_{\lambda}^*)^2 \geq \epsilon \right\}
$$
$D_{\epsilon}$ is a closed set.
\begin{align*}
	\P \Big(\exists v \in B_{\infty}(0,1),& \ W_2(\what{\mu}_{v},\nu^*_{\lambda})^2 \geq \epsilon \ \text{and} \ \mathcal{V}_{\lambda}(v) \geq \max \mathcal{V}_{\lambda} - \gamma \epsilon \Big)
	= \P\Big(\max_{v \in D_{\epsilon}} \mathcal{V}_{\lambda}(v) \geq \max_{\|v\|_{\infty} \leq 1}\mathcal{V}_{\lambda}(v) - \gamma \epsilon \Big)
	\\
	&\leq 2 \P\Big(\max_{v \in D_{\epsilon}} V_{\lambda}(v) \geq \max_{\|v\|_{\infty} \leq 1}V_{\lambda}(v) - 3 \gamma \epsilon \Big) + C e^{-cn\epsilon}
	\leq \frac{C}{\epsilon} e^{-cn\epsilon^2 \log(\epsilon)^{-2}} \,,
\end{align*}
where we used successively Proposition~\ref{prop:apply_gordon_v} and Lemma~\ref{lem:strong_concave_v_law}.
\end{proof}

\subsubsection{Proof of Theorem~\ref{th:key_1_v}}\label{sec:proof_th_key_1_v}

\begin{lemma}\label{lem:strong_concave_v_1}
	There exists constants $\gamma,c,C > 0$ that depend only on $\Omega$, such that for all $\epsilon \in (0,1]$ we have
	$$
	\P\left(
		\max_{v \in D_{\epsilon}} V_{\lambda}(v)
		\geq \max_{\|v\|_{\infty} \leq 1}V_{\lambda}(v) - 3 \gamma \epsilon^3
	\right)
	\leq 
	\frac{C}{\epsilon^3} e^{-cn\epsilon^6} \,,
	$$
	where $ D_{\epsilon} = \Big\{ v \in B_{\infty}(0,1) \, \Big|  \frac{1}{N} \# \big\{ i \, \big| \, |v_i| \geq 1 - \epsilon \big\} > s_*(\lambda) + 2 (1+\alpha_{\rm max}) \epsilon \Big\}$.
\end{lemma}
\begin{proof}
	Let $\epsilon \in (0,1]$ and define
	$$
	s_{\epsilon} = \frac{1}{N} \# \Big\{ i \in \{1,\dots,N\} \, \Big| \, |\mathsf{v}_{\lambda,i}| \geq 1- 2 \epsilon \Big\} \,.
	$$
	$s_{\epsilon}$ is the mean of independent Bernoulli random variables. By Hoeffding's inequality we have
	$$
	\P \Big(s_{\epsilon} \leq \P\big(|\tau_*^{-1}\Theta +  Z| \geq \alpha_* - 2 \alpha_* \epsilon \big) + \epsilon \Big) \geq 1 - e^{-2 N\epsilon^2} \,.
	$$
	Compute
	\begin{align*}
		\P\big(|\tau_*^{-1}\Theta + Z| \geq \alpha_* - 2\alpha_* \epsilon \big)
		&=
		\E \left[\Phi\Big(\frac{\Theta}{\tau_*(\lambda)} - \alpha_*(\lambda)+ 2 \alpha_*(\lambda)\epsilon \Big) + \Phi\Big(-\frac{\Theta}{\tau_*(\lambda)} - \alpha_*(\lambda)+ 2 \alpha_*(\lambda)\epsilon \Big)\right]
		\\
		&\leq s_*(\lambda) + 2 \alpha_{\rm max}\epsilon \,.
	\end{align*}
	We obtain $\P \Big( s_{\epsilon} \leq s_*(\lambda) + 2 \alpha_{\rm max}\epsilon + \epsilon \Big) \geq 1 - e^{-2 N\epsilon^2}$. By Theorem~\ref{th:gordon_aux_v} there exists a constant $\gamma >0$ such that the event
	\begin{equation*}
		\Big\{ \forall v \in B_{\infty}(0,1),  \ V_{\lambda}(v) \geq \max_{\|v\|_{\infty} \leq 1} V_{\lambda}(v) - 3 \gamma \epsilon^3  \implies \frac{1}{N} \|v-\mathsf{v}_{\lambda}\|^2 < \epsilon^3 \Big\}
		\bigcap
		\Big\{  s_{\epsilon} \leq s_*(\lambda) + \alpha_{\rm max} \epsilon + \epsilon \Big\}
	\end{equation*}
	has probability at least $1-\frac{C}{\epsilon^3} e^{-cn\epsilon^6}$. We have on this event, for all $v \in D_{\epsilon}$, $\frac{1}{N} \|v - \mathsf{v}_{\lambda} \|^2 \geq \epsilon^3$.
	Therefore, on the above event we have $\displaystyle \max_{v \in D_{\epsilon}} V_{\lambda}(v) < \max_{\|v\|_{\infty} \leq 1} V_{\lambda}(v) - 3 \gamma \epsilon^3$,
which concludes the proof.
\end{proof}
\\

\begin{proof}[of Theorem~\ref{th:key_1_v}]
Let $\gamma >0$ be a constant that verify the statement of Lemma~\ref{lem:strong_concave_v_1}.
Let $\epsilon \in (0,1]$ and define 
$$ D_{\epsilon} = \left\{ v \in B_{\infty}(0,1) \, \middle| \,  \frac{1}{N} \# \big\{ i \, \big| \, |v_i|\geq 1 - \epsilon \big\} \geq s_*(\lambda) + 2 (1+\alpha_{\rm max})\epsilon \right\} \,.
$$
$D_{\epsilon}$ is a closed set.
\begin{align*}
	\P \Big( \frac{1}{N} \# \big\{ i \, \big| \, |\what{v}_{\lambda,i}|\geq 1 - \epsilon \big\} \geq s_*(\lambda) &+ 2 (1 + \alpha_{\rm max})\epsilon  \Big)
								 \leq \P\Big(\max_{v \in D_{\epsilon}} \mathcal{V}_{\lambda}(v) \geq \max_{\|v\|_{\infty} \leq 1}\mathcal{V}_{\lambda}(v) - \gamma \epsilon \Big)
	\\
	&\leq 2 \P\Big(\max_{v \in D_{\epsilon}} V_{\lambda}(v) \geq \max_{\|v\|_{\infty} \leq 1}V_{\lambda}(v) - 3 \gamma \epsilon \Big) + \frac{C}{\epsilon} e^{-cn\epsilon}
	\leq \frac{C}{\epsilon} e^{-cn\epsilon^2} \,,
\end{align*}
where we used successively Proposition~\ref{prop:apply_gordon_v} and Lemma~\ref{lem:strong_concave_v_1}.
\end{proof}

\subsubsection{Uniform control over \texorpdfstring{$\lambda$}{lambda}: proof of Theorems~\ref{th:unif_law_v},~\ref{th:unif_norm_v} and~\ref{th:unif_1_v}}\label{sec:proof_unif_v}

Theorems~\ref{th:unif_law_v},~\ref{th:unif_norm_v} and~\ref{th:unif_1_v} are deduced from Theorems~\ref{th:law_v},~\ref{th:norm_v} and~\ref{th:key_1_v} by an $\epsilon$-net argument, as we did to deduce Theorems~\ref{th:unif_lambda_law} and~\ref{th:unif_lambda_risk} from Theorems~\ref{th:key_law} and~\ref{th:key_risk}. Since the ideas are the same, we only present here the key argument:

\begin{proposition}
	Assume that $\mathcal{D}$ is $\cF_0(s)$ or $\cF_1(\xi)$ for some $s < s_{\rm max}(\delta)$ and $\xi \geq 0, p > 0$.
	Let $q=0$ if $\mathcal{D} = \cF_0(s)$ and $q=(1/p-1)_+$ if $\mathcal{D} = \cF_p(\xi)$.
	Then there exists constants $K,C,c>0$ that depend only on $\Omega$ such that for all $\theta^{\star} \in \mathcal{D}$
	\begin{equation}\label{eq:lambda_lip_v}
		\P \Big(
			\forall \lambda,\lambda' \in [\lambda_{\rm min},\lambda_{\rm max}],
			\quad
			\frac{1}{N} \|\what{v}_{\lambda} - \what{v}_{\lambda'}\|^2 \leq KN^q |\lambda - \lambda'|
		\Big)
		\geq 1 - C e^{-cn} \,.
	\end{equation}
\end{proposition}
\begin{proof}
	By Proposition~\ref{prop:lipschitz_u_hat}, there exists a constant $K$ such that with probability at least $1 - Ce^{-cn}$ we have
	$$
			\forall \lambda,\lambda' \in [\lambda_{\rm min},\lambda_{\rm max}],
			\quad
			\frac{1}{n} \|\what{u}_{\lambda} - \what{u}_{\lambda'}\|^2 \leq K N^q |\lambda - \lambda'| \,.
	$$
	Notice now that $\what{v}_{\lambda} = -\frac{1}{\lambda} X^{\sT} \what{u}_{\lambda}$ and that with probability at least $1 - 2 e^{-n/4}$, $\sigma_{\rm max}(X)\leq \delta^{-1/2} + 2$ (by Proposition~\ref{prop:max_singular}) which combined with the above inequality, prove the Proposition.
\end{proof}

\section{Some auxiliary results and proofs}

\subsection{Proof of Remark~\ref{rmk:ell0ball}}
\label{app:ell0ball_remark}

Let $k\le N$ and define the  vector  $\theta^\star = (N, 2N, \dots, k N,0,\dots,0)$.
With the definitions given in  Remark~\ref{rmk:ell0ball}, we claim that $W_2(\what{\mu}_{(\htheta_{\lambda},\theta^{\star})},\mu^*_{\lambda})\ge
\sqrt{k/N}$ with probability at least $1-e^{-c k}$, for some constant $c>0$.
Indeed, consider the case $\lambda = 0$, $\tau_* =1$, and let $\P$, $\E$ denote probability and expectation with respect to the coupling that achieves the
Wasserstein distance. This is a coupling for a  triple of random variables
$(I,\Theta,Z)$, with $I\sim {\rm Unif}\big(\{1, \dots, N\}\big)$, $(\Theta,Z)\sim \what{\mu}_{\theta^{\star}}\otimes \cN(0,1)$, with
\begin{align}
W_2(\what{\mu}_{(\htheta_{\lambda},\theta^{\star})},\mu^*_{\lambda})^2 &
= \E\big\{(\theta_I-\Theta)^2\big\}+ \E\big\{(\theta_I+z_I-\Theta-Z)^2\big\}\equiv A+B\, .
\end{align}
We will proceed to bound separately the two terms above. Define $\delta_i \equiv \P(\Theta\neq \theta^{\star}_i|I=i)$, 
and $\delta_{\max} \equiv \max_{i\le k}\delta_i$. Since $\Theta\in \{0,N,\dots,kN\}$ with probability one, we have
\begin{align}
A & \ge \frac{1}{N}\sum_{i=1}^k \E\big\{(\theta^{\star}_I-\Theta)^2\big| I = i\big\}\, \delta_i
\ge  N\sum_{i=1}^k \delta_i\ge N\delta_{\max}\, .
\end{align}
For the second term, we have
\begin{align}
B & \ge\sum_{i=1}^k \E\big\{(\theta^{\star}_I+z_I-\Theta-Z)^2\bfone_{I=i}\bfone_{\Theta=\theta^{\star}_i}\big\}
=  \sum_{i=1}^k \E\big\{(z_i-Z)^2\bfone_{I=i}\bfone_{\Theta=\theta^{\star}_i}\big\}
\\
  &=  \sum_{i=1}^k \E\big\{(z_i-Z)^2\bfone_{\Theta=\theta^{\star}_i}\big\} -\sum_{i=1}^k \E\big\{(z_i-Z)^2\bfone_{I\neq i}\bfone_{\Theta=\theta^{\star}_i}\big\}  \, .
\end{align}
Note that, by the coupling definition, $\P\big(I = i\big|\Theta=\theta^{\star}_i\big) = N \P\big(I = i\, ;\Theta=\theta^{\star}_i\big) = 1-\delta_i$.
Using the fact that $\Theta$ and $Z$ are independent random variables, together with Cauchy-Schwartz inequality,we get
\begin{align}
B & \ge
  \sum_{i=1}^k \E\big\{(z_i-Z)^2\}\,\P(\Theta=\theta^{\star}_i)-\sum_{i=1}^k \E\big\{(z_i-Z)^4 \bfone_{\Theta=\theta^{\star}_i}\big\}^{1/2}\P\big(I\neq i;\Theta=\theta^{\star}_i\big)^{1/2}\\
& \ge \frac{1}{N}\sum_{i=1}^k \E\big\{(z_i-Z)^2\}-\frac{1}{N}\sum_{i=1}^k \E\big\{(z_i-Z)^4 \big\}^{1/2}\P\big(I\neq i\big|\Theta=\theta^{\star}_i\big)^{1/2}\\
& \ge \frac{1}{N}\sum_{i=1}^k \E\big\{(z_i-Z)^2\}-\frac{1}{N}\sum_{i=1}^k \delta_i^{1/2}\E\big\{(z_i-Z)^4 \big\}^{1/2}\\
& \ge \frac{1}{N}\sum_{i=1}^k (1+z_i^2)-\frac{1}{N}\sum_{i=1}^k \delta_i^{1/2}\big(3+6z_i^2+z_i^4\big)^{1/2}\, .
\end{align}
where in the last step we used the fact that $Z\sim\cN(0,1)$. Using $(3+6x^2+x^4)\le 4(1+x^2)^2$,
we thus conclude
\begin{align}
B\ge \frac{1-2\delta^{1/2}_{\max}}{N}\sum_{i=1}^k (1+z_i^2)\, .
\end{align}
By concentration properties of chi-squared random variables, for  any $\eps>0$, there exists $c(\eps)>0$ such that, with probability at least $1-e^{-ck}$ we have
$\frac{1}{k} \sum_{i=1}^k z_i^2 \ge 1-2\eps$. Hence, with the same probability
\begin{align}
W_2(\what{\mu}_{(\htheta_{\lambda},\theta^{\star})},\mu^*_{\lambda})^2 &\ge N\delta_{\max} + \frac{2k}{N}(1-2\delta^{1/2}_{\max})(1-\eps)\\
&\ge \frac{k}{N}\, .
\end{align}
The last inequality follows by lower bounding the first term for  $\delta_{\max}>1/N$, and the second for $\delta_{\max}\le 1/N$, and fixing 
$\eps$ a sufficiently small constant.

\subsection{Concentration properties of \texorpdfstring{$\mathsf{w}_{\lambda}$}{w}}\label{sec:concentration_w}

We prove in this section concentrations of the norms and some scalar product of $\mathsf{w}_{\lambda}$.

\begin{lemma}\label{lem:conc_w_RS}
	There exists constants $c,C >0$ that only depend on $\Omega$ such that for all $t \geq 0$ the event
	$$
	\left\{
		\left|\frac{1}{n} g^{\sT} \mathsf{w}_{\lambda}  - \E\left[ \frac{1}{n} g^{\sT} \mathsf{w}_{\lambda}\right]\right| \leq t
		\,, \quad
		\left| \frac{\|\mathsf{w}_{\lambda}\|^2}{n} - \E \left[ \frac{\|\mathsf{w}_{\lambda}\|^2}{n}\right] \right| \leq t
		\quad \text{and} \quad
		\left| \frac{|\mathsf{w}_{\lambda} + \theta^{\star}|}{n} - \E \left[ \frac{|\mathsf{w}_{\lambda} + \theta^{\star}|}{n} \right] \right| \leq t
	\right\}
	$$
	has probability at least $1-C e^{-c t^2 n} - C e^{-c t n}$.
\end{lemma}
\begin{proof}
	The function
	$$
	g \mapsto \mathsf{w}_{\lambda} = \left( \eta(\theta^{\star}_i + \tau_*(\lambda) g_i, \alpha_*(\lambda) \tau_*(\lambda) ) - \theta^{\star}_i \right)_{1 \leq i\leq N}
	$$
	is $\tau_{\rm max}$-Lipschitz. 
	Consequently:
	\begin{itemize}
		\item $g \mapsto \frac{|\mathsf{w}_{\lambda}+\theta^{\star}|}{n}$ is $\delta^{-1/2}n^{-1/2} \tau_{\rm  max}$-Lipschitz. Therefore $\frac{|\mathsf{w}_{\lambda} + \theta^{\star}|}{n}$ is $\tau_{\rm max}^2 \delta^{-1} n^{-1}$ sub-Gaussian: for all $t \geq 0$,
			$$
			\P \left( \left| \frac{|\mathsf{w}_{\lambda} + \theta^{\star}|}{n} - \E \left[ \frac{|\mathsf{w}_{\lambda} + \theta^{\star}|}{n} \right] \right| > t \right) \leq 2 e^{-n t^2 \delta / \tau_{\rm max}^2} \,.
			$$
		\item $g \mapsto \frac{\|\mathsf{w}_{\lambda}\|}{\sqrt{n}}$ is $n^{-1/2} \tau_{\rm  max}$-Lipschitz. Therefore $\frac{\|\mathsf{w}_{\lambda}\|}{\sqrt{n}}$ is $\tau_{\rm max}^2 n^{-1}$ sub-Gaussian. Its expectation is bounded by $\E\frac{\|\mathsf{w}_{\lambda}\|}{\sqrt{n}} \leq (\E \frac{\|\mathsf{w}_{\lambda}\|^2}{n})^{1/2} = \tau_* \leq \tau_{\rm max}$. By Proposition~\ref{prop:square_sub_gaussian}, we obtain that $\frac{\|\mathsf{w}_{\lambda}\|^2}{n}$ is $(C n^{-1},Cn^{-1})$-sub-Gamma for some constant $C$ and therefore for all $t\geq 0$,
			$$
			\P\left(
			\left| \frac{\|\mathsf{w}_{\lambda}\|^2}{n} - \E \left[ \frac{\|\mathsf{w}_{\lambda}\|^2}{n}\right] \right| > t \right) \leq 2 e^{-cnt^2} + 2 e^{-cnt} \,.
			$$
	\end{itemize}
	Now for $i \in \{1, \dots, N\}$,
	$$
	g_i \mathsf{w}_{\lambda,i} = \tau_* g_i^2 + g_i \big( \mathsf{w}_{\lambda,i} - \tau_* g_i\big) \,.
	$$
	$|\mathsf{w}_{\lambda,i} - \tau_* g_i| \leq \alpha_* \tau_*$ and $g_i$ is $1$-sub-Gaussian and $\E[|g_i|] = \sqrt{\frac{2}{\pi}} \leq 1$. Consequently, Lemma~\ref{lem:conc_prod_bounded} gives that $g_i (\mathsf{w}_{\lambda,i} -\tau_* g_i)$ is $48 \tau_*^2 \alpha_*^2$-sub-Gaussian. This gives that 
	$\frac{1}{n} g^{\sT} (\mathsf{w}_{\lambda} - \tau_* g)$ concentrates exponentially fast around its mean. So does $\frac{1}{n} \|g\|^2$.
\end{proof}

\begin{lemma} \label{lem:opti_w_star}
	$$
	\frac{1}{n}\E \|\mathsf{w}_{\lambda}\|^2 + \sigma^2 = \tau_*^2(\lambda) \qquad \text{and} \qquad \frac{1}{n} \E \big[g^{\sT} \mathsf{w}_{\lambda}\big] = \tau_*(\lambda) - \beta_*(\lambda) = \frac{1}{\delta} s_*(\lambda) \,.
	$$
\end{lemma}
\begin{proof}
	The first equality comes from Lemma~\ref{lem:tau_star}: since $(g_i) \iid \cN(0,1)$, we have $\E \|\mathsf{w}_{\lambda}\|^2 = N \E \big[w^*(\alpha_*,\tau_*)^2 \big]$. The second equality comes from the optimality condition of $\beta_*$, see Lemma~\ref{lem:der_Psi}, and the definition~\eqref{eq:def_s_star} of $s_*(\lambda)$.
\end{proof}
\\

The next proposition simply follows from Lemma~\ref{lem:conc_w_r} and standard concentration arguments, so we omit its proof.
\begin{proposition}\label{prop:concentration_cost_w_star}
	There exists constant $C,c>0$ that only depend on $\Omega$ such that for all $\epsilon \in [0,1]$,
	$$
	\P \Big(\big|  L_{\lambda}(\mathsf{w}_{\lambda})- \psi_{\lambda}(\beta_*(\lambda),\tau_*(\lambda))\big| > \epsilon \Big) \leq C e^{-cn\epsilon^2} \,.
	$$
\end{proposition}

\subsection{Concentration of the empirical distribution}

\begin{proposition}\label{prop:concentration_empirical_distribution}
	Let $\theta^{\star} \in \cF_p(\xi)$, where $p,\xi >0$.
	Let $\mu = \what{\mu}_{\theta^{\star}} \otimes \cN(0,1)$
	and let $\what{\mu}$ be the empirical distribution of the entries of $\big(  \theta_i^{\star} ,  g_i \big)_{1 \leq i \leq N}$, where $g_1,\dots,g_N \iid \cN(0,1)$.
	Then there exists constants $C,c > 0$ that only depends on $\xi^p$, such that for all $\epsilon \in (0,\frac{1}{2}]$,
	$$
	\P\left( W_2 (\what{\mu},\mu)^2 > \epsilon \right) \leq C \epsilon^{-a} \exp\left( -cN \epsilon^2 \epsilon^a  \log(\epsilon)^{-2}\right) \,,
	$$
	where $a = \frac{1}{2} + \frac{1}{p}$.
\end{proposition}

Before proving Proposition~\ref{prop:concentration_empirical_distribution}, we will need two simple lemmas.
	For $r \geq 0$ and $x\in \R$ we use the notation
	$$
	x_{|r} = 
	\begin{cases}
		x & \text{if} \quad -r \leq x \leq r\,, \\
		r & \text{if} \quad x \geq r\,, \\
		-r & \text{if} \quad x \leq -r \,.
	\end{cases}
	$$
	Let $\mu_{| r}$ be the law of $ \big(  \Theta,  Z_{|r} \big)$ where $(\Theta,Z) \sim \what{\mu}_{\theta^{\star}} \otimes \cN(0,1)$.
	\begin{lemma}\label{lem:conc_Z_r}
		$$
		W_2(\mu,\mu_{| r})^2 \leq e^{-r^2/2} \,.
		$$
	\end{lemma}
	\begin{proof}
		We have
		$\displaystyle
			W_2(\mu,\mu_{|r})^2 
			\leq 
			\E\left[
				(Z - Z_{|r})^2
			\right]
			\leq \frac{2}{\sqrt{2\pi}} \int_r^{+\infty} (z-r)^2 e^{-z^2/2} dz
			\leq e^{-r^2/2}$.
	\end{proof}
	\\

	Let $\what{\mu}_{|r}$ be the empirical distribution of the entries of $\big( \theta^{\star}_i, \, g_{i|r} \big)_{1 \leq i \leq N}$.

	\begin{lemma}\label{lem:conc_w_r}
		With probability at least $1 - e^{-\frac{1}{128} N\epsilon^2}$, we have
		$$
		W_2(\what{\mu},\what{\mu}_{|r})^2 \leq \epsilon + e^{-r^2/2} \,.
		$$
	\end{lemma}
	\begin{proof}
		Obviously $ W_2(\what{\mu},\what{\mu}_{|r})^2 
			\leq 
			\frac{1}{N} \sum_{i=1}^N 
			(g_i - g_{i|r})^2$.  
		The function $x \mapsto x - x_{|r}$ is $1$-Lipschitz, so the variables $(g_i - g_{i|r})^2$ are i.i.d.\ $(16,4)$-sub-Gamma. Therefore for all $\epsilon \in [0,1]$,
		$$
		\P \left(
			\frac{1}{N} \sum_{i=1}^N \left(
				g_{i} - g_{i|r}
			\right)^2
			> \E \left(Z - Z_{|r}\right)^2 + \epsilon
		\right) \leq e^{-\frac{1}{128}N\epsilon^2} \,.
		$$
		And we conclude using $\E \big(Z-Z_{|r}\big)^2 \leq e^{-r^2/2}$, which we proved in the lemma above.
	\end{proof}
	\\

	We need now some concentration results for empirical measures, in Wasserstein distance. The next proposition follows from a direct application of Theorem~2 from~\cite{fournier2015rate} to distributions with bounded support. Notice that the results from~\cite{fournier2015rate} are much more general than this.
	\begin{proposition}
		Let $A_1, \dots A_m \iid \nu$ be a collection of i.i.d.\ random variables, bounded by some constant $r>0$. Let 
		$$ 
		\what{\nu}_{m} = \frac{1}{m} \sum_{i=1}^m \delta_{A_i}
		$$
		be the empirical distribution of $A_1 ,\dots,A_m$. Then there exists two absolute constants $c,C >0$ such that for all $t \geq 0$
		$$
		\P\left(W_2(\nu,\what{\nu}_m)^2 \geq r^2 t \right) \leq C \exp(-cm t^2 ) \,.
		$$
	\end{proposition}

	\begin{proof}[of Proposition~\ref{prop:concentration_empirical_distribution}]
	We are now going to couple $\mu_{|r}$ with $\what{\mu}_{|r}$. 
	Let $R>0$.
	Let $k \geq 1$ and let $\delta = 2R/k$.
	Define
	$$
	B_l = \big[-R + (l-1) \delta, -R + l \delta\big) \,,
	$$
	for $l=1, \dots, k$. We define also $B_0 = (-\infty, R) \cup [R, + \infty)$. For $l=0, \dots k$ we write
	$$
	I_l = \{ i \, | \, \theta^{\star}_i \in B_l \} \quad \text{and} \quad N_l = \# I_l \,.
	$$
	Let $t>0$. Let $l \in \{1,\dots,k\}$. The random variables $(g_{i|r})_{i \in I_l}$ are i.i.d.\ and bounded by $r$. By the proposition above, one can couple $i_l \sim {\rm Unif}(I_l)$ with $Z_l \sim \cN(0,1)$ such that we have with probability at least $1 - Ce^{-c t^2 N }$. 
	$$
	\mathsf{E}\left[(Z_{l|r} - g_{i_l|r})^2 \right] \leq tr^2 \sqrt{\frac{N}{N_l}} \,,
	$$
	where $\mathsf{E}$ denotes the expectation with respect to $i_l$ and $Z_l$.
	Let $j_l \sim {\rm Unif}(I_l)$ independently of everything else.

	For $l=0$, we define $(i_0,Z_0) \sim {\rm Unif}(I_0) \otimes \cN(0,1)$, independently of everything else.
	We have with probability at least $1 - Ce^{-ct^2 N}$:
	$$
	\mathsf{E}\left[(Z_{0|r} - g_{i_0|r})^2\right] 
	=
	\mathsf{E}\big[Z_{0|r}^2\big] 
	+
	\mathsf{E}\big[g_{i_0|r}^2\big] 
	\leq 2 + t r^2 \sqrt{\frac{N}{N_0}} \,,
	$$
	where $\mathsf{E}$ denotes the expectation with respect to $Z_0$ and $i_0$. 
	Indeed, $\E\big[g_{i_0|r}^2\big] = \frac{1}{N_0} \sum_{i \in I_0} g_{i|r}^2 \leq 1 + t r^2 \sqrt{\frac{N}{N_0}}$ with probability at least $1 - Ce^{-ct^2 N}$.
	The equality comes from the fact that $Z_0$ and $i_0$ are independent. Finally, we define $j_0 = i_0$.
	\\

	Let us now define the random variable $L$ whose law is given by $\P(L = l) = \frac{N_l}{N}$, independently of everything else.
	Define 
	$$
	\begin{cases}
		Y_1 = \big( \theta^{\star}_{j_L}, \, Z_{L|r}\big)\,, \\
		Y_2 = \big( \theta^{\star}_{i_L}, \, g_{i_L|r}\big) \,.
	\end{cases}
	$$
	$(Y_1,Y_2)$ is a coupling of $(\mu_{| r},\what{\mu}_{|r})$. Let $\mathsf{E}$ denote the expectation with respect to $(i_l,Z_l)_{0 \leq l \leq k}$ and $L$. Then
	\begin{align*}
		\mathsf{E} \left\|Y_1 - Y_2\right\|^2
		&=
		\sum_{l=0}^k \frac{N_l}{N}
		\mathsf{E} \left[ \
			\left(\theta^{\star}_{i_l} - \theta^{\star}_{j_l} \right)^2
			+ \left(Z_{l|r} - g_{i_l|r} \right)^2
		\right]
		\\
		&\leq
		\sum_{l=1}^k \frac{N_l}{N}
		\left( 
			 \sqrt{\frac{N}{N_l}} t r^2 
			+
			\delta^2
		\right)
		+ \frac{N_0}{N} \left( 2 + tr^2 \sqrt{\frac{N}{N_0}}\right)
		\\
		&\leq
		\delta^2 +  \sqrt{k} t  r^2  + 2 \frac{N_0}{N}
		\leq
		\delta^2 +  \sqrt{k} t  r^2  + 2 \frac{\xi^p}{R^p} \,,
	\end{align*}
	with probability at least $1 - C (k+1) e^{-c t^2 N}$, where the last inequality comes from Markov's inequality, since $\theta^{\star} \in \cF_p(\xi)$.
	\\

	Let now $\epsilon \in (0,\frac{1}{2}]$. Let us chose
	$$
	r = \sqrt{- 2 \log(\epsilon)}, \qquad R = \epsilon^{-1/p} \qquad \text{and} \qquad k = \ceil{\epsilon^{-1/2 -1/p}} \leq 2\epsilon^{-1/2 -1/p},
	$$
	so that $\delta = 2R/k \leq 2 \sqrt{\epsilon}$.
	Consequently
	$$
\E \big\| Y_1 - Y_2 \|^2 \leq (4 + 2 \xi^p) \epsilon + 2\sqrt{2} \epsilon^{-1/4 - 1/(2p)} |\log(\epsilon)| t \,.
	$$
	So if we chose $t=|\log(\epsilon)|^{-1} \epsilon^{\frac{5}{4} + \frac{1}{2p}}$ we obtain
	$$
	\P \Big( W_2(\mu_{|r},\what{\mu}_{|r})^2 \leq  (4+2\xi^p + 2\sqrt{2}) \epsilon \Big)
	\geq 1 - C \epsilon^{-1/p - 1/2} \exp(-c N \epsilon^2 \epsilon^{1/2 + 1/p}/ \log(\epsilon)^2).
	$$
	Combining this with Lemmas~\ref{lem:conc_Z_r} and~\ref{lem:conc_w_r} proves the proposition.
\end{proof}

\subsection{Sparsity of the Lasso estimator}\label{sec:proof_sparsity}

The goal of this section is to prove:

\begin{theorem}\label{th:sparsity_uniform}
	Assume here that $\mathcal{D}$ is either $\cF_0(s)$ or $\cF_p(\xi)$ for some $0 \leq s < s_{\rm max}(\delta)$ and $\xi >0, p >0$.
	There exists constants $C,c>0$ that only depend on $\Omega$, such that
	for all $\epsilon \in (0,1)$
	$$
	\sup_{\theta^{\star} \in \mathcal{D}} \ \P \left(
		\sup_{\lambda \in [\lambda_{\rm min},\lambda_{\rm max}]} 
		\Big| \frac{1}{N}\| \what{\theta}_{\lambda} \|_0 - s_*(\lambda)\Big| \geq \epsilon 
	\right) 
	\leq \frac{C}{\epsilon^6} N^q e^{-cN\epsilon^6} \,,
	$$
	where $q = 0$ if\, $\mathcal{D} = \cF_0(s)$ and $q=(1/p - 1)_+$ if\, $\mathcal{D} = \cF_p(\xi)$.
\end{theorem}

Since $\# \{ i \, | \, |\what{v}_{\lambda,i}| =1 \} \geq \|\what{\theta}_{\lambda}\|_0$,
Theorem~\ref{th:unif_1_v} gives that
\begin{equation}\label{eq:unif_upper_bound_sparsity}
	\sup_{\theta^{\star} \in \mathcal{D}} \ \P \left(
		\exists \lambda \in [\lambda_{\rm min},\lambda_{\rm max}], \quad
		\frac{1}{N}\| \what{\theta}_{\lambda} \|_0 \geq s_*(\lambda) + \epsilon 
	\right) 
	\leq \frac{C}{\epsilon^6} N^q e^{-cN\epsilon^6} \,.
\end{equation}


	It remains to prove the converse lower bound in order to get Theorem~\ref{th:sparsity_uniform}.
	We start with the following `local stability' property of the Lasso cost:
\begin{proposition}\label{prop:key_sparsity1}
	There exists constants $C,c,\gamma >0$ that only depend on $\Omega$ such that
	for all $\epsilon \in (0,1]$
	$$
	\sup_{\lambda \in [\lambda_{\rm min},\lambda_{\rm max}]} \
	\sup_{\theta^{\star} \in \mathcal{D}} \
	\P \left(\exists \theta \in \R^N, \quad 
		\frac{1}{N}\|\theta \|_0  < s_*(\lambda) - \epsilon
		\quad \text{and} \quad
	\mathcal{L}_{\lambda}(\theta) \leq \min \mathcal{L}_{\lambda} + \gamma \epsilon^3 \right) \leq \frac{C}{\epsilon^3} e^{-cN\epsilon^6} \,.
	$$
\end{proposition}
Proposition~\ref{prop:key_sparsity1} is a consequence of Proposition~\ref{prop:apply_gordon} and Lemma~\ref{lem:strong_convex_W_L_sparsity} below.

\begin{lemma}\label{lem:strong_convex_W_L_sparsity}
	There exists constants $\gamma,c,C > 0$ that only depend on $\Omega$ such that for all $\epsilon \in (0,1]$ we have
	$$
	\P\left(
		\min_{w \in D_{\epsilon}} L_{\lambda}(w)
		\leq \min_{w\in \R^N}L_{\lambda}(w) + 3\gamma \epsilon^3
	\right)
	\leq \frac{C}{\epsilon^3} e^{-cn\epsilon^6} \,,
	$$
	where $D_{\epsilon} = \left\{ w \in \R^N \, \middle| \frac{1}{N} \| w + \theta^{\star} \|_0 < s_*(\lambda) - \epsilon \right\}$.
\end{lemma}
\begin{proof}
	Define $\mathsf{x}_{\lambda} = \mathsf{w}_{\lambda} + \theta^{\star} = \big( \eta\big(\theta^{\star}_i + \tau_* g_i, \alpha_* \tau_* \big) \big)_{1 \leq i \leq N}$, and for $r > 0$
	$$
	s_{r} = \frac{1}{N} \# \Big\{ i \in \{1,\dots,N\} \, \Big| \, |\mathsf{x}_{\lambda,i}| \geq r \Big\} \,.
	$$
	$s_r$ is a mean of independent Bernoulli random variables, by Hoeffding's inequality we have:
	$$
	\P \Big(s_{r} \geq \P\big(|\Theta + \tau_* Z| \geq \alpha_* \tau_* + r\big) - \frac{\epsilon}{4}\Big) \geq 1 - e^{-N\epsilon^2/8} \,.
	$$
	Compute
	\begin{align*}
		\P\big(|\Theta + \tau_* Z| \geq \alpha_* \tau_* + r\big)
		&=
		\E \left[\Phi\Big(\frac{\Theta}{\tau_*(\lambda)} - \alpha_*(\lambda)- \frac{r}{\tau_*(\lambda)}\Big) + \Phi\Big(-\frac{\Theta}{\tau_*(\lambda)} - \alpha_*(\lambda)- \frac{r}{\tau_*(\lambda)}\Big)\right]
		\\
		&\geq s_*(\lambda) - \frac{r}{\sigma} \,.
	\end{align*}
	Let us chose $r = \sigma \epsilon /4$. We have then $\P \Big(s_{r} \geq s_*(\lambda) - \frac{\epsilon}{2}\Big) \geq 1 - e^{-N\epsilon^2/8}$.
	By Theorem~\ref{th:gordon_aux} there exists a constant $\gamma >0$ such that the event
	\begin{equation}\label{eq:event_good_s}
		\Big\{ \forall w \in \R^N, \ L_{\lambda}(w) \leq \min_{v \in \R^N} L_{\lambda}(v) + 3 \gamma \epsilon^3  \implies \frac{1}{N} \|w-\mathsf{w}_{\lambda}\|^2 < \frac{\sigma^2 \epsilon^3}{32} \Big\}
		\bigcap
		\Big\{  s_{r} \geq s_*(\lambda) - \frac{\epsilon}{2}\Big\}
	\end{equation}
	has probability at least $1-\frac{C}{\epsilon^3} e^{-cn\epsilon^6}$. We have on this event, for all $w \in D_{\epsilon}$
	$$
	\frac{1}{N} \|w - \mathsf{w}_{\lambda} \|^2 
	=
	\frac{1}{N} \|w + \theta^{\star} - \mathsf{x}_{\lambda} \|^2 \geq \frac{\epsilon}{2} r^2 = \frac{\sigma^2 \epsilon^3}{32} \,.
	$$
	Therefore, on the event~\eqref{eq:event_good_v} we have $\displaystyle \min_{w \in D_{\epsilon}} L_{\lambda}(w) > \min_{w \in \R^N} L_{\lambda}(w) + 3 \gamma \epsilon^3$.
	We conclude
	$$
	\P \left(\min_{w \in D_{\epsilon}} L_{\lambda}(w) \leq \min_{w \in \R^N} L_{\lambda}(w) + \gamma \epsilon^3 \right)
	\leq \frac{C}{\epsilon^3} e^{-cn\epsilon^6} \,.
	$$
\end{proof}

Using the same arguments that we use to deduce Theorems~\ref{th:unif_lambda_law} and~\ref{th:unif_lambda_risk}~\eqref{eq:unif_risk} from Theorem~\ref{th:key_law} and Theorem~\ref{th:key_risk} in Section~\ref{sec:proof_unif}, we deduce from Proposition~\ref{prop:key_sparsity1} that for all $\epsilon \in (0,1]$
$$ \sup_{\theta^{\star} \in \mathcal{D} } \P \Big(\exists \lambda \in [\lambda_{\rm min},\lambda_{\rm max}],\quad  \frac{1}{N} \|\what{\theta}_{\lambda} \|_0 < s_*(\lambda) - \epsilon  \Big) \leq \frac{C}{\epsilon^6}N^q e^{-cn\epsilon^6} \,.
	$$
	This proves, together with~\eqref{eq:unif_upper_bound_sparsity}, Theorem~\ref{th:sparsity_uniform}.




	\subsection{Proof of Theorem~\ref{th:law_debiased}}\label{sec:proof_law_debiased}

	Recall that the distributions $\mu_{\lambda}^*$ and $\nu_{\lambda}^*$ are respectively defined by Definition~\ref{def:mu_star} and~\eqref{eq:def_nu_star}.
	Let $\epsilon \in (0,1]$. From now, we will work on the event
	\begin{align*}
		\mathcal{E} = &\Big\{ \forall \lambda \in [\lambda_{\rm min},\lambda_{\rm max}], \quad
			W_2(\what{\mu}_{(\what{\theta}_{\lambda},\theta^{\star})},\mu^*_{\lambda})^2
			+
			W_2(\what{\mu}_{(\what{v}_{\lambda},\theta^{\star})},\nu^*_{\lambda})^2 
			\leq \epsilon^6
		\Big\}
		\\
		& \qquad \bigcap \
		\Big\{ \forall \lambda \in [\lambda_{\rm min},\lambda_{\rm max}], \quad
		\Big| \frac{1}{N}\|\what{\theta}_{\lambda}\|_0 - s_*(\lambda) \Big|
		+
		\Big| \frac{1}{N} \# \big\{ i \, \big| \, |\what{v}_{\lambda,i}| =1 \big\} - s_*(\lambda) \Big|
			\leq \epsilon^2
		\Big\} \,,
	\end{align*}
	which has probability at least $1-C \epsilon^{-12} e^{-cN\epsilon^{17}}$ from what we have just seen, and Theorems~\ref{th:unif_lambda_law},~\ref{th:unif_law_v},~\ref{th:sparsity_uniform} and~\ref{th:unif_1_v}.
	From now, $\mathsf{E}$ and $\mathsf{P}$ will denote the probability with respect to the empirical distributions of the entries of the vectors we study, and the variables that we couple with them.
	Let $\lambda \in [\lambda_{\rm min},\lambda_{\rm max}]$. 
	On the event $\mathcal{E}$ one can couple $(\Theta^{x},Z^{x})\sim \hat{\mu}_{\theta^{\star}} \otimes \cN(0,1)$ and  $(\Theta^{v},Z^{v})\sim \hat{\mu}_{\theta^{\star}} \otimes \cN(0,1)$ with
	$(\Theta,\what{\Theta}_{\lambda},\what{V}_{\lambda},\what{\Theta}_{\lambda}^d)$ which is sampled from the empirical distribution of the entries of $(\theta^{\star},\what{\theta}_{\lambda}, \what{v}_{\lambda}, \what{\theta}_{\lambda}^d)$,
	such that
	 \begin{align*}
		 &\mathsf{E} \Big[
		\big(\what{\Theta}_{\lambda} - \eta(\Theta^{x} + \tau_* Z^x, \alpha_* \tau_*)\big)^2
		+
\big(\Theta - \Theta^x\big)^2
	\Big] \leq \epsilon^6 \,,
		\\
		 &
		 \mathsf{E} \Big[ 
			 \big(\what{V}_{\lambda} + \frac{1}{\alpha_* \tau_*} \big( \eta(\Theta^{v} + \tau_* Z^v, \alpha_* \tau_*) - \Theta^v - \tau_* Z_0^v \big) \big)^2
		+
		 \big(\Theta - \Theta^v\big)^2 \Big]
		 \leq \epsilon^6 \,.
	 \end{align*}
	 Let 
	 $$
	 E_1 = \Big\{
		\big|\what{\Theta}_{\lambda} - \eta(\Theta^{x} + \tau_* Z^x, \alpha_* \tau_*)\big| \leq \epsilon^2, \quad
		\big| \what{V}_{\lambda}  + \frac{1}{\alpha_* \tau_*}\big( \eta(\Theta^{v} + \tau_* Z^v, \alpha_* \tau_*) - \Theta^v - \tau_* Z^v \big)\big| \leq \epsilon^2
	 \Big\} \,.
	 $$
	 By Chebychev's inequality, $\mathsf{P}(E_1) \geq 1 - C \epsilon^2$, for some constant $C>0$.
	 Let us also define the event
	 $$
	 E_2 = \Big\{ \Theta^x + \tau_* Z^x \neq \alpha_* \tau_* \quad \text{and} \quad \Theta^v + \tau_* Z^v \neq \alpha_* \tau_* \Big\} \,.
	 $$
$\Theta^x + \tau_* Z^x$ and $\Theta^v + \tau_* Z^v$ admit a density with respect to Lebesgue's measure. Therefore $\mathsf{P}(E_2)=1$.
	 \begin{lemma} The event
		 $$
		 E_3 = \Big\{ |\what{\Theta}_{\lambda}| \not\in (0,\epsilon^2]  \quad \text{and} \quad | \what{V}_{\lambda}| \not\in [1- \epsilon^2, 1) \Big\}
		 $$
		 has probability at least $1 - C \epsilon^2$.
	 \end{lemma}
\begin{proof}
	We denote here by $O(\epsilon^2)$ quantities that are bounded by $C\epsilon^2$, from some constant $C$.
	Since $\Theta^x + \tau_* Z^x$ admits a density with respect to Lebesgue's measure we have
	 $$
	 \mathsf{P} \Big(
		 \big|\eta(\Theta^{x} + \tau_* Z^x, \alpha_* \tau_*)\big| \not\in (0,2\epsilon^2]
	 \Big) = 1 - O(\epsilon^2)\,.
	 $$
	Consequently, since the events $E_1$ has probability at least $1-O(\epsilon^2)$, we have
	\begin{align*}
		&\mathsf{P}(|\what{\Theta}_{\lambda}| \in [0,\epsilon^2]) 
				  \\
				  &=
		\mathsf{P}\Big(|\what{\Theta}_{\lambda}| \in [0,\epsilon^2] \quad \text{and} \quad |\eta(\Theta^x+\tau_*Z^x,\alpha_*\tau_*)| \not\in (0,2\epsilon^2]
			\quad \text{and} \quad \big|\what{\Theta}_{\lambda} - \eta(\Theta^{x} + \tau_* Z^x, \alpha_* \tau_*)\big| \leq \epsilon^2
		\Big) + O(\epsilon^2)
		\\
		&=
		\mathsf{P}\Big( \eta\big(\Theta^x+\tau_*Z^x,\alpha_*\tau_*\big) = 0 \Big) + O(\epsilon^2) = s_*(\lambda) + O(\epsilon^2) \,.
	\end{align*}
	Since $\mathsf{P}(\what{\Theta}_{\lambda} = 0) = s_*(\lambda) + O(\epsilon^2)$ because we are working on $\mathcal{E}$, we conclude that $\mathsf{P}\big(|\what{\Theta}_{\lambda}| \in (0,\epsilon^2]\big) = O(\epsilon^2)$.
	One can prove the same way that $\mathsf{P}\big(|\what{V}_{\lambda}| \in [1- \epsilon^2,1)\big) = O(\epsilon^2)$, which gives the desired result.
\end{proof}
\\

\begin{lemma}
	The event
	$$
	E_4 = \Big\{\what{\Theta}_{\lambda} \neq 0 \iff \what{V}_{\lambda} = {\rm sign}(\what{\Theta}_{\lambda}) \Big\}
	$$
	has probability at least $1-C\epsilon^2$, for some constant $C>0$.
\end{lemma}
\begin{proof}
	Since $\what{v}_{\lambda} \in \partial | \what{\theta}_{\lambda}|$, $\what{\theta}_{\lambda,i} > 0$ implies that $\what{v}_{\lambda,i} = {\rm sign}(\what{\theta}_{\lambda,i})$. Thus $\mathsf{P}\big(\what{\Theta}_{\lambda} \neq 0 \implies \what{V}_{\lambda} = {\rm sign}(\what{\theta}_{\lambda,i}) \big) = 1$. We have thus
	\begin{equation}\label{eq:inclusion}
		\big\{\what{\Theta}_{\lambda} \neq 0\big\} \subset \big\{|\what{V}_{\lambda}| = 1 \big\} \,.
	\end{equation}
		On the event $\mathcal{E}$ we have $
		\Big| \frac{1}{N}\|\what{\theta}_{\lambda}\|_0 - s_*(\lambda) \Big|
		+
		\Big| \frac{1}{N} \# \big\{ i \, \big| \, |\what{v}_{\lambda,i}| =1 \big\} - s_*(\lambda) \Big| \leq \epsilon^2$ which gives
		$$
		\mathsf{P}\big(\what{\Theta}_{\lambda} \neq 0\big) = s_*(\lambda) + O(\epsilon^2)
		\qquad \text{and} \qquad
		\mathsf{P}\big(|\what{V}_{\lambda}| = 1\big) = s_*(\lambda) + O(\epsilon^2).
		$$
		We deduce then from~\eqref{eq:inclusion} that $\mathsf{P}\big(|\what{V}_{\lambda}| = 1 \ \text{and} \ \what{\Theta}_{\lambda} = 0 \big) = O(\epsilon^2)$ and finally
$
	\mathsf{P}\big( \what{V}_{\lambda} = {\rm sign}(\what{\Theta}_{\lambda}) \implies \what{\Theta}_{\lambda} \neq 0 \big) \geq 1 - C \epsilon^2 \,.
		$
\end{proof}
\\

\begin{lemma}\label{lem:equiv_debiased}
	Let $E = E_1 \cap E_2 \cap E_3 \cap E_4$. The event $E$ has probability at least $1 - C\epsilon^2$ and on $E$ we have
	$$
	\Theta^v + \tau_*Z^v \geq  \alpha_* \tau_* \iff
	\Theta^x + \tau_*Z^x \geq  \alpha_* \tau_*  \,,
	$$
	and
	$$
	\Theta^v + \tau_*Z^v \leq  -\alpha_* \tau_* \iff
	\Theta^x + \tau_*Z^x \leq  -\alpha_* \tau_*  \,. 
	$$
\end{lemma}
\begin{proof}
	Since $E_1,E_2,E_3$ and $E_4$ have all a probability greater than $1 - O(\epsilon^2)$, the event $E = E_1 \cap E_2 \cap E_3 \cap E_4$ has probability at least $1- O(\epsilon^2)$. On $E$ we have
\begin{align*}
	\Theta^v + \tau_*Z^v \geq  \alpha_* \tau_* &\iff
		\Theta^v + \tau_* Z_0^v - \eta(\Theta^{v} + \tau_* Z^v, \alpha_* \tau_*) = \alpha_* \tau_* 
		\\
		&\iff \what{V}_{\lambda} \geq 1 - \epsilon^2 \qquad \text{(because we are on the event} \ E_1 \text{)}
		\\
		&\iff \what{V}_{\lambda} = 1
\qquad \text{(because we are on the event} \ E_3 \text{)}
		\\
		&\iff \what{\Theta}_{\lambda} > 0 \qquad \text{(because we are on the event} \ E_4 \text{)}\\
		&\iff \what{\Theta}_{\lambda} > \epsilon^2
\qquad \text{(because we are on the event} \ E_3 \text{)}
			\\
			&\iff \eta(\Theta^x + \tau_* Z^x,\alpha_* \tau_*) > 0
\qquad \text{(because we are on the event} \ E_1 \text{)}
\\
&\iff \Theta^x + \tau_* Z^x \geq \alpha_* \tau_* \qquad \text{(because we are on the event} \ E_2 \text{)}\,.
\end{align*}
The second equivalence is proved exactly the same way.
\end{proof}
\\

Let us define 
$$
X^d = 
		 \eta(\Theta^{x} + \tau_* Z^x, \alpha_* \tau_*) +  \Theta^v + \tau_* Z_0^v - \eta(\Theta^{v} + \tau_* Z^v, \alpha_* \tau_*) \,.
$$
We have 
\begin{align*}
\mathsf{E} (\what{\Theta}_{\lambda}^d - X^d)^2 
&= 
2\mathsf{E} \Big[\big(\what{\Theta}_{\lambda} - \eta(\Theta^{x} + \tau_* Z^x, \alpha_* \tau_*)\big)^2\Big]
+
2\mathsf{E} \Big[\big(\frac{\lambda}{1-\frac{1}{n}\|\what{\theta}_{\lambda}\|_0}\what{V}_{\lambda} - \Theta^v - \tau_* Z^v + \eta(\Theta^{v} + \tau_* Z^v, \alpha_* \tau_*)\big)^2\Big]
\\
&\leq C \epsilon^4 \,,
\end{align*}
for some constant $C>0$, because on the event $\mathcal{E}$, $\frac{1}{N} \big| \| \what{\theta}_{\lambda} \|_0 - s_*(\lambda) \big| \leq \epsilon^2$, so 
$$
\frac{\lambda}{1 - \frac{1}{n} \| \what{\theta}_{\lambda} \|_0}
=
\frac{\lambda}{1 - \frac{1}{\delta}s_*(\lambda)}
+ O(\epsilon^2)
= \alpha_* \tau_* + O(\epsilon^2) \,.
$$
By Lemma~\ref{lem:equiv_debiased} above, we have on the event $E$,
$$
X^d = 
\begin{cases}
	\Theta^x + \tau_* Z^x & \text{if} \quad \Theta^x +\tau_* Z^x \geq \alpha_* \tau_* \quad \text{or} \quad \Theta^x +\tau_* Z^x \leq -\alpha_* \tau_* \,, \\
	\Theta^v + \tau_* Z^v & \text{otherwise}.
\end{cases}
$$
Let us denote $T^x = (\Theta^x + \tau_* Z^x, \Theta^x)$ and $T^v = (\Theta^v + \tau_* Z^v, \Theta^v)$. 

Since $\Theta^x + \tau_*^x Z$ and $\Theta^v + \tau_* Z^v$ have the same law and $\mathsf{P}\big(\Theta^x + \tau_* Z^x \geq \alpha_* \tau_* \big| E\big) = \mathsf{P}\big( \Theta^v + \tau_* Z^v \geq \alpha_* \tau_*\big|E\big)$ (by Lemma~\ref{lem:equiv_debiased}), we have $\mathsf{P}\big( \Theta^x + \tau_* Z^x \geq \alpha_* \tau_* \big| E^c\big) = \mathsf{P}\big(\Theta^v + \tau_* Z^v \geq \alpha_* \tau_*\big|E^c\big)$. Similarly we have $\mathsf{P}\big( \Theta^x + \tau_* Z^x \leq - \alpha_* \tau_* \big| E^c\big) = \mathsf{P}\big(\Theta^v + \tau_* Z^v \leq - \alpha_* \tau_*\big|E^c\big)$.
\\

One can therefore define two random variables $\widetilde{T}^x = (\widetilde{\Theta}^x + \tau_* \widetilde{Z}^x, \, \widetilde{\Theta}^x)$ and $\widetilde{T}^v=(\widetilde{\Theta}^v + \tau_* \widetilde{Z}^v, \, \widetilde{\Theta}^v)$ such that
\begin{itemize}
	\item conditionally on $E^c$, $\widetilde{T}^x$ (respectively $\widetilde{T}^v$) and $T^x$ (respectively $T^v$) have the same law.
	\item On the event $E^c$, $\widetilde{\Theta}^x + \tau_* \widetilde{Z}^x \geq \alpha_* \tau_* \iff  \widetilde{\Theta}^v + \tau_* \widetilde{Z}^v \geq \alpha_* \tau_*$ and $\widetilde{\Theta}^x + \tau_* \widetilde{Z}^x \leq -\alpha_* \tau_* \iff  \widetilde{\Theta}^v + \tau_* \widetilde{Z}^v \leq - \alpha_* \tau_*$.
\end{itemize}

We define then
$$
\big( \widetilde{X}^d, \widetilde{\Theta} \big)
=
\begin{cases}
	T^x & \text{on the event} \ E \ \text{provided that} \quad |\Theta^x + \tau_*Z^x| \geq \alpha_* \tau_* \,, \\
	T^v & \text{on the event} \ E \ \text{provided that} \quad |\Theta^x + \tau_* Z^x| < \alpha_* \tau_* \,, \\
	\widetilde{T}^x & \text{on the event} \ E^c \ \text{provided that} \quad |\widetilde{\Theta}^x + \tau_* \widetilde{Z}^x| \geq \alpha_* \tau_* \,, \\
	\widetilde{T}^v & \text{on the event} \ E^c \ \text{provided that} \quad |\widetilde{\Theta}^x + \tau_* \widetilde{Z}^x| < \alpha_* \tau_* \,.
\end{cases}
$$
$(\widetilde{X}^d,\widetilde{\Theta}) \sim \mu_{\lambda}^d$ which is the law of $(\Theta + \tau_* Z,\Theta)$ where $(\Theta,Z) \sim \what{\mu}_{\theta^{\star}} \otimes \cN(0,1)$. Indeed, for every continuous bounded function $f$ we have
\begin{align*}
	\mathsf{E}[f(\widetilde{X}^d,\widetilde{\Theta})]	
	&=
	\mathsf{E} \Big[\bbf{1}_{E} \bbf{1}_{|\Theta^x + \tau_* Z^x|\geq \alpha_* \tau_*}f(T^x)\Big]
	+\mathsf{E} \Big[\bbf{1}_{E} \bbf{1}_{|\Theta^v + \tau_* Z^v|< \alpha_* \tau_*}f(T^v)\Big]
	\\
	&\qquad+\mathsf{E} \Big[\bbf{1}_{E^c} \bbf{1}_{|\widetilde{\Theta}^x + \tau_* \widetilde{Z}^x| \geq \alpha_* \tau_*}f(\widetilde{T}^x)\Big]
	+\mathsf{E} \Big[\bbf{1}_{E^c} \bbf{1}_{| \widetilde{\Theta}^v + \tau_* \widetilde{Z}^v|< \alpha_* \tau_*}f(\widetilde{T}^v)\Big]
	\\
	&=
	\mathsf{E} \Big[\bbf{1}_{E} \bbf{1}_{|\Theta^x + \tau_* Z^x|\geq \alpha_* \tau_*}f(T^x)\Big]
+\mathsf{E} \Big[\bbf{1}_{E} \bbf{1}_{|\Theta^v + \tau_* Z^v|< \alpha_* \tau_*}f(T^v)\Big]
\\
&\qquad+\mathsf{E} \Big[\bbf{1}_{E^c} \bbf{1}_{|\Theta^x + \tau_* Z^x| \geq \alpha_* \tau_*}f(T^x)\Big]
+\mathsf{E} \Big[\bbf{1}_{E^c} \bbf{1}_{| \Theta^v + \tau_* Z^v|< \alpha_* \tau_*}f(T^v)\Big]
\\
&=
\mathsf{E} \Big[\bbf{1}_{|\Theta^x + \tau_* Z^x|\geq \alpha_* \tau_*}f(T^x)\Big]
+\mathsf{E} \Big[\bbf{1}_{|\Theta^v + \tau_* Z^v| < \alpha_* \tau_*}f(T^v)\Big]
\\
&=
\mathsf{E} \Big[\bbf{1}_{|\Theta^x + \tau_* Z^x|\geq \alpha_* \tau_*}f(T^x)\Big]
+\mathsf{E} \Big[\bbf{1}_{|\Theta^x + \tau_* Z^x| < \alpha_* \tau_*}f(T^x)\Big]
=
\mathsf{E} \Big[f(T^x)\Big]\,.
\end{align*}

Let us now compute
\begin{align*}
\mathsf{E} \Big[\big(X^d - \widetilde{X}^d\big)^2\Big]
=
\mathsf{E} \Big[\bbf{1}_{E^c}\big(X^d - \widetilde{X}^d\big)^2\Big]
\leq
C \sqrt{\mathsf{P}(E^c)} \leq C \epsilon \,,
\end{align*}
and 
\begin{align*}
\mathsf{E}\Big[\big(\widetilde{\Theta} - \Theta\big)^2\Big]
&\leq
\mathsf{E}\Big[\bbf{1}_E \big(\Theta^x - \Theta\big)^2\Big]
+\mathsf{E}\Big[\bbf{1}_E \big(\Theta^v - \Theta\big)^2\Big]
+\mathsf{E}\Big[\bbf{1}_{E^c} \big(\widetilde{\Theta}^x - \Theta\big)^2\Big]
+\mathsf{E}\Big[\bbf{1}_{E^c} \big(\widetilde{\Theta}^v - \Theta\big)^2\Big]
\\
&\leq 2 \epsilon^6 + 2 C \sqrt{\mathsf{P}(E^c)} \leq C \epsilon \,.
\end{align*}

Therefore $\mathsf{E} \big\| (\what{\Theta}_{\lambda}^d,\Theta) - (\widetilde{X}^d,\widetilde{\Theta}) \big\|^2 \leq C\epsilon$ and consequently $W_2(\what{\mu}_{(\what{\theta}_{\lambda}^d,\theta^{\star})},\mu_{\lambda}^d)^2 \leq C \epsilon$, on the event $\mathcal{E}$ which has probability at least $1-C \epsilon^{-12} e^{-cN\epsilon^{17}}$.

\subsection{Proof of Corollary~\ref{cor:estim_risk}}\label{sec:proof_estim_risk}
	Let $\epsilon \in (0,1]$.
	Let us work on the intersection the events of Theorem~\ref{th:sparsity_uniform},Corollary~\ref{cor:estim_tau} and~\ref{th:norm_v}, which as probability at least $1 - \frac{C}{\epsilon^6} e^{-cN\epsilon^6}$. Let $\lambda \in [\lambda_{\rm min},\lambda_{\rm max}]$.
	$$
	\frac{1}{N}\big\| X^{\sT}(y-X\what{\theta}_{\lambda}) \|^2 = \frac{\lambda^2}{N} \|\what{v}_{\lambda}\|^2 = \lambda^2 \kappa_*(\lambda) + O(\epsilon) \,.
	$$
	We have also $1 - \frac{1}{n} \|\what{\theta}_{\lambda}\|_0 = 1 - \frac{1}{\delta} s_*(\lambda) + O(\epsilon)= \beta_*(\lambda)/\tau_*(\lambda) + O(\epsilon)$. Therefore
	\begin{align*}
		\frac{\big\|X^{\sT}(y-X\what{\theta}_{\lambda})\big\|^2}{N(1-\frac{1}{n}\|\what{\theta}_{\lambda}\|_0)^2} = \Big(\frac{\lambda \tau_*(\lambda)}{\beta_*(\lambda)}\Big)^2 \kappa_*(\lambda) + O(\epsilon)
		= \tau_*(\lambda)^2 \big(1 + \delta - 2 s_*(\lambda)\big) - \delta \sigma^2 + O(\epsilon) \,.
	\end{align*}
	Now we have $\what{\tau}(\lambda) = \tau_*(\lambda) + O(\epsilon)$ and $\frac{1}{N} \|\what{\theta}_{\lambda}\|_0 = s_*(\lambda) + O(\epsilon)$. Consequently
	$$
	\what{\tau}(\lambda)^2 \Big(\frac{2}{N}\|\what{\theta}_{\lambda}\|_0 -1 \Big) = \tau_*(\lambda)^2 (2 s_*(\lambda) - 1) + O(\epsilon) \,.
	$$
	Putting all together we obtain $\what{R}(\lambda) = \delta \tau_*(\lambda)^2 - \delta \sigma^2 + O(\epsilon)= R_*(\lambda) + O(\epsilon)$, and we conclude using Theorem~\ref{th:unif_lambda_risk}.

	\subsection{Proof of Proposition~\ref{prop:riskCV}}\label{app:proof_riskCV}

	Let $n' \in \{1, \dots, n\}$. We consider a random $n' \times N$ matrix $X'$ and a random vector $z' = (z_1', \dots, z_{n'}')$ such that $X'_{i,j} \iid \cN(0, 1 / n)$ and $z'_i \iid \cN(0,1)$ are independent and independent of everything else.

\begin{lemma}
	There exists constants $\gamma, c,C > 0$ that only depend on $\Omega$ such that for all $\theta^{\star}$ in $\mathcal{D}$ and all $\lambda \in [\lambda_{\rm min},\lambda_{\rm max}]$
	such that for all $\epsilon \in (0,1]$,
	$$
	\P \Big(
		\exists w \in \R^N, \quad \Big| \frac{1}{n'} \| X' w - \sigma z' \|^2 - \frac{1}{n} \|w\|^2 - \sigma^2 \Big| \geq \sqrt{\epsilon} \frac{n}{n'}
		\quad \text{and} \quad 
		L_{\lambda}(w) \leq \min_{v \in \R^N} L_{\lambda}(v) + \gamma \epsilon
	\Big) \leq \frac{C}{\epsilon} e^{-cn\epsilon^2}.
	$$
\end{lemma}
\begin{proof}
	The vector $\mathsf{w}_{\lambda}$ is independent from $X'$, $z'$. Hence
	\begin{align}\label{eq:chi1}
		\| X' \mathsf{w}_{\lambda} - \sigma z' \|^2
		= \Big(\frac{1}{n} \| \mathsf{w}_{\lambda} \|^2 + \sigma^2 \Big) \chi
	\end{align}
	where $\chi$ is independent from $\mathsf{w}_{\lambda}$ and follows a $\chi$-squared distribution with $n'$ degrees of freedom. We have therefore for all $t \geq 0$
	\begin{equation}\label{eq:chi2}
		\P \Big(\big| \chi - n' \big| \geq t n' \Big) \leq C e^{- c n' t } + C e^{- c n' t^2},
	\end{equation}
	for some constants $c,C > 0$.
	We know by Lemma~\ref{lem:conc_w_RS} and Lemma~\ref{lem:opti_w_star} that $\frac{1}{n} \| \mathsf{w}_{\lambda}\|^2$ concentrates exponentially fast around $\tau_*(\lambda)^2 - \sigma^2$, which is (Theorem~\ref{th:control_beta_tau}) bounded by some constant. There exists therefore constants $C, c > 0$ such that 
	\begin{equation}\label{eq:bound_norm_hp}
	\P \Big(\frac{1}{n} \| \mathsf{w}_{\lambda} \|^2 + \sigma^2 > C\Big) \leq C e^{-c n}.
	\end{equation}
	From~\eqref{eq:chi1}-\eqref{eq:chi2} and~\eqref{eq:bound_norm_hp} above, we deduce that for all $t \geq 0$
	\begin{equation}\label{eq:concentration_sf}
	\P \Big(
		\big| \frac{1}{n'}\| X' \mathsf{w}_{\lambda} - \sigma z' \|^2 - \frac{1}{n} \|\mathsf{w}_{\lambda}\|^2 - \sigma^2 \big| > t
	\Big) \leq C e^{-c n' t} + C e^{- cn't^2} + C e^{-cn},
	\end{equation}
	for some constants $c,C>0$. By Proposition~\ref{prop:max_singular}, we know that $\P\big(\sigma_{\rm max}(X') > \delta^{-1/2} + \sqrt{n'/n} + 1 \big) \leq e^{-n /2}$. 
	Let $\epsilon \in (0,1]$.
	Let $w \in \R^N$ such that $\|w - \mathsf{w}_{\lambda}\|^2 \leq \epsilon N$.
	\begin{align*}
\frac{1}{\sqrt{n'}}
\| X' \mathsf{w}_{\lambda} - \sigma z' \|
+
\frac{1}{\sqrt{n'}} \| X' w - \sigma z' \|
\leq 
\frac{2\sigma}{\sqrt{n'}} \| z' \|
+
\frac{\sigma_{\rm max}(X')}{\sqrt{n'}}\big( \| \mathsf{w}_{\lambda} \| + \| w \| \big)
\leq C\sqrt{\frac{n}{n'}}
	\end{align*}
	for some constant $C>0$, with probability at least $1 - C e^{-cn}$. Consequently
	\begin{align*}
\Big| 
\frac{1}{n'}\| X' \mathsf{w}_{\lambda} - \sigma z' \|^2
-
\frac{1}{n'}\| X' w - \sigma z' \|^2
\Big|
&\leq \frac{C\sqrt{n}}{n'}
\Big| 
\| X' \mathsf{w}_{\lambda} - \sigma z' \|
-
\| X' w - \sigma z' \|
\Big|
\leq \frac{C\sqrt{n}}{n'} \|X' (\mathsf{w}_{\lambda} - w) \|
\\
&\leq C \sigma_{\rm max}(X')\sqrt{\delta^{-1} \epsilon} \frac{n}{n'}
\leq C \sqrt{\epsilon} \frac{n}{n'}
	\end{align*}
	with probability at least $1 - C e^{-cn}$ for some constant $C > 0$. Similarly, we have with probability at least $1 - C e^{-cn}$, 
$$
\Big| \frac{1}{n} \| \mathsf{w}_{\lambda} \|^2  - \frac{1}{n} \| w \|^2 \Big|
\leq C \sqrt{\epsilon}.
$$
We conclude that with probability at least $1 - C e^{-cn}$ we have for all $w \in \R^N$ such that $\|w - \mathsf{w}_{\lambda} \|^2 \leq N \epsilon$,
$$
		\Big| 
		\frac{1}{n'}\| X' \mathsf{w}_{\lambda} - \sigma z' \|^2 - \frac{1}{n} \|\mathsf{w}_{\lambda}\|^2
		-
		\frac{1}{n'}\| X' w + \sigma z' \|^2 + \frac{1}{n} \|w\|^2
		\Big| 
		\leq C \sqrt{\epsilon} \big(1 + \frac{n}{n'} \big) 
		\leq 2 C \sqrt{\epsilon} \frac{n}{n'}.
$$
Combining this with~\eqref{eq:concentration_sf}, we get that 
for all $\epsilon \in (0,1]$,
$$
\P \Big( \exists w \in \R^N, \quad \| w - \mathsf{w}_{\lambda}\|^2 \leq N \epsilon \quad \text{and} \quad  \Big| \frac{1}{n'}\| X' w + \sigma z' \|^2 - \frac{1}{n} \|w\|^2 - \sigma^2 \Big| > 2 C \sqrt{\epsilon}\frac{n}{n'}
\Big) \leq  C e^{-cn \epsilon}.
$$
We conclude using Theorem~\ref{th:gordon_aux} that
	$$
	\P \Big(
		\exists w \in \R^N, \quad \Big| \frac{1}{n'} \| X' w - \sigma z' \|^2 - \frac{1}{n} \|w\|^2 - \sigma^2 \Big| \geq \sqrt{\epsilon} \frac{n}{n'}
		\quad \text{and} \quad 
		L_{\lambda}(w) \leq \min_{v \in \R^N} L_{\lambda}(v) + \gamma \epsilon
	\Big) \leq  \frac{C}{\epsilon} e^{-cn\epsilon^2}
	$$
	for some constants $c,C, \gamma > 0$.
\end{proof}
\\

Using Proposition~\ref{prop:apply_gordon}, we deduce
\begin{lemma}\label{lem:chi}
	There exists constants $\gamma, c,C > 0$ that only depend on $\Omega$ such that for all $\theta^{\star}$ in $\mathcal{D}$ and all $\lambda \in [\lambda_{\rm min},\lambda_{\rm max}]$
	such that for all $\epsilon \in (0,1]$,
	$$
	\P \Big(
		\exists \theta \in \R^N, \ 
		\Big| \frac{1}{n'} \| X' \theta^{\star} + \sigma z' - X' \theta \|^2 - \frac{1}{n} \|\theta - \theta^{\star} \|^2 - \sigma^2 \Big| \geq \sqrt{\epsilon} \frac{n}{n'}
		\quad \text{and} \quad 
		\mathcal{L}_{\lambda}(\theta) \leq \min \mathcal{L}_{\lambda} + \gamma \epsilon
	\Big) \leq \frac{C}{\epsilon} e^{-cn\epsilon^2}.
	$$
\end{lemma}

\vspace{2mm}

	We have
	\begin{align*}
\what{\theta}_{\lambda}^{i} 
&= \underset{\theta \in \R^N}{\text{arg\,min}} \left\{
	\frac{1}{2 n_k} \Big\|y^{(\sminus i)} - X^{(\sminus i)}\theta\Big\|^2 + \frac{\lambda}{n} |\theta|
\right\} 
= \underset{\theta \in \R^N}{\text{arg\,min}} \left\{
\frac{1}{2 n_k} \Big\|X^{(\sminus i)} \theta^{\star} + \sigma z^{(\sminus i)} - X^{(\sminus i)}\theta\Big\|^2 + \frac{\lambda}{n} |\theta|
\right\} 
\\
&= \underset{\theta \in \R^N}{\text{arg\,min}} \left\{
\frac{1}{2 n} \Big\|\sqrt{\frac{k}{k-1}}X^{(\sminus i)} \theta^{\star} + \sqrt{\frac{k}{k-1}} \sigma z^{(\sminus i)} - \sqrt{\frac{k}{k-1}}X^{(\sminus i)}\theta\Big\|^2 + \frac{\lambda}{n} |\theta|
\right\}.
	\end{align*}
	$\what{\theta}^i_{\lambda}$ is thus the minimizer of the Lasso cost~\eqref{eq:def_lasso_est} for $\delta^{(k)} = \frac{k-1}{k} \delta$ and $\sigma^{(k)} = \sqrt{k/(k-1)} \sigma$. Let $\tau_*^{(k)}(\lambda)$ be the $\tau_*$ defined by Theorem~\ref{th:scalar_max_min}, but with $\delta^{(k)}$ instead of $\delta$ and $\sigma^{(k)}$ instead of $\sigma$. Define the corresponding `risk':
$$
R_*^{(k)}(\lambda) = \delta^{(k)}\big(\tau^{(k)}_*(\lambda)^2 - (\sigma^{(k)})^2 \big) \,.
$$
	It is not difficult to verify that the bounds on $\tau_*,\beta_*$ of Section~\ref{sec:control_tau_beta} are uniform with respect to $\delta$ and $\sigma$. More precisely
	$$
	\sup_{\delta \in [\delta_{\rm min},\delta_{\rm max}]}
	\sup_{\sigma \in [\sigma_{\rm min},\sigma_{\rm max}]}
	\sup_{\lambda \in [\lambda_{\rm min},\lambda_{\rm max}]}
	\sup_{\theta^{\star} \in \mathcal{D}}
	\big\{ \tau_*(\lambda,\delta,\sigma) + \beta_*(\lambda,\delta,\sigma) \big\} < +\infty  \,,
	$$
	where $\delta_{\rm max},\delta_{\rm min},\sigma_{\rm max},\sigma_{\rm min} > 0$ such that $s_{\rm max}(\delta_{\rm min}) > s$ if we are in the case $\mathcal{D} = \cF_0(s)$.
	This gives that under the assumptions of Proposition~\ref{prop:riskCV}, $\tau_*^{(k)}$ and $R_*^{(k)}$ are bounded for all $k \geq 2$ (that verify $s_{\rm max}(\delta(k-1)/k) > s$ in the case $\mathcal{D} = \cF_0(s)$) by some constant that depends only on $\Omega$.

	\begin{lemma}\label{lem:chi2}
		There exists constants $C,c >0$ that only depend on $\Omega$ such that for all $\theta^{\star} \in \mathcal{D}$, for all $i\in \{1, \dots,k \}$ and for all $\epsilon \in (0,1]$,
	$$
	\P \Big(\sup_{\lambda \in [\lambda_{\rm min},\lambda_{\rm max}]}
		\Big| \frac{k}{n} \| y^{(i)} - X^{(i)} \what{\theta}^{i}_{\lambda} \|^2 - \frac{1}{n} \|\what{\theta}_{\lambda}^{i} - \theta^{\star} \|^2 - \sigma^2 \Big| \geq k \epsilon 
	\Big) \leq C N^q \epsilon^{-4} e^{-cn \epsilon^4}.
	$$
\end{lemma}
\begin{proof}
	Let $i \in \{1, \dots, k\}$.
	Let us define for $\theta \in \R^N$,
	$$
	\mathcal{L}_{\lambda}^{(\sminus i)}(\theta) = \frac{1}{2 n_k} \big\| y^{(\sminus i)} - X^{(\sminus i)} \theta \big\|^2 + \frac{\lambda}{n} | \theta |.
	$$
	Let $\epsilon \in (0,1]$.
	Let $\eta = \frac{\gamma \epsilon}{K N^q}$ and $M = \ceil{(\lambda_{\rm max} - \lambda_{\rm min}) / \eta}$. Define for $j \in \{0, \dots, M \}$, define $\lambda_j = \min \big( \lambda_{\rm min} + j \eta, \lambda_{\rm max} \big)$. 
	We apply Lemma~\ref{lem:chi} with $n' = n / k$, $X' = X^{(i)}$ and $z' = z^{(i)}$ to obtain that the event
\begin{align*}
	E_1 = 
	\Big\{\forall j \in \{1, \dots, M\}, \quad 
		\forall \theta \in \R^N, \quad
		&\mathcal{L}_{\lambda_j}^{(\sminus i)}(\theta) \leq \min \mathcal{L}_{\lambda_j}^{(\sminus i)} + \gamma \epsilon
	\\
													   &\implies
		\Big| \frac{k}{n} \| y^{(i)} - X^{(i)} \theta \|^2 - \frac{1}{n} \|\theta - \theta^{\star} \|^2 - \sigma^2 \Big| < \sqrt{\epsilon} k
	\Big\}
\end{align*}
has probability at least $1 - M \frac{C}{\epsilon} e^{-c \epsilon^2 n}$. By Lemma~\ref{lem:lambda_lip} the event
\begin{equation}
	E_2 = 
	\Big\{ \forall \lambda,\lambda' \in [\lambda_{\rm min},\lambda_{\rm max}], \ 
		\mathcal{L}_{\lambda'}^{(\sminus i)}(\what{\theta}^i_{\lambda}) \leq \min_{x \in \R^N} \mathcal{L}^{(\sminus i)}_{\lambda'}(x) + K N^q | \lambda - \lambda'|
	\Big\}
\end{equation}
has probability at least $1 - C e^{-cn}$.
On the event $E_2$, we have for all $j \in \{1, \dots, k\}$ and all $\lambda \in [\lambda_{j-1},\lambda_{j}]$
\begin{align*}
	\mathcal{L}_{\lambda_j}^{(\sminus i)}(\what{\theta}_{\lambda}^i)
	\leq \min_{x \in \R^N} \mathcal{L}^{(\sminus i)}_{\lambda_j}(x) + K N^q \eta
	\leq \min_{x \in \R^N} \mathcal{L}^{(\sminus i)}_{\lambda_j}(x) + \gamma \epsilon.
\end{align*}
We obtain that on $E_1 \cap E_2$, which has probability at least $1 - C N^q \epsilon^{-2} e^{-c n \epsilon^2}$
$$
\forall \lambda \in [\lambda_{\rm min}, \lambda_{\rm max} ], \quad
\Big| \frac{k}{n} \| y^{(i)} - X^{(i)} \what{\theta}^i_{\lambda} \|^2 - \frac{1}{n} \|\what{\theta}^i_{\lambda} - \theta^{\star} \|^2 - \sigma^2 \Big| < \sqrt{\epsilon}k.
$$
\end{proof}

\begin{proposition}\label{prop:rcv}
	There exists constants $c,C > 0$ that only depend on $\Omega$, such that for all $\theta^{\star} \in \mathcal{D}$ and for all $i \in \{1, \dots, k \}$
	$$
	\P \Big(
		\sup_{\lambda \in [\lambda_{\rm min},\lambda_{\rm max}]} 
		\Big|
		\frac{1}{N} \| \what{\theta}^i_{\lambda} - \theta^{\star} \|^2 - R_*(\lambda)
		\Big| \leq \frac{C}{\sqrt{k}}
	\Big) \leq C N^q k^4 e^{-c N / k^4}.
	$$
\end{proposition}
\begin{proof}
	Let us fix $i \in \{1,\dots,k\}$. 
	By Proposition~\ref{prop:risk_lipschitz}, $\lambda \mapsto R_*(\lambda)$ is $K_1$-Lipschitz on $[\lambda_{\rm min},\lambda_{\rm max}]$, for some constant $K_1 > 0$.
	By Propositions~\ref{prop:lipschitz_lambda_l1} and~\ref{prop:lipschitz_lambda_s} there exists a constant $K_2>0$ such that the event
	\begin{equation}
		E_1 = 
		\Big\{ \forall \lambda \in [\lambda_{\rm min},\lambda_{\rm max}], \quad
		\frac{1}{N} \big| |\what{\theta}^i_{\lambda}| - | \theta^{\star} | \big| \leq K_2 N^q 
		\quad \text{and} \quad
		\frac{1}{N} \big| |\what{\theta}_{\lambda}| - | \theta^{\star} | \big| \leq K_2 N^q 
	\Big\}
	\end{equation}
	has probability at least $1- Ce^{-cn}$. Let us define $\eta = \min \Big(\frac{\delta}{2 N^q k K_2}, \frac{1}{K_1 \sqrt{k}}\Big)$ and $M = \ceil{(\lambda_{\rm max} - \lambda_{\rm min}) / \eta}$. For all $j \in \{0, \dots, M\}$, we write $\lambda_j = \min\big(\lambda_{\rm min} + j \eta, \, \lambda_{\rm max}\big)$.

By Theorem~\ref{th:unif_lambda_risk} the event
\begin{equation}\label{eq:appli_risk1}
	E_2 =
	\Big\{
		\sup_{\lambda \in [\lambda_{\rm min},\lambda_{\rm max}]} 
		\Big| \frac{1}{N}\| \what{\theta}^i_{\lambda} - \theta^{\star} \|^2 - R^{(k)}_*(\lambda) \Big| \leq 1
	\Big\}
		\ \bigcap \
		\Big\{
		\sup_{\lambda \in [\lambda_{\rm min},\lambda_{\rm max}]} 
		\Big| \frac{1}{n} \big\| y - X \what{\theta}_{\lambda}\big\|^2 - \beta_*(\lambda) \Big| \leq 1
	\Big\}
\end{equation}
has probability at least $1-C N^q e^{-cN}$.
By Lemma~\ref{lem:chi2}, applied with $\epsilon = k^{-1}$,
	$$
	E_3 = 
	\Big\{\sup_{\lambda \in [\lambda_{\rm min},\lambda_{\rm max}]}
		\Big| \frac{k}{n} \| y^{(i)} - X^{(i)} \what{\theta}^{i}_{\lambda} \|^2 - \frac{1}{n} \|\what{\theta}_{\lambda}^{i} - \theta^{\star} \|^2 - \sigma^2 \Big| \leq 1
	\Big\}
	$$
	has probability at least $1 -C k^4 N^q e^{-cn/k^4}$.
	On the event $E_2 \cap E_3$, we have, for all $\lambda \in [\lambda_{\rm min},\lambda_{\rm max}]$,
	\begin{align*}
		\mathcal{L}_{\lambda}(\what{\theta}_{\lambda}^i) &= 
		\frac{1}{2n} \big\| y^{(i)} - X^{(i)} \what{\theta}_{\lambda}^i \big\|^2
		+
		\frac{1}{2n} \big\| y^{(\text{-}i)} - X^{(\text{-}i)} \what{\theta}_{\lambda}^i \big\|^2
		+ \frac{\lambda}{n} | \what{\theta}_{\lambda}^i |
		\\
		&\leq
		\frac{1}{k} \big( \sigma^2  + \frac{1}{n} \| \theta^{\star} - \what{\theta}_{\lambda}^i \|^2 +1 \big)
		+
		\frac{1}{2n_k} \big\| y^{(\text{-}i)} - X^{(\text{-}i)} \what{\theta}_{\lambda}^i \big\|^2
		+ \frac{\lambda}{n} | \what{\theta}_{\lambda}^i |
		\\
		&\leq
		\frac{1}{k} \big( 1 + \sigma^2  + \delta^{-1}(R_*^{(k)}(\lambda) + 1) \big)
		+
		\frac{1}{2n_k} \big\| y^{(\text{-}i)} - X^{(\text{-}i)} \what{\theta}_{\lambda} \big\|^2
		+ \frac{\lambda}{n} | \what{\theta}_{\lambda} | 
		\\
		&\leq
		\frac{1}{k} \big( 1+  \sigma^2  + \delta^{-1}(R_*^{(k)}(\lambda) + 1) \big)
		+
		\mathcal{L}_{\lambda}(\what{\theta}_{\lambda})
		+\Big(\frac{1}{2 n_k} - \frac{1}{2n}\Big) \big\| y - X \what{\theta}_{\lambda} \big\|^2
		\\
		&\leq 
		\mathcal{L}_{\lambda}(\what{\theta}_{\lambda}) + \frac{C}{k}
	\end{align*}
	for some constant $C>0$. Let $j \in \{1 ,\dots , M\}$.
	We have $\mathcal{L}_{\lambda}(\what{\theta}^i_{\lambda}) = \mathcal{L}_{\lambda_j}(\what{\theta}^i_{\lambda}) - \frac{\lambda_j - \lambda}{n} | \what{\theta}^i_{\lambda}|$ and
	$$
	\mathcal{L}_{\lambda}(\what{\theta}_{\lambda})
		\leq
		\mathcal{L}_{\lambda}(\what{\theta}_{\lambda_j}) 
	=
		\mathcal{L}_{\lambda_j}(\what{\theta}_{\lambda_j})
		+ \frac{\lambda - \lambda_j}{n} |\what{\theta}_{\lambda_j}|
		= \min_{\theta \in \R^N} \mathcal{L}_{\lambda_j}(\theta)
		+ \frac{\lambda - \lambda_j}{n} |\what{\theta}_{\lambda_j}|.
	$$
	So we get that on the event $E_1 \cap E_2 \cap E_3$, for all $j \in \{1, \dots, M\}$ and all $\lambda \in [\lambda_{j-1},\lambda_j]$,
	\begin{align*}
	\mathcal{L}_{\lambda_j}(\what{\theta}^i_{\lambda}) 
	&\leq
		\min_{\theta \in \R^N} \mathcal{L}_{\lambda_j}(\theta)
		+ \frac{\lambda_j - \lambda}{n}\big( |\what{\theta}^i_{\lambda_j}| - |\what{\theta}_{\lambda_j}| \big) + \frac{C}{k}
	\leq
		\min_{\theta \in \R^N} \mathcal{L}_{\lambda_j}(\theta)
		+ \frac{\lambda_j - \lambda}{\delta}2 K_2 N^q + \frac{C}{k}
		\\
	&\leq
		\min_{\theta \in \R^N} \mathcal{L}_{\lambda_j}(\theta)
		+ \frac{\eta}{\delta}2 K_2 N^q + \frac{C}{k}
	\leq
		\min_{\theta \in \R^N} \mathcal{L}_{\lambda_j}(\theta)
		+ \frac{C_0}{k},
	\end{align*}
	for some constant $C_0>0$,
	because on $E_1$ we have $\forall \lambda \in [\lambda_{\rm min},\lambda_{\rm max}], \frac{1}{N} \big| |\what{\theta}^i_{\lambda}| - |\what{\theta}_{\lambda}| \big| \leq 2 K_2 N^q$.
	By Theorem~\ref{th:key_risk}, there exists constants $C,c,\gamma > 0$ such that for all $\epsilon \in (0,1]$ the event
	$$
	E_4 = 
	\Big\{ \forall j \in \{1, \dots, M\},  \quad
		\forall \theta \in \R^N, \quad
		\mathcal{L}_{\lambda_j}(\theta) \leq \min \mathcal{L}_{\lambda_j} + \gamma \epsilon
		\implies 
	\Big| \frac{1}{N} \| \theta - \theta^{\star} \|^2 - R_*(\lambda_j) \Big| \leq \sqrt{\epsilon} \Big\}
	$$
	has probability at least $1 - C M \epsilon^{-1} e^{-cN\epsilon^2}$. 
	Consider the constant $\kappa = \frac{C_0}{\gamma}$.
	If $k \geq \kappa$, then $\epsilon = \frac{C_0}{\gamma k} \leq 1$ and the event $E_4$ has probability at least $1 - C M k e^{-cN / k^2}$. So we obtain that on the event $E_1 \cap E_2 \cap E_3 \cap E_4$, which has probability $1 -C N^q k^4 e^{-cn / k^4} $,
	$$
	\forall j \in \{1, \dots, M\}, \quad \forall \lambda \in [\lambda_{j-1},\lambda_j], \quad
	\Big| \frac{1}{N} \| \what{\theta}^i_{\lambda} - \theta^{\star} \|^2 - R_*(\lambda_j) \Big|
	\leq \frac{C}{\sqrt{k}},
	$$
	for some constant $C>0$.
	If now $k < \kappa$. Then on the event $E_2$ we have
\begin{align*}
	\forall j \in \{1, \dots, M\}, \quad \forall \lambda \in [\lambda_{j-1},\lambda_j], \quad
	\Big| \frac{1}{N} \| \what{\theta}^i_{\lambda} - \theta^{\star} \|^2 - R_*(\lambda_j) \Big|
	&\leq 
	\sup_{\lambda \in [\lambda_{\rm min},\lambda_{\rm max}]} \!\! R_*^{(k)}(\lambda)
	+
	\sup_{\lambda \in [\lambda_{\rm min},\lambda_{\rm max}]} \!\! R_*(\lambda) + 1
	\\
	&\leq C \leq \frac{C \sqrt{\kappa}}{\sqrt{k}},
\end{align*}
	where $C$ is a constant.
	We conclude that (in both cases) there exists a constant $C> 0$ such that
	$$
	\forall j \in \{1, \dots, M\}, \quad \forall \lambda \in [\lambda_{j-1},\lambda_j], \quad
	\Big| \frac{1}{N} \| \what{\theta}^i_{\lambda} - \theta^{\star} \|^2 - R_*(\lambda_j) \Big|
	\leq \frac{C}{\sqrt{k}},
	$$
	holds with probability at least $1 -  C N^q k^4 e^{-cN/k^4}$
	Proposition~\ref{prop:rcv} follows from the fact that for all $\lambda \in [\lambda_{j-1},\lambda_j]$, $|R_*(\lambda) - R_*(\lambda_j)| \leq K_1 |\lambda - \lambda_j| \leq \frac{1}{\sqrt{k}}$.
\end{proof}
\\

\begin{proof}[of Proposition~\ref{prop:riskCV}]
	We apply Lemma~\ref{lem:chi2} with $\epsilon = k^{-3/2}$ to obtain that with probability at least $1 - C k^6 N^q e^{-c n /k^6}$ we have
	$$
	\forall \lambda \in [\lambda_{\rm min},\lambda_{\rm max}], \quad
	\forall i \in \{1, \dots, k \}, \quad
	\Big| \frac{k}{n} \| y^{(i)} - X^{(i)} \what{\theta}^{i}_{\lambda} \|^2 - \frac{1}{n} \|\what{\theta}_{\lambda}^{i} - \theta^{\star} \|^2 - \sigma^2 \Big| \leq \frac{1}{\sqrt{k}}.
	$$
	By summing these inequalities for $i = 1 \dots k$ and using the triangular inequality, we get
	$$
	\Big| \frac{k}{n} \sum_{i=1}^k \| y^{(i)} - X^{(i)} \what{\theta}^{i}_{\lambda} \|^2 - \frac{1}{n} \sum_{i=1}^k \|\what{\theta}_{\lambda}^{i} - \theta^{\star} \|^2 - k \sigma^2 \Big| \leq \frac{k}{\sqrt{k}}.
	$$
	and then
	\begin{equation}\label{eq:bb}
	\Big| \frac{1}{N} \sum_{i=1}^k \| y^{(i)} - X^{(i)} \what{\theta}^{i}_{\lambda} \|^2 - \frac{1}{k} \sum_{i=1}^k \frac{1}{N}\|\what{\theta}_{\lambda}^{i} - \theta^{\star} \|^2 - \delta \sigma^2 \Big| \leq \frac{\delta}{\sqrt{k}}.
	\end{equation}
	By Proposition~\ref{prop:rcv}, we have with probability at least $1 - C N^q k^4 e^{-cN/k^3}$,
	$$
	\forall \lambda \in [\lambda_{\rm min},\lambda_{\rm max}], \quad
	\forall i \in \{1, \dots, k \}, \quad
		\Big|
		\frac{1}{N} \| \what{\theta}^i_{\lambda} - \theta^{\star} \|^2 - R_*(\lambda)
		\Big| \leq \frac{C}{\sqrt{k}}.
	$$
	This implies (again by summing and using the triangular inequality) that
	$$
	\Big|
	\frac{1}{k} \sum_{i=1}^k \frac{1}{N}\|\what{\theta}_{\lambda}^{i} - \theta^{\star} \|^2
	- R_*(\lambda) \Big| \leq \frac{C}{\sqrt{k}},
	$$
	which, combined with~\eqref{eq:bb} proves Proposition~\ref{prop:riskCV}.
\end{proof}

\subsection{The scalar lasso}

In this section we study
\begin{equation}\label{eq:scalar_lasso}
	\ell_{\alpha}(y)=
	\min_{x \in \R} \left\{
		\frac{1}{2}(y-x)^2 + \alpha |x|
	\right\} \,.
\end{equation}
\begin{lemma} 
	The minimum~\eqref{eq:scalar_lasso} is achieved at an unique point $x^* = \eta(y,\alpha)$ and 
	$$
	\ell_{\alpha}(y) = 
	\left\{
		\begin{array}{ccc}
			\frac{1}{2} y^2 & \text{if} & -\alpha \leq y \leq \alpha \\
			\alpha y - \frac{1}{2} \alpha^2 & \text{if} & y \geq \alpha \\
			-\alpha y - \frac{1}{2} \alpha^2 & \text{if} & y \leq -\alpha
		\end{array}
	\right.
	$$
\end{lemma}
Suppose now that
$$
y = x + Z \,,
$$
for some $x \in \R$ and $Z \sim \cN(0,1)$.
\begin{lemma}\label{lem:prop_Delta}
	Define
	\begin{align*}
		\Delta_{\alpha}(x) = \E \big[ \ell_{\alpha}(x+Z) - \alpha |x| \big] \,.
	\end{align*}
	The function $\Delta_{\alpha}$ is continuous, even, decreasing on $\R_{\geq 0}$, $\alpha$-Lipschitz. Moreover
	$$
	\begin{cases}
		\Delta_{\alpha}(0) = \frac{1}{2} + \alpha \phi(\alpha) - (1+\alpha^2) \Phi(-\alpha)\\
		\lim\limits_{x \to \pm \infty} \Delta_{\alpha}(x) = -\frac{\alpha^2}{2}
	\end{cases}
	$$
	and $\Delta_{\alpha}'(0^+)=-\Delta_{\alpha}'(0^-) = -\alpha$.
\end{lemma}

\begin{proof}
	Since $Z$ and $-Z$ have the same law, one verify easily that $\Delta_{\alpha}$ is an even function. We have for all $x >0$
	$$
	\Delta_{\alpha}'(x) = \E \left[\ell_{\alpha}'(x+Z) - \alpha \right]
	= \E \left[\bbf{1}(x+Z \in [-\alpha,\alpha]) (x+Z -\alpha) \right] \leq 0 \,.
	$$
	$\ell_{\alpha}$ is convex, therefore $x \mapsto \E[\ell'_{\alpha}(x+Z)]$ is non-decreasing. $\E[\ell_{\alpha}'(Z)] = 0$ because $\ell_{\alpha}'$ is an odd function. Consequently, for all $x > 0$
	$$
	-\alpha \leq \E[\ell_{\alpha}'(x+Z) -\alpha]  = \Delta_{\alpha}'(x) \,.
	$$
	This gives (recall that $\Delta_{\alpha}$ is even and continuous over $\R$) that $\Delta_{\alpha}$ is $\alpha$-Lipschitz. From what we have seen above, we have also $\Delta_{\alpha}'(0^+) = - \Delta_{\alpha}'(0^-) = -\alpha$. Compute now, using the fact that $\ell_{\alpha}$ is even:
	$$
	\Delta_{\alpha}(0) = \E [\ell_{\alpha}(Z)]
	=
	\int_{0}^{\alpha} z^2 \phi(z)dz
	+
	\int_{\alpha}^{+\infty} (2 \alpha z - \alpha^2)\phi(z)dz \,.
	$$
	By integration by parts
	$$
	\int_0^{\alpha} z^2 \phi(z) dz
	=
	\Big[- z \phi(z) \Big]_0^{\alpha}
	+\int_0^{\alpha} \phi(z) dz
	= -\alpha \phi(\alpha) + \frac{1}{2} - \Phi(-\alpha)
	$$
	$$
	\int_{\alpha}^{+\infty} (2 \alpha z - \alpha^2)\phi(z)dz
	=
	-\alpha^2 \Phi(-\alpha) + 2 \alpha \phi(\alpha) \,.
	$$
	Therefore $\Delta_{\alpha}(0) = \frac{1}{2} + \alpha \phi(\alpha) - (1+\alpha^2) \Phi(-\alpha)$. We have almost-surely
	$$
	\ell_{\alpha}(x+Z) - \alpha |x| \xrightarrow[x \to \pm \infty]{} -\frac{\alpha^2}{2} \,.
	$$
	Thus, by dominated convergence $\lim\limits_{x \to \pm \infty} \Delta_{\alpha}(x) = -\frac{\alpha^2}{2}$.
\end{proof}

\subsection{A convexity lemma}

\begin{lemma}\label{lem:convex_sqrt}
	The function
	$$
	f : x \in \R^N \mapsto \sqrt{\frac{\|x\|^2}{n} + \sigma^2}
	$$
	is $\frac{\sigma^2}{n(R^2 + \sigma^2)^{3/2}}$-strongly convex on $B(0,\sqrt{n}R)$.
\end{lemma}
\begin{proof}
	Let $x,y \in B(0,\sqrt{n}R)$ and define for $t\in[0,1]$, $g(t) = f(z_t)$, where $z_t = (t x + (1-t) y)$.  Compute
	$$
	g'(t)
	= \frac{\frac{1}{n}(x-y)^{\sT} z_t}{\sqrt{\frac{\|z_t\|^2}{n} + \sigma^2}} \,,
	$$
	and 
	\begin{align*}
		g''(t) 
		&= \frac{\frac{1}{n}\|x-y\|^2}{\sqrt{\frac{\|z_t\|^2}{n} + \sigma^2}}
		- \frac{\big(\frac{1}{n}(x-y)^{\sT}z_t\big)^2}{\big(\frac{\|z_t\|^2}{n} + \sigma^2\big)^{3/2}}
		\\
		&= 
		\frac{1}{\big(\frac{\|z_t\|^2}{n} + \sigma^2\big)^{3/2}}
		\left(
			\frac{1}{n}\|x-y\|^2 \Big(\frac{\|z_t\|^2}{n} + \sigma^2\Big) -
			\left(\frac{1}{n}(x-y)^{\sT}z_t\right)^2
		\right)
		\\
		&\geq
		\frac{\sigma^2}{\big(\frac{\|z_t\|^2}{n} + \sigma^2\big)^{3/2}}
		\left(
			\frac{1}{n}\|x-y\|^2
		\right)
		\geq
		\frac{1}{n}\|x-y\|^2
		\frac{\sigma^2}{(R^2 + \sigma^2)^{3/2}} \,.
	\end{align*}
	Consequently
	\begin{align*}
		tf(x) + (1-t) f(y)
		= tg(1) + (1-t) g(0)
		&\geq g(t) + \frac{1}{2}t(1-t)
		\frac{1}{n}\|x-y\|^2
		\frac{\sigma^2}{(R^2 + \sigma^2)^{3/2}}
		\\
		&= f(tx + (1-t)y) + \frac{1}{2}t(1-t)
		\frac{1}{n}\|x-y\|^2
		\frac{\sigma^2}{(R^2 + \sigma^2)^{3/2}} \,.
	\end{align*}
\end{proof}

\section{Toolbox}

\subsection{Notations recap}\label{sec:recap}

Recall that $X$ is a $n \times N$ random matrix with entries $X_{i,j} \iid \cN(0,1/n)$. The random vectors $z \in \R^n$, $g \in \R^N$ and $h \in \R^n$ are standard Gaussian random vectors. The following table displays the main cost (or objective) functions used in this paper and their corresponding optimizers.
\begin{table}[H]
\vspace{-3mm}
$$
\arraycolsep=10.8pt\def\arraystretch{1.5}
\begin{array}{|l|c|}
	  \hline
	  \text{Definition} & \text{Optimizer}\\
  \hline
  \mathcal{L}_{\lambda}(\theta) = \frac{1}{2n}\left\| X \theta - y \right\|^2 + \frac{\lambda}{n} |\theta| & \what{\theta}_{\lambda} \\
  \mathcal{C}_{\lambda}(w) = \frac{1}{2n}\| X w - \sigma z \|^2 + \frac{\lambda}{n} (|w + \theta^{\star}| - |\theta^{\star}|) & \what{w}_{\lambda} \\
\mathcal{U}_{\lambda}(u) = \min\limits_{w \in \R^N} \big\{ u^{\sT} X w - \sigma u^{\sT} z - \frac{1}{2} \|u\|^2 + \lambda (|\theta^{\star} + w| - |\theta^{\star}|) \big\} & \what{u}_{\lambda} \\
\mathcal{V}_{\lambda}(v)= \min\limits_{w \in B}  \left\{
\frac{1}{2n} \| X w - \sigma  z\|^2 + \frac{\lambda}{n} v^{\sT}(\theta^{\star}+w) - \frac{\lambda}{n}|\theta^{\star}| \right\} & \what{v}_{\lambda} \\
	L_{\lambda}(w) = 
	\frac{1}{2}\left(\sqrt{\frac{\|w\|^2}{n} + \sigma^2}\ \frac{\|h\|}{\sqrt{n}}
	-	\frac{1}{n} g^{\sT} w  + \frac{g'\sigma}{\sqrt{n}} \right)_{\!\!+}^2
	+ \frac{\lambda}{n} |w+\theta^{\star}|
	- \frac{\lambda}{n} |\theta^{\star}| & w^*_{\lambda} \\
U_{\lambda}(u) = \min\limits_{w \in \R^N} \big\{
		\frac{-1}{n^{3/2}} \|u\| g^{\sT} w 
		+\frac{1}{n^{3/2}} \|w\| h^{\sT} u
		- \frac{\sigma}{n} u^{\sT} z - \frac{1}{2n} \|u\|^2 + 
	\frac{\lambda}{n} \big(  |w+\theta^{\star}| - |\theta^{\star}| \big) \big\} & u^*_{\lambda} \\
	V_{\lambda}(v) 
	= 
	\min\limits_{w \in B}
	\left\{
		\frac{1}{2} \left(
			\sqrt{ \frac{1}{n}\|w\|^2  + \sigma^2 } \frac{\|h\|}{\sqrt{n}}
		-	\frac{1}{n} g^{\sT} w  + \frac{g'\sigma}{\sqrt{n}} \right)_{\!+}^2
		+ \frac{\lambda}{n} v^{\sT}(w+\theta^{\star})
		- \frac{\lambda}{n} |\theta^{\star}|
	\right\} & v^*_{\lambda} \\
  \hline
\end{array}
$$
\vspace{-3mm}
\caption{Main cost/objective functions}
\label{tab:cost}
\end{table}

In the definition of $V_{\lambda}$ above, $B = \big\{w\in\R^N \, | \, |w| \leq 2 |\theta^{\star}| + 5 \sigma^2 \lambda_{\rm min}^{-1} n +K \big\}$, where $K >0$ is the constant given by Lemma~\ref{lem:kappa}.
The functions $L_{\lambda}$, $U_{\lambda}$ and $V_{\lambda}$ are the ``corresponding cost/objective functions'' to $\mathcal{C}_{\lambda}$, $\mathcal{U}_{\lambda}$ and $\mathcal{V}_{\lambda}$. A main part of the analysis is to show that $w^*_{\lambda}$, $u^*_{\lambda}$ and $v^*_{\lambda}$ are approximately equal to $\mathsf{w}_{\lambda}$, $\mathsf{u}_{\lambda}$ and $\mathsf{v}_{\lambda}$ given by:

\begin{table}[H]
$$
\arraycolsep=10.8pt\def\arraystretch{1.5}
\begin{array}{|l|}
	\hline
	\mathsf{w}_{\lambda} = \eta\big(\theta^{\star} + \tau_*(\lambda) g, \alpha_*(\lambda) \tau_*(\lambda)\big) - \theta^{\star} \\
	\mathsf{u}_{\lambda} = 
	\frac{\beta_*(\lambda)}{\tau_*(\lambda)} \Big(\sqrt{\tau_*(\lambda)^2 - \sigma^2} \frac{h}{\sqrt{n}} - \frac{\sigma}{\sqrt{n}} z\Big) \\
	\mathsf{v}_{\lambda} = -\alpha_*(\lambda)^{-1} \tau_*(\lambda)^{-1} \Big( \eta\big(\theta^{\star} + \tau_*(\lambda) g,\alpha_*(\lambda) \tau_*(\lambda)\big) - \theta^{\star} - \tau_*(\lambda) g\Big)\\
	\hline
\end{array}
$$
\vspace{-3mm}
\caption{``Asymptotic optimizers''}
\label{tab:opt}
\end{table}

\subsection{Convex analysis lemmas}

\begin{proposition}[Corollary 37.3.2 from~\cite{rockafellar2015convex}]\label{prop:rockafellar}
Let $C$ and $D$ be non-empty closed convex sets in $\R^m$ and $\R^n$, respectively, and let $f$ be a continuous finite concave-convex function on $C \times D$. If either $C$ or $D$ is bounded, one has
$$
\inf_{v \in D} \sup_{u \in C} f(u,v) = \sup_{u \in C} \inf_{v \in D} f(u,v) \,.
$$
\end{proposition}

\begin{definition}\label{def:smooth}
A convex function $f$ over $\R^n$ is said to be
\begin{itemize}
	\item $\gamma$-strongly convex if $x \mapsto f(x) - \frac{\gamma}{2} \|x\|^2$ is convex.
	\item $L$-strongly smooth is $f$ is differentiable everywhere and for all $x,y \in \R^n$ we have
		$$
		f(y) \leq f(x) + (y-x)^{\sT} \nabla f(x) + \frac{L}{2} \|x-y\|^2 \,.
		$$
\end{itemize}
\end{definition}
\begin{remark}
If $f$ is convex, differentiable over $\R^n$, and $\nabla f$ is $L$-Lipschitz, then $f$ is $L$-strongly smooth. Indeed, if we take $x,y \in \R^n$ and if we define $h(t)=f((1-t)x  +t y)$ we have
\begin{align*}
	f(y)-f(x) &= h(1)-h(0) = \int_0^1 h'(t) dt = \int_0^1 (y-x)^{\sT} \nabla f ((1-t) x + ty) dt
	\\
	&= (y-x)^{\sT} f(x) +  \int_0^1 (y-x)^{\sT} \big( \nabla f ((1-t) x + ty) - \nabla f(x) \big) dt
	\\
	&\leq 
(y-x)^{\sT} f(x) +  \int_0^1 t L \|x-y\|^2 dt
\leq
(y-x)^{\sT} f(x) +  \frac{L}{2} \|x-y\|^2 \,.
\end{align*}
\end{remark}

\begin{proposition}\label{prop:duality}
Let $f$ be a closed convex function over $\R^n$. Then $f$ is $\gamma$-strongly convex if and only if $f^*$ is $\frac{1}{\gamma}$-strongly smooth.
\end{proposition}

This result can be found in the book~\cite{zalinescu2002convex}, see Corollary 3.5.11 on page 217 and the Remark 3.5.3 below. A more accessible presentation of this result can be found in~\cite{kakade2012regularization}.

\subsection{Gaussian min-max Theorem}\label{sec:proof_gordon}

In this section, we reproduce the proof of the tight Gaussian min-max comparison theorem from~\cite{thrampoulidis2015regularized} for completeness, but also because we need a slightly more general version of this result.

We recall the classical Gordon's min-max Theorem from~\cite{gordon1985some} (see also Corollary 3.13 from~\cite{ledoux2013probability}):
\begin{theorem}\label{th:gordon_classic}
	Let $X_{i,j}$ and $(Y_{i,j})$, $1 \leq i \leq n$, $1 \leq j \leq m$ be two (centered) Gaussian random vectors such that
	$$
	\begin{cases}
		\E X_{i,j}^2 = \E Y_{i,j}^2 & \text{for all} \ i,j \,, \\
		\E X_{i,j}X_{i,k} \geq \E Y_{i,j} Y_{i,k} & \text{for all} \ i,j,k \,, \\
		\E X_{i,j}X_{l,k} \leq \E Y_{i,j} Y_{l,k} & \text{for all} \ i\neq l \ \text{and} \ j,k  \,.
	\end{cases}
	$$
	Then, for all real numbers $\lambda_{i,j}$:
	$$
	\P \left( \bigcap_{i=1}^n\ \bigcup_{j=1}^m \big\{ X_{i,j} > \lambda_{i,j} \big\} \right)
	\leq
	\P \left( \bigcap_{i=1}^n\ \bigcup_{j=1}^m \big\{ Y_{i,j} > \lambda_{i,j} \big\} \right) \,.
	$$
\end{theorem}

\begin{theorem}\label{th:gordon0}
Let $D_u \subset \R^n$ and $D_v \subset \R^m$ be two compact sets.
Let $Q: D_u \times D_v \to \R$ be a continuous function.
Let $\big(X(u,v)\big)_{(u,v)\in D_u \times D_v}$ and $\big(Y(u,v)\big)_{(u,v)\in D_u \times D_v}$ be two centered Gaussian processes. Suppose that the functions
$$
(u,v) \mapsto X(u,v) \qquad \text{and} \qquad (u,v) \mapsto Y(u,v)
$$
are continuous on $D_u \times D_v$ almost surely. Assume that
$$
\begin{cases}
	\E \big[ X(u,v)^2\big] = \E \big[ Y(u,v)^2\big] & \text{for all} \quad (u,v) \in D_u \times D_v\,, \\
	\E \big[ X(u,v) X(u,v') \big] \geq \E \big[ Y(u,v) Y(u,v')\big] & \text{for all} \quad u \in D_u, v,v'\in D_v \,, \\
	\E \big[ X(u,v) X(u',v') \big] \leq \E \big[ Y(u,v) Y(u',v')\big] & \text{for all} \quad u,u' \in D_u, v,v'\in D_v \quad \text{such that} \quad u \neq u' \,.
\end{cases}
$$
		\vspace{-0.1cm}
		Then for all $t \in \R$
		\vspace{-0.1cm}
		$$
		\P\Big(\min_{u \in D_u}\max_{v \in D_v} Y(u,v) + Q(u,v) \leq t\Big)
		\leq \P\Big(\min_{u \in D_u}\max_{v \in D_v} X(u,v) + Q(u,v) \leq t\Big) \,.
		$$
		\vspace{-0.4cm}
\end{theorem}
\begin{proof}
Define the random variable
$$
d_0 = \sup \Big\{ d \in \Q_+ \, \Big| \, \forall (z,z') \in (D_u \times D_v)^2, \ \|z-z'\| \leq d \implies \big(|X(z)-X(z')| \leq \epsilon \ \text{and} \ |Y(z)-Y(z')| \leq \epsilon \big)  \Big\} \,.
$$
$X$ and $Y$ are continuous on the compact set $D_u \times D_v$ and are therefore uniformly continuous on this set: $d_0 >0$ almost surely. 
Let $\epsilon>0$. 
By tightness there exists a constant $d>0$ such that
$$
\P(d_0 \geq d) \geq 1-\epsilon \,.
$$
$Q$ is continuous and thus uniformly continuous on $D_u \times D_v$: there exists $\delta \in (0,d]$ such that for all $z,z' \in D_u \times D_v,$ $\|z-z'\| \leq \delta \implies |Q(z)-Q(z')| \leq \epsilon$.

Let $D_u^{\delta}$ (respectively $D_v^{\delta}$) be a $\delta/\sqrt{2}$-net of $D_u$ (respectively $D_v$). $D_u^{\delta} \times D_v^{\delta}$ is thus a $\delta$-net of $D_u \times D_v$. By Theorem~\ref{th:gordon_classic} we have for all $t \in \R$
		$$
		\P\Big(\min_{u \in D_u^{\delta}}\max_{v \in D_v^{\delta}} X(u,v) + Q(u,v) > t\Big)
		\leq \P\Big(\min_{u \in D_u^{\delta}}\max_{v \in D_v^{\delta}} Y(u,v) + Q(u,v) > t\Big) \,,
		$$
		which gives by taking the complementary:
		$$
		\P\Big(\min_{u \in D_u^{\delta}}\max_{v \in D_v^{\delta}} Y(u,v) + Q(u,v) \leq t\Big)
		\leq \P\Big(\min_{u \in D_u^{\delta}}\max_{v \in D_v^{\delta}} X(u,v) + Q(u,v) \leq t\Big) \,.
		$$
		By construction of $\delta$ we have with probability at least $1-\epsilon$
		$$
\Big|
\min_{u \in D_u^{\delta}}\max_{v \in D_v^{\delta}} X(u,v) + Q(u,v)
- \min_{u \in D_u}\max_{v \in D_v} X(u,v) + Q(u,v) \Big| \leq 2 \epsilon \,,
		$$
		and similarly for $Y$. We have therefore, for all $t\in\R$
		$$
		\P\Big(\min_{u \in D_u}\max_{v \in D_v} Y(u,v) + Q(u,v) \leq t -2 \epsilon \Big) - \epsilon
		\leq \P\Big(\min_{u \in D_u}\max_{v \in D_v} X(u,v) + Q(u,v) \leq t + 2\epsilon \Big) + \epsilon \,,
		$$
		and thus
		$$
		\P\Big(\min_{u \in D_u}\max_{v \in D_v} Y(u,v) + Q(u,v) \leq t  \Big) 
		\leq \P\Big(\min_{u \in D_u}\max_{v \in D_v} X(u,v) + Q(u,v) \leq t + 4\epsilon \Big) + 2\epsilon \,,
		$$
		which proves the theorem by taking $\epsilon \to 0$.

\end{proof}
\\

\begin{corollary}\label{cor:gordon0}
	Let $D_u \subset \R^{n_1+n_2}$ and $D_v \subset \R^{m_1 + m_2}$ be compact sets and let $Q: D_u \times D_v \to\R$ be a continuous function. Let $G = (G_{i,j}) \iid \cN(0,1)$, $g \sim \cN(0,\bbf{I}_{n_1})$ and $h \sim \cN(0,\bbf{I}_{m_1})$ be independent standard Gaussian vectors. 
	For $u \in \R^{n_1+n_2}$ and $v \in \R^{m_1+m_2}$ we define $\tilde{u} = (u_1,\dots,u_{n_1})$ and $\tilde{v} = (v_1,\dots,v_{m_1})$.
	Define 
	$$
	\begin{cases}
		C^*(G) = \min\limits_{u \in D_u} \max\limits_{v \in D_v} \tilde{v}^{\sT} G \tilde{u} + Q(u,v) \,, \\
		L^*(g,h) = \min\limits_{u \in D_u} \max\limits_{v \in D_v} \|\tilde{v}\| g^{\sT} \tilde{u} + \|\tilde{u}\| h^{\sT} \tilde{v} + Q(u,v) \,.
	\end{cases}
	$$
	Then we have:
	\begin{itemize}
		\item For all $t \in \R$,
			$$
			\P \Big(C^*(G) \leq t \Big) \leq 2 \P \Big(L^*(g,h) \leq t \Big) \,.
			$$
		\item If $D_u$ and $D_v$ are convex and if $Q$ is convex concave, then for all $t \in \R$
			$$
			\P \Big(C^*(G) \geq t \Big) \leq 2 \P \Big(L^*(g,h) \geq t \Big) \,.
			$$
	\end{itemize}
\end{corollary}
\begin{proof}
	Let us consider the Gaussian processes:
	$$
	\begin{cases}
		X(u,v) = \|\tilde{v}\| g^{\sT} \tilde{u} + \|\tilde{u}\| h^{\sT} \tilde{v} \,, \\
		Y(u,v) = \tilde{v}^{\sT} G \tilde{u} + \|\tilde{u}\| \|\tilde{v} \| z \,,
	\end{cases}
	$$
	where $z \sim \cN(0,1)$ is independent from $G$. Let $(u,v), (u',v') \in D_u \times D_v$ and compute
	\begin{align*}
		\E \big[ Y(u,v) Y(u',v') \big] &- \E \big[ X(u,v) X(u',v')\big]
		\\
	&= \|\tilde{u}\| \|\tilde{v}\| \|\tilde{u}'\| \|\tilde{v}'\|  + (\tilde{u}^{\sT} \tilde{u}')(\tilde{v}^{\sT} \tilde{v}')
	- \|\tilde{v}\| \|\tilde{v}'\|(\tilde{u}^{\sT} \tilde{u}') - \|\tilde{u}\| \|\tilde{u}'\|(\tilde{v}^{\sT} \tilde{v}')
		\\
		&= 
		\big( \|\tilde{u}\| \|\tilde{u}'\| - (\tilde{u}^{\sT} \tilde{u}') \big)
		\big( \|\tilde{v}\| \|\tilde{v}'\| - (\tilde{v}^{\sT} \tilde{v}') \big) \geq 0 \,.
	\end{align*}
	Therefore $X$ and $Y$ verify the covariance inequalities of Theorem~\ref{th:gordon0}: one can apply Theorem~\ref{th:gordon0}:
	$$
	\P\Big(\min_{u \in D_u} \max_{v \in D_v} Y(u,v) + Q(u,v) \leq t \Big) \leq \P \Big(\min_{u \in D_u} \max_{v \in D_v} Y(u,v) + Q(u,v) \leq t \Big) \,,
	$$
	We have then
	\begin{align*}
\P \Big(\min_{u \in D_u} \max_{v \in D_v} Y(u,v) + Q(u,v) \leq t \Big)
&\geq \frac{1}{2}
\P \Big(\min_{u \in D_u} \max_{v \in D_v} Y(u,v) + Q(u,v) \leq t \, \Big| \, z \leq 0 \Big)
\\
&\geq \frac{1}{2}
\P \Big(\min_{u \in D_u} \max_{v \in D_v} \tilde{v}^{\sT} G \tilde{u} + Q(u,v) \leq t \, \Big| \, z \leq 0 \Big)
= \frac{1}{2}
\P \Big(C^*(G) \leq t \Big) \,,
	\end{align*}
	which proves that
	$$
	\P\Big(\min_{u \in D_u} \max_{v \in D_v} \tilde{v}^{\sT} G \tilde{u} + Q(u,v)  \leq t \Big)
	\leq 2 \P \Big(\min_{u \in D_u} \max_{v \in D_v} \|\tilde{v}\| g^{\sT} \tilde{u} + \|\tilde{u}\| h^{\sT} \tilde{v} + Q(u,v) \leq t \Big) \,.
	$$
	Let us suppose now that $D_u$ and $D_v$ are convex and that $G$ is convex-concave. We now apply the inequality we just proved, but with the role of $u$ and $v$ being switched (and $-Q$ and $-t$ instead of $Q$ and $t$):
	$$
	\P\Big( \min_{v \in D_v} \max_{u \in D_u} \tilde{v}^{\sT} G \tilde{u} - Q(u,v)  \leq -t \Big)
	\leq 2 \P \Big( \min_{v \in D_v} \max_{u \in D_u} \|\tilde{v}\| g^{\sT} \tilde{u} + \|\tilde{u}\| h^{\sT} \tilde{v} - Q(u,v) \leq -t \Big) \,,
	$$
	which gives (using the fact that $(G,g,h)$ and $(-G,-g,-h)$ have the same law):
	$$
	\P\Big( \max_{v \in D_v} \min_{u \in D_u} \tilde{v}^{\sT} G \tilde{u} + Q(u,v)  \geq t \Big)
	\leq 2 \P \Big( \max_{v \in D_v} \min_{u \in D_u} \|\tilde{v}\| g^{\sT} \tilde{u} + \|\tilde{u}\| h^{\sT} \tilde{v} + Q(u,v) \geq t \Big) \,.
	$$
	By Proposition~\ref{prop:rockafellar}, one can switch the min-max of the left-hand side, because $Q$ is convex-concave and we are working on convex sets $D_u$ and $D_v$. For the right-hand side, we simply use the fact that:
	$$
\max_{v \in D_v} \min_{u \in D_u} \|\tilde{v}\| g^{\sT} \tilde{u} + \|\tilde{u}\| h^{\sT} \tilde{v} + Q(u,v) 
\leq
\min_{u \in D_u} \max_{v \in D_v}  \|\tilde{v}\| g^{\sT} \tilde{u} + \|\tilde{u}\| h^{\sT} \tilde{v} + Q(u,v) \,,
	$$
	to conclude the proof.
\end{proof}

\subsection{Basic concentration results}

We recall in this section some elementary concentration results, see Chapter~2 from~\cite{boucheron2013concentration} for a more detailed presentation of these facts.

\begin{definition}
A real random variable $X$ is said to be 
\begin{itemize}
	\item $\sigma^2$-sub-Gaussian if for every $s \in \R$, \ 
		$\displaystyle
		\log \E e^{s(X - \E[X])} \leq \frac{s^2 \sigma^2}{2} \,,
		$
	\item $(v,c)$-sub-Gamma if for every $s \in (-1/c, 1/c)$, \ 
		$\displaystyle
		\log \E e^{s(X - \E[X])} \leq \frac{s^2 v}{2(1-c|s|)} \,.
		$
\end{itemize}
\end{definition}

One deduces immediately from the above definition:
\begin{proposition}
Let $(X_1, \dots, X_n)$ be independent real random variables. Define $S = \sum_{i=1}^n X_i$.
\begin{itemize}
	\item Suppose that for all $i \in \{1,\dots,n\}$, $X_i$ is $\sigma_i^2$-sub-Gaussian. Then $S$ is $\sum_{i=1}^n \sigma_i^2$-sub-Gaussian.
	\item Suppose that for all $i \in \{1,\dots,n\}$, $X_i$ is $(v_i,c_i)$-sub-Gamma. Then $S$ is $\big(\sum_{i=1}^n v_i, \max c_i \big)$-sub-Gamma.
\end{itemize}
\end{proposition}

\begin{proposition}
Let $X$ be a real random variable.
\begin{itemize}
	\item if $X$ is $\sigma^2$-sub-Gaussian, then for all $t >0$
			$$
			\P( X - \E[X] \geq t ) \ \vee \ \P(X - \E[X] \leq -t) \leq e^{-\frac{t^2}{2\sigma^2}} \,,
			$$
		\item if $X$ is $(v,c)$-sub-Gamma, then for all $t >0$
			$$
			\P( X - \E[X] \geq \sqrt{2ct} + vt ) \ \vee \ \P(X - \E[X] \leq -(\sqrt{2ct} + vt)) \leq e^{-t} \,.
			$$
	\end{itemize}
\end{proposition}

\begin{remark}
	The bound $\P(X > \sqrt{2vt} + ct) \leq e^{-t}$ implies that
	$$
	\P(X>t) \leq 
	\begin{cases}
		\exp\left( -\frac{t^2}{8v}\right) & \text{for} \quad 0 < t \leq \frac{2v}{c^2} \,, \\
		\exp\left( -\frac{t}{2c}\right) & \text{for} \quad  t \geq \frac{2v}{c^2}  \,.
	\end{cases}
	$$
\end{remark}

\begin{proposition}\label{prop:square_sub_gaussian}
	If $X$ is $\sigma^2$-sub-Gaussian and has mean $\mu$, then $X^2$ is a sub-Gamma random variable with parameters 
	$$
	\begin{cases}
		v=16 \sigma^2 + 4 \mu^2 \sigma^2 \,, \\
		c = 4 \sigma^2 \,.
	\end{cases}
	$$
\end{proposition}

\begin{proof}
	Let $\mu = \E[X]$ and $Y = X -\mu$. $X^2 = Y^2 + 2 \mu Y + \mu^2$. 
	$$
	\E \left[(Y^2)^2\right] = \E[(X-\mu)^4] \leq 16 \sigma^4  \,,
	$$
	\begin{align*}
		\E\left[(Y^2)^q\right]=
		\E\left[(X-\mu)^{2q}\right] \leq 2 q! (2\sigma^2)^q
		=\frac{1}{2} q! (16\sigma^2) (2\sigma^2)^{q-2} \,.
	\end{align*}
	By Bernstein's inequality (see for instance Theorem~2.10 in~\cite{boucheron2013concentration})
	$$
	\log \E e^{s(Y^2 - \E[Y^2])} \leq \frac{16\sigma^2 s^2}{2(1-2\sigma^2 |s|)} \,.
	$$
	$2\mu Y$ is $4 \mu^2 \sigma^2$-sub-Gaussian, therefore $\log \E e^{2\mu s Y} \leq 2 \mu^2 \sigma^2 s^2$ and
	\begin{align*}
		\log \E e^{s (X^2 - \E[X^2])} &=
		\log \E e^{s (Y^2 - \E[Y^2]) + 2\mu s Y} 
		\leq \frac{1}{2} \log \E e^{2 s (Y^2 - \E[Y^2])} 
		+ \frac{1}{2} \log \E e^{4\mu s Y} 
		\\
		&\leq \frac{16 \sigma^2 s^2}{1-4\sigma^2 s} + 4 \mu^2 \sigma^2 s^2
		\leq \frac{(16 \sigma^2 + 4 \mu^2 \sigma^2) s^2}{1-4\sigma^2 |s|}  \,.
	\end{align*}
	$X^2$ is therefore a Sub-Gamma random variable with variance factor $v = 16\sigma^2 + 4 \mu^2 \sigma^2$ and scale parameter $c=4\sigma^2$.
\end{proof}

\begin{lemma}\label{lem:conc_prod_bounded}
	Let $X$ be a $\sigma^2$-sub-Gaussian random variable. Define $m = \E[|X|]$. Let $Y$ be a random variable bounded by $1$. Then 
	$XY$ is $16(m^2 + 2 \sigma^2)$-sub-Gaussian.
\end{lemma}
\begin{proof}
	We have $| \E[XY] | \leq m$, therefore
	\begin{align*}
		\E \left[(XY - \E[XY])^{2q}\right] 
		&\leq 2^{2q-1} \E[X^{2q}] + 2^{2q-1} m^{2q}
		\leq  q! (8\sigma^2)^q + q! (4m^2)^{q}
		\leq  q! (8\sigma^2 + 4 m^2)^q  \,.
	\end{align*}
		\vspace*{-0.3cm}
\end{proof}

\subsection{Largest singular value of a Gaussian matrix}

The largest singular value of a $n \times N$ matrix $A$ is defined as
$$
\sigma_{\rm max}(A) = \max_{\|x\| \leq 1} \|A x\| \,.
$$

The next classical result is a simple consequence of Slepian's Lemma (see for instance~\cite{ledoux2013probability}, Section~3.3) and the classical Gaussian concentration inequality (see for instance~\cite{boucheron2013concentration}, Theorem~5.6).
\begin{proposition}\label{prop:max_singular}
	Let $G$ be a $n \times N$ random matrix, whose entries are i.i.d.\ $\cN(0,1)$.
	For all $t\geq0$ we have
	$$
	\P(\sigma_{\rm max}(G) > \sqrt{N} + \sqrt{n} + t) \leq e^{-t^2/2} \,.
	$$
	\end{proposition}

\end{appendices}

\bibliographystyle{plain}
\bibliography{references}

\end{document}